\documentclass[a4paper,11pt]{article}
\usepackage{amsmath}
\usepackage{latexsym}
\usepackage{amsfonts}
\usepackage{a4wide}
\usepackage{xspace}
\usepackage{color}
\usepackage{accents}

\usepackage[all]{xy}

\newcommand{\hem}{\hspace*{1em}}

\newcommand{\hmm}{\hspace*{2em}}

\newcommand{\hfl}{\hspace*{\fill}}

\newlength{\dede}     
\newcommand{\tab}[1]{
\settowidth{\dede}{#1}
\hspace*{\dede}}      

\newcommand{\nhp}{\hspace*{-\parindent}}
\newcommand{\noi}{\hspace*{-\parindent}}

\newcommand{\re}{\mbox{$\mathbb R$}}  
\newcommand{\N}{\mbox{$\mathbb N$}}   
\newcommand{\Zh}{\mbox{$\mathbb Z$}} 

\newcommand{\fr}[2]{\mbox{$\frac{\raisebox{2pt}{$#1$}}
{\raisebox{-3pt}{$#2$}}$}}

\newcommand{\lra}{\mbox{$\longrightarrow$}}  

\newcommand{\se}[2]{\mbox{$ #1_{\raisebox{-3pt}{$\scriptstyle #2$}}$}}

\newcommand{\D}[1]{\mbox{$\se D {\!{#1}}$}} 
\newcommand{\DG}{\D{G}} 
\newcommand{\Dthr}{\mbox{$\se D {\bf 3}$}} 
\newcommand{\Dt}[1]{\mbox{$\se {D^t} {\!{#1}}$}} 
\newcommand{\Dtthr}{\Dt {\bf 3}} 

\newcommand{\Spec}{\mbox{$\mathrm {Spec\,}$}} 
\newcommand{\Sper}{\mbox{$\mathrm {Sper\,}$}} 
\newcommand{\supp}{\mbox{$\mathrm {supp\,}$}} 
\newcommand{\sub}{\mbox{$\subseteq$}}

\newcommand{\sth}{\, \vert \,}
\newcommand{\0}{\mbox{$\emptyset$}} 

\newcommand{\Ga}{\mbox{$\Gamma$}}

\newcommand{\si}{\mbox{$\sigma$}}

\newcommand{\De}{\mbox{$\Delta$}}

\newcommand{\al}{\mbox{$\alpha$}}
\newcommand{\bt}{\mbox{$\beta$}}

\newcommand{\py}{\mbox{$\pi$}}

\newcommand{\cC}{\mbox{$\cal C$}}

\newcommand{\cH}{\mbox{$\cal H$}}

\newcommand{\fM}{\mbox{$\mathfrak M$}}

\newcommand{\cl}[1]{\mbox{$\overline {#1}$}}

\newcommand{\con}{\mathrm{con}}
\newcommand{\spez}{\,\mbox{$\rightsquigarrow$\,}}
\newcommand{\ex}{\mbox{$\exists$}}
\newcommand{\fa}{\mbox{$\forall$}}
\newcommand{\w}{\mbox{$\wedge$}} 
\newcommand{\jo}{\mbox{$\vee$}} 

\newcommand{\HTS}[1]{\mbox{${\mathrm {Hom}}${\raisebox{-5pt}{$\scriptstyle{\mathrm {TS}}$}
{\!$(#1$,\,{\bf 3}$)$}}}}
\newcommand{\HRS}[1]{\mbox{${\mathrm {Hom}}${\raisebox{-5pt}{$\scriptstyle{\mathrm {RS}}$}
{\!$(#1$,\,{\bf 3}$)$}}}}

\newcommand{\Ra}{\mbox{$\Rightarrow$}}
\newcommand{\Lra}{\mbox{$\Leftrightarrow$}}

\newcommand{\La}{\mbox{$\Leftarrow$}}

\newcommand{\und}[1]{\raisebox{-.2ex}{\underline{\raisebox{.2ex}{#1}}}}  
\newcommand{\y}[1]{\mbox{$#1$}}

\newcommand{\bc}{\begin{center}}
\newcommand{\ec}{\end{center}}

\newcommand{\lbr}{\linebreak}

\newcommand{\fm}[1]{\mbox{$\langle \, #1 \,  \rangle$}}

\newcommand{\h}[1]{\mbox{$\widehat{#1}$}}

\newtheorem{Th}{Theorem}[section]
\newtheorem{Co}[Th]{Corollary}
\newtheorem{Df}[Th]{Definition}
\newtheorem{Pro}[Th]{Proposition}
\newtheorem{Le}[Th]{Lemma}
\newtheorem{Exa}[Th]{Example}
\newtheorem{Exas}[Th]{Examples}
\newtheorem{Rem}[Th]{Remark}
\newtheorem{Dfre}[Th]{Definition and Remarks}
\newtheorem{Dfno}[Th]{Definition and Notation}
\newtheorem{Rems}[Th]{Remarks}
\newtheorem{Remno}[Th]{Remarks and Notation}
\newtheorem{Fa}[Th]{Fact}
\newtheorem{No}[Th]{Notation}
\newtheorem{Prel}[Th]{Preliminaries}
\newtheorem{Prelnot}[Th]{Preliminaries and Notation}
\newtheorem{Ct}[Th]{\hspace*{-4pt}}

\newcommand{\bre}{\begin{Rem} \em}
\newcommand{\bdfr}{\begin{Dfre}}
\newcommand{\edfr}{\end{Dfre}}
\newcommand{\bdfn}{\begin{Dfno} \em}
\newcommand{\edfn}{\end{Dfno}}
\newcommand{\bdf}{\begin{Df} \em}
\newcommand{\edf}{\end{Df}}
\newcommand{\bremn}{\begin{Remno} \em}
\newcommand{\eremn}{\end{Remno}}
\newcommand{\bth}{\begin{Th}}
\newcommand{\eth}{\end{Th}}
\newcommand{\bco}{\begin{Co}}
\newcommand{\eco}{\end{Co}}
\newcommand{\ble}{\begin{Le}}
\newcommand{\ele}{\end{Le}}
\newcommand{\bpr}{\begin{Pro}}
\newcommand{\epr}{\end{Pro}}
\newcommand{\bex}{\begin{Exa} \em}
\newcommand{\eex}{\end{Exa}}
\newcommand{\bexs}{\begin{Exas} \em}
\newcommand{\eexs}{\end{Exas}}
\newcommand{\ere}{\end{Rem}}
\newcommand{\bres}{\begin{Rems} \em}
\newcommand{\eres}{\end{Rems}}
\newcommand{\bfa}{\begin{Fa}}
\newcommand{\efa}{\end{Fa}}
\newcommand{\bno}{\begin{No} \em}
\newcommand{\eno}{\end{No}}
\newcommand{\bprel}{\begin{Prel} \em}
\newcommand{\eprel}{\end{Prel}}
\newcommand{\bprelnot}{\begin{Prelnot} \em}
\newcommand{\eprelnot}{\end{Prelnot}}
\newcommand{\bct}{\begin{Ct} \em}
\newcommand{\ect}{\end{Ct}}

\newcommand{\llab}[1]{\ar@{}[l]|-<<{\txt{\footnotesize #1}}} 
\newcommand{\sllab}[1]{\ar@{}[l]|-<<{\txt{\scriptsize #1}}}  
\newcommand{\ulab}[1]{\ar@{}[u]|-<<{\txt{\footnotesize #1}}} 
\newcommand{\dlab}[1]{\ar@{}[d]|-<<{\txt{\footnotesize #1}}} 
\newcommand{\rlab}[1]{\ar@{}[r]|-<<<{\txt{\footnotesize #1}}} 
\newcommand{\srlab}[1]{\ar@{}[r]|-<<{\txt{\scriptsize #1}}}  
\newcommand{\lulab}[1]{\ar@{}[l]_<<{\txt{\footnotesize #1}}} 
\newcommand{\rulab}[1]{\ar@{}[r]^<<{\txt{\footnotesize #1}}} 
\newcommand{\ldlab}[1]{\ar@{}[l]^<<<{\txt{\footnotesize #1}}} 
\newcommand{\rdlab}[1]{\ar@{}[r]_>>>>{\txt{\footnotesize #1}}} 
\newcommand{\edge}[1]{\ar@{-}[#1]} 
\newcommand{\node}{*{\bullet}} 

\def\NewFont#1#2#3#4#5{%
\expandafter\font\csname #1display\endcsname =#1 at #2%
\expandafter\font\csname #1normal\endcsname =#1 at #3%
\expandafter\font\csname #1script\endcsname =#1 at #4%
\expandafter\font\csname #1scriptscript\endcsname =#1 at #5%
}
\def\NewFontLetter#1#2{{%
\mathchoice%
{{\expandafter\hbox{\csname #1display\endcsname\char"#2}}}%
{{\expandafter\hbox{\csname #1normal\endcsname\char"#2}}}%
{{\expandafter\hbox{\csname #1script\endcsname\char"#2}}}%
{{\expandafter\hbox{\csname #1scriptscript\endcsname\char"#2}}}%
}}%

\ifx\undefined\fontshape   
\else
\fi

\begin{document}

\title{Fans in the Theory of Real Semigroups \\
I. Algebraic Theory}

\author{M. Dickmann \and A. Petrovich}

\date{March 2017}

\maketitle

\pagestyle{myheadings}

\vspace*{-1cm}

\begin{abstract} In \cite{DP1} we introduced the notion of a {\it real semigroup} (RS) 
as an axiomatic framework to study diagonal quadratic forms with arbitrary entries over
(commutative, unitary) semi-real rings. (For the axioms of RS, cf. \ref{RS}.)
Two important classes of RSs were studied at length in \cite{DP2}, \cite{DP3}. In this paper 
we introduce and develop the algebraic theory of {\it RS-fans}, a third class of RSs providing a vast 
generalization of homonymous notions previously existing in field theory and in the theories
of abstract order spaces and of reduced special groups; for a background on fans, see 
paragraph A of the Introduction, below. The contents of this paper are briefly reviewed in
paragraph B of the Introduction. The combinatorial theory of the structures dual to RS-fans, 
called {\it ARS-fans}, is the subject of \cite{DP5b}, a continuation of the present paper.
\end{abstract}

\vspace*{-0.5cm}

\section*{Introduction} \label{fan-intro} 

\vspace*{-0.4cm}

In \cite{DP1} we introduced the notion of a {\it real semigroup} (henceforth abbreviated
RS), an axiomatic framework aimed at studying diagonal quadratic forms with arbitrary
entries over commutative, unitary rings\,\footnote{\,In this paper simply referred to  as 
{\it rings}.} admitting a minimum of orderability.\,\footnote{\,Namely, having a non-empty
real spectrum or, equivalently, that $-1$ is not a sum of squares.} For ready reference we
have included below the axioms defining RSs (\ref{RS}) and their underlying structures, 
the {\it ternary semigroups} (abbreviated TS, see \ref{ts}). The basic properties of these
structures are proved in \cite{DP1}, \S\S\,1,2, pp. 100-112, and \cite{DP2}, \S\,2, pp. 57-59. 
We also proved (\cite{DP1}, Thm. 4.1, p. 115) that the RSs are categorically dual to the 
{\it abstract real spectra} (ARS), previously introduced in \cite{M}, Chs. 6 -- 9, with a similar goal.
\vspace{-0.1cm}

In \cite{DP2}, \cite{DP3} we introduced and studied two outstanding classes of RSs, the Post 
algebras and the spectral real semigroups, and their dual ARSs. The aim of this paper and its
continuation [DP5b] is to present a third natural class of RSs (and their dual ARSs), 
namely {\it fans}, and develop their theory.
\vspace{-0.1cm}

{\bf A. Background on fans.} Fans were introduced by Becker and K\"opping \cite{BK}\,\footnote{\,See
\cite{ABR}, p. 84, and \cite{La}, Notes on \S\,5, p. 48.} as a distinguished class of 
preorders in fields, and further investigated by several authors. Chapter 5 of the 
monograph \cite{La} gives a quite complete picture of the role of fans in the context 
of fields, and contains many bibliographical references.
\vspace{-0.1cm}

A further step was taken by Marshall, see \cite{M}, Ch. 3, who generalized the notion 
of a fan to the context of {\it abstract spaces of orderings} (AOS), an axiomatic
framework extending the field case. In \cite{Li} (see also \cite{DM1}, Ex. 1.7, pp. 8-9,
and pp. 89-90) this notion was treated in the framework of {\it reduced special groups} 
(RSG), and its functorial duality with the corresponding notion of fan in the category 
of AOSs proved.
\vspace{-0.1cm}

Fans surfaced again in \cite{ABR}, Chs. 3, 5, under the still more general clothing of
{\it spaces of signs}, a framework equivalent to that of ARSs. This book extensively witnesses 
the key role that fans play in real algebraic and real analytic geometry; see, e.g., \cite{ABR}, 
Thms. IV.7.3 and V.1.4 (the ``generation formulae''), and \cite{AR}, pp. 1-7, where further 
references can be found. However, the notion of a fan used in \cite{ABR} (cf. Def. 3.12, 
p. 75) is that of an AOS-fan suitably {\it embedded} in an ARS; see also \cite{M}, p. 162.
\vspace{-0.1cm}

A sufficiently general and {\it intrinsic} notion of fan in the categories of ARSs 
and of RSs\,\footnote{\,The categories of ARSs, RSs, AOSs, RSG's, under natural morphisms, 
will be denoted by boldfacing the corresponding acronyms.} does not exist at present. The 
aim of this paper and its continuation [DP5b], is to present and study such a 
notion. The key leading to this goal consists in bringing into play the (enriched) semigroup 
structures underlying the real semigroups, namely the {\it ternary semigroups} (\cite{DP1}, 
Def. 1.1, p. 100; see also \ref{ts} below). This task \und{is not} a straightforward extrapolation
of the situation in the categories {\bf {RSG}} and {\bf {AOS}}, on two accounts. Firstly, owing 
to the rather complex algebraic and topological structure of the TSs, far more involved than the 
(trivial) structures underlying the RSGs (namely groups of exponent 2 with a distinguished 
element $-1$). Secondly, because the natural topology on ARSs is spectral ---i.e., has non-trivial 
specialization--- while the corresponding topology on AOSs is Boolean, and therefore has a trivial 
specialization order. While the first of these factors plays a central role in the present paper, the 
second will be crucial in its sequel, [DP5b], where we develop the combinatorial theory of ARS-fans.
\vspace{-0.1cm}

To motivate the main ideas presented in this paper, we begin by briefly reviewing the 
definition of a fan in the (dual) categories {\bf AOS} and {\bf RSG} (for more 
details, see \cite{M}, Ch. 3; \cite{Li}, Ex. 1.1.6, pp. 30-31; \cite{DM1}, Ex. 1.7, pp. 8-9).
\vspace{-0.15cm}

\noindent--- A  fan in the category {\bf AOS} (henceforth called an 
{\bf AOS-fan}) is an abstract space of orders $(X,G)$ where ``$X$ is biggest possible"; 
there are two equivalent ways of making sense of this idea\,: 
\vspace{-0.15cm}

\noindent (1) $X$ consists of all group homomorphisms $h:G \, \lra \, \{\, \pm 1 \}$ 
such that $h(-1) = -1$.
\vspace{-0.15cm}

\noindent (2) $(X,G)$ is an AOS and $X$ is closed under the product of any 
three of its members.
\vspace{-0.15cm}

\noindent--- A  fan in the category {\bf RSG} (henceforth an 
{\bf RSG-fan}) is a reduced special group $G$ whose binary representation relation is 
``smallest posible"; there is only one way of making sense of this:
\vspace{-0.15cm}

\noindent [RSG-fan]\hfl $a \in D_{G}(b,c)$  \  \  iff  \  \ either \ $b = -c$ \ or 
\ ($b \neq -c$ \ and \ $a \in \{b,c \}$). \hfl
\vspace{-0.4cm}

\bres \label{remksAOSFans} (a) While condition (1) above implies that  $(X,G)$ is 
an AOS, the last requirement in (2) alone is not sufficient to guarantee that 
$(X,G)$ is an AOS; in addition, one must require that:
\vspace{-0.15cm}

\nhp (i)\, \y X separates points in $G$, i.e., $\bigcap_{\,\sigma \in X}$ 
\mbox{\textrm {ker}}($\sigma) = \{\, 1 \}$.
\vspace{-0.15cm}

\nhp (ii) \y X verifies the following maximality condition (see \cite{M}, 
axiom [AX2] for AOSs, p. 22): for every group homomorphism 
$\sigma: G\, \lra \, \{\, \pm 1 \}$, if\; $\si(-1) = -1$\; and\; 
$a, b \in {\mathrm {ker}}(\sigma) \; \Ra $
$\se D X(a,b) \, \sub \, {\mathrm {ker}}(\si)$,\; then $\sigma \in X$.
\vspace{-0.1cm}

\nhp (b) The definition of binary representation given by condition [RSG-fan] above 
(together with $1 \neq -1$) implies that $G$ is a RSG (\cite{Li}, Prop. 1.1.14, pp.  
34-36). Thus, every non-trivial group of exponent 2 is endowed with a structure of
RSG. This \und{is not} the case for ternary semigroups, where an additional condition
(called condition [Z]) ought to be satisfied for a TS to be endowed with a structure of RS, 
cf. Fact \ref{zeros}, Proposition \ref{charfan} and Theorem \ref{fanRS}.  \hfl $\Box$
\eres
\vspace{-0.45cm}

We define the notion of a fan in the category {\bf ARS} of abstract 
real spectra by postulating the analogs of conditions (1) and (2) above, 
{\it upon replacing} the underlying notion of a group of exponent 2 with a 
distinguished element $-1$ {\it by that of a ternary semigroup} and, of course, the 
target group $\{\, \pm 1 \}$ by the ternary semigroup ${\bf 3} = \{ -1,0,1 \}$:

\vspace{-0.45cm}

\bdf \label{fans-various} Given a ternary semigroup $G$ and a non-empty set
$X \subseteq \HTS G$\,\footnote{\;\HTS G denotes the set of all TS-homomorphisms 
from the TS $G$ into the TS ${\bf 3}$.},
\vspace{-0.15cm}

\noindent $(1)$ $(X,G)$ is a\, {\bf $\se {\textrm{fan}} {1}$} iff $X$ consists of 
\underline{all} TS-homomorphisms from $G$ to ${\bf 3} = \{ -1,0,1 \}$, i.e., 
$X = \HTS G$.
\vspace{-0.15cm}

\noindent $(2)$ $(X,G)$ is a {\bf $\se {\textrm{fan}} {2}$} iff it is an ARS and 
$X$ is closed under the product of \underline{any three} of its members.
\vspace{-0.15cm}

\nhp We shall frequently use in the sequel the following weaker notion to 
which we give a name:
\vspace{-0.15cm}

\nhp $(3)$ $(X,G)$ is a\, {\bf q-fan} $($quasi-fan$)$ iff \y X is closed 
under the product of any three of its members and $X$ separates points in 
$G$, i.e., for every $a,b \in G, \ a \neq b$, there is \ $h \in X$ such that \ 
$h(a) \neq h(b)$. \hfl $\Box$
\edf
\vspace{-0.7cm}

\bres \label{ex:q-fan} (a) The set ${\bf 3} = \{ -1,0,1 \}$ under obvious operations
has a unique structure of TS. In fact, endowed with suitable ternary representation
and transversal representation relations, see \ref{Ex3}, it has a {\it unique} structure 
of RS. It obviously is a $\se {\textrm{fan}} {1}$. 
\vspace{-0.15cm}

\noi (b) In \ref{fans-various}\,(2) we allow products of type 
$h_1^{2}h_2$; as opposed to the case of special groups, squaring a TS-homomorphism
does not produce a map constantly equal to $1$. Note also that $h^{3} = h$, and that 
the product of any three TS-homomorphisms is again a TS-homomorphism.
\vspace{-0.15cm}

\nhp (c) An obvious example of q-fan over a TS, $G$, is $(\HTS G, G)$: 
$\HTS G$ is closed under product of any three members, and separates points 
in $G$ by the separation theorem for TSs, \cite{DP1}, Thm. 1.9, pp. 103-104. 
\hfl $\Box$
\eres
\vspace{-0.35cm}

\noi {\bf B. Contents of the paper.}
\vspace{-0.15cm}

\noi In Section 1 we include for ready reference the axioms defining ternary
semigroups (\ref{ts}) and real semigroups (\ref{RS}), as well as the basic example
\ref{Ex3}. Other than these, the section reviews briefly some algebraic and topological 
notions and results used throughout the paper which do not appear in print elsewhere. 
\vspace{-0.1cm}

In \S\,\ref{TStoRS} we prove (Theorem \ref{RSaxioms}) that the ternary representation relation 
induced on any TS by a non-empty set of TS-characters separating points, automatically satisfies all
axioms for real semigroups {\it with the exception of} the strong associativity axiom [RS3], see
Definition \ref{RS}.
\vspace{-0.1cm}

In \S\,\ref{RSfan} we give a purely algebraic characterization of the representation and the transversal
representation relations naturally occurring in a q-fan satisfying condition [Z] (Theorem 
\ref{fanrepr}). We also show that any ternary semigroup endowed with a ternary relation 
satisfying these algebraic requirements is automatically a real semigroup (Theorem \ref{fanRS}).
A number of important consequences follow from these results; notably
\vspace{-0.15cm}

\nhp $\bullet$ The identity of the two notions of fan defined in \ref{fans-various} above
(Proposition \ref{equivfan}).
\vspace{-0.15cm}

\nhp $\bullet$ Various properties substantiating the fact that both representation relations 
in RS-fans are 
\vspace{-0.2cm}

\tab{\nhp $\bullet$ }``smallest possible'' (Corollaries \ref{fanmaps}\,--\,\ref{fanideals}).
\vspace{-0.15cm}

Section \ref{fan-examp} gives a few examples of finite RS- and ARS-fans constructed from 
ternary semigroups on one and three generators. For each of these examples we draw the 
graph of the representation partial order of the corresponding RS-fan (cf. Definition \ref{reprord}), 
and that of the specialization root system of its dual ARS-fan. We prove that, under the representation 
partial order, each RS-fan is a bounded lattice (not modular, in general), Theorem \ref{fanlattice}.
\vspace{-0.15cm}

After a brief presentation of a general notion of congruence in real semigroups, in \S\,\ref{fan-quot} 
we prove that the class of RS-fans is closed under the formation of arbitrary quotients (Proposition
\ref{pr:fanQuot-1}), and that quotients of RS-fans by congruences defined by ideals are (upon
omitting zero) fans in the category of reduced special groups (Proposition \ref{fanQuot-2}).
\vspace{-0.15cm}

Section \ref{sect:char-fans} contains a characterization of RS-fans in terms of specialization and quotients
(Theorem \ref{thm:char-fan}). It follows that certain geometric configurations of the character 
space of RSs give rise to RS-fans (Corollaries \ref{TotOrder=Fan} and \ref{2chains=Fan}).
\vspace{-0.15cm}

In \S\,\ref{fan-prings} we apply the theory previously developed to the real semigroups 
associated to preordered rings. Proposition \ref{TS-charInP-rings} gives a characterization of the 
{\it ternary semigroup} characters of the RS $\se G {\!A,T}$ associated to a preordered ring 
$\fm{A,T}$ in terms of the algebraic operations of $A$ and the preorder $T$. This yields a 
natural extension of the original definition of fans by \cite{BK} (cf. \cite{La}, Def. 5.1, p. 39) to 
the context of rings and, hence, a characterization of those RSs $\se G {\!A,T}$ that are fans in 
terms of $\fm{A,T}$ (Proposition \ref{T-convPrimesTotOrd}\,(2)).
\vspace{-0.15cm}

A {\it total preorder} $T$ of a ring $A$ is a (proper) preorder such that $T \cup -T = A$.
Theorem \ref{TotPreordsAndFans} proves that if $T$ is either a total preorder of $A$ or the
intersection of two total preorders such that the set of $T$-convex prime ideals of $A$ is
totally ordered under inclusion, then $\se G {\!A,T}$ is a RS-fan. This gives a ring-theoretic
analog of the notion of a {\it trivial fan}, well-known in the field case (\cite{La}, Prop. 5.3,
p. 39).

\noi {\bf Acknowledgements.} First author partially supported by a ``Projet de mobilit\'e 
international Sorbonne Paris Cit\'e\,$-$\,Argentine\,$-$\,Br\'esil 2014''. The authors wish 
to thank F. Miraglia for his careful reading of and fruitful comments on the final draft of 
this paper.

\vspace{-0.35cm}

\section{Preliminaries\vspace{-0.3cm}} \label{fan-prel}

\paragraph{A. Ternary\vspace{-0.3cm} semigroups.} \label{TS}

\bdf \label{ts} A {\bf ternary semigroup} $($abbreviated {TS}$)$ is a 
structure $\langle S, \cdot\,,1,0,-1 \rangle$ with individual constants 
$1,0,-1,$ and a binary operation $``\cdot"$ such that:
\vspace{-0.1cm}

\nhp [TS1]\; $\langle S, \cdot\, ,1 \rangle$ is a commutative semigroup with unit.
\vspace{-0.1cm}

\nhp [TS2]\:\, $x^{3} = x$\, for all $x \in S$.
\vspace{-0.1cm}

\nhp [TS3]\,\, $-1 \neq 1$ and $(-1)(-1) = 1$.
\vspace{-0.1cm}

\nhp [TS4]\;\, $x \cdot 0 = 0$\, for all $x \in S$.
\vspace{-0.1cm}

\nhp [TS5]\:\, For all $x \in S$, $x = -1 \cdot x\; \Rightarrow\; x =\, 0$.
\vspace{-0.1cm}

\nhp We shall write $-x$ for $-1 \cdot x$. \hfl $\Box$
\edf
\vspace{-0.3cm}

In \S\,1 of \cite{DP1}, pp. 100-105, the reader will find an account of basic results
on ternary semigroups. In particular, the separation theorem \cite{DP1}, Thm. 1.9,
pp. 103-104, and notation and results on the spectral and constructible topologies
on the set $\HTS T$ of characters of a TS, $T$, with values in ${\bf 3} = \{1, -1, 0\}$,
cf. \cite{DP1}, pp. 104-105, are repeatedly used in this paper. 

\vspace{-0.1cm}

\nhp {\bf Warning.} Throughout this paper {\it the default topology} on all character spaces 
{\it is the spectral topology}. Whenever the associated constructible topology is used, the 
modifier $(.)_{\con}$ will be attached to the name of the space. \hfl $\Box$

Beyond the results in \cite{DP1}, \S\,1, we shall need the following results, which do not appear
therein. The next Lemma gives several characterizations of the specialization order of the spectral
topology in ternary semigroups.
\vspace{-0.4cm}

\ble \label{char-specializ} Let $T$ be a TS, and let $g,h \in \se X T$. The following are equivalent:
\vspace{-0.15cm}

\nhp $(1)$ $g \spez \, h$ $($i.e., \y h is an specialization of \y g$)$.
\vspace{-0.15cm}

\nhp $(2)$ $h^{-1}[1]\, \sub \, g^{-1}[1]$\, $($equivalently,\, $h^{-1}[-1]\, \sub \, g^{-1}[-1])$.
\vspace{-0.15cm}

\nhp $(3)$ $g^{-1}[\{0,1\}]\, \sub \, h^{-1}[\{0,1\}]$.
\vspace{-0.15cm}

\nhp $(4)$ $Z(g)\, \sub \, Z(h)$ and $\fa\, a \in G\, (a \not\in Z(h) \;\; \Ra \;\; g(a) = h(a))$.
\vspace{-0.15cm}

\nhp $(5)$ $h = h^2g$ $($equivalently, $h^2 = hg)$.
\ele 

\vspace{-0.4cm}
\nhp {\bf Proof.} The specialization partial order in any spectral space is defined by:
\vspace{-0.15cm}

\nhp \hfl $g \spez \, h$ \, iff \, $h \in \cl {\{g\}}$\; iff \, for every 
subbasic open \y U,\; $h \in U \; \Ra \; g \in U$, \hfl
\vspace{-0.15cm}

\nhp  Since the subbasic opens of $\se X T$ are the sets $\{\, h \in \se X T \,\vert \, h(a) = 1 \}$ 
for $a \in G$, we get at once the equivalence of (1) and (2). By taking complements and replacing 
\y a by $-a$, (3) is equivalent to (2).
\vspace{-0.17cm}

\nhp (1)/(3) $\Ra$ (4). For the first assertion, if $g(a) = 0$, (3) 
gives $h(a) \in \{0,1\}$, but (2) precludes $h(a) = 1$. For the second, 
(2) gives $h(a) = 1\; \Ra \; g(a) = 1$; if $h(a) = -1$, just replace \y a by $-a$.
\vspace{-0.17cm}

\nhp (4) $\Ra$ (5). The identity $h = h^2g$ obviously holds if $h(a) = 0$; 
if $h(a) \neq 0$, it follows from the second assertion in (4) and $h^2(a) = 1$.
\vspace{-0.17cm}

\nhp (5) $\Ra$ (2). $h = h^2g$ and $h(a) = 1$ clearly imply $g(a) = 1$. \hfl $\Box$

\vspace{-0.1cm}

We also register the following algebraic characterizations of inclusion 
and equality of zero-sets of elements of $\se X T$.

\vspace{-0.4cm}

\ble \label{char-zeroset} Let $T$ be a TS, and let $u,g,h \in \se X T$. Then:
\vspace{-0.15cm}

\nhp $(1)$ $Z(g)\, \sub\, Z(h)\; \Lra\; h = h g^2$.
\vspace{-0.15cm}

\nhp $(2)$ $Z(g) = Z(h)\; \Lra\; g^2 = h^2$.
\vspace{-0.15cm}

\nhp $(3)$ If\, $u \spez \, g, h$, then\, $Z(g)\, \sub \, Z(h)$\, if and only if\, 
$g \spez \, h$.
\ele

\vspace{-0.5cm}
\nhp {\bf Proof.} (1) and (2) are straightforward.

\vspace{-0.15cm}

\nhp (3) Lemma \ref{char-specializ}\,(4) proves the implication ($\La$) and that
$u \spez \, g, h$ implies\, $Z(u)\, \sub \, Z(g)\, \cap\, Z(h)$. 

\vspace{-0.15cm}

\nhp ($\Ra$) Assuming\, $Z(g)\, \sub \, Z(h)$, it suffices to verify 
the second clause of \ref{char-specializ}\,(4). Let $a \not\in Z(h)$. 
Then, $a \not\in Z(g)$, and the equivalence of items (1) and (4) in 
\ref{char-specializ}, together with\, $u \spez \, g$\, and\, $u \spez \, h$ 
yields\, $u(a) = g(a)$\, and\, $u(a) = h(a)$, respectively. Thus, 
$g(a) = h(a)$, as required. \hfl $\Box$ 

\vspace{-0.1cm}

\noi {\bf Remark.} Contrary to the case of real semigroups (\cite{M}, Prop. 6.4.1, p. 114), 
the specialization order of the character space of arbitrary ternary semigroups may not be 
a root system. A counterexample is given in \cite{DP1}, Ex. 1.14, p. 105. However, it is a
normal space in the usual topological sense of this notion (\cite{DP4}, Prop. I.1.21). \hfl $\Box$

\vspace{-0.4cm}

\bfa \label{zeros} Let $G$ be a ternary semigroup, let $X \subseteq$ \HTS G, 
and assume that $(X,G)$\vspace{0.06cm} is a q-fan. A necessary condition for 
$(X,G)$ to be an ARS is that for all \,$a,b \in G$, either $Z(a) \subseteq Z(b)$ 
or $Z(b) \subseteq Z(a)$. Here, $Z(a) = \{\, h \in X\, \vert \, h(a) = 0 \}$.
\efa
\vspace{-0.4cm}

\noindent {\bf Proof.} Assume $(X,G) \models$\,ARS but there are $a,b \in G$ so that 
$Z(a) \not \subseteq Z(b)$ and $Z(b) \not \subseteq Z(a)$, i.e., $\se h 1(a) = 0,\: 
\se h 1(b) \neq 0,\: \se h 2(b) = 0,\: \se h 2(a) \neq\, 0$, for some $\se h 1, 
\se h 2 \in X$. Since $(X,G)$ is a q-fan, $\se {h^2} 1 \se h 2 \in X$. By \cite{M}, Prop.
6.1.5, p. 103, $ \Dt G(a^2,b^2) = \{c^2\}$ for some $c \in G$, and 
$Z(a) \cap Z(b) = Z(c)$. This entails $\se h 1(c) = \se h 1(b) \neq 0$ and 
$\se h 2(c) = \se h 2(a) \neq 0$, whence $\se {h^2} 1 \se h 2(c) \neq 0$, contradicting
that $\se {h^2} 1 \se h 2 \in Z(a) \cap Z(b) = Z(c)$.  \hfl $\Box$
\vspace{-0.3cm}

\bfa \label{zerosetsElts} Let $T$ be a ternary semigroup and let $a,b \in T$. The followng are
equivalent:
\vspace{-0.15cm}

\noi $(1)\;  Z(a) \subseteq Z(b)$;  \hspace*{0.4cm} $(2)\;  a^2b^2 = b^2$;
\hspace{0.4cm}  $(3)\;  a \mid b$ $($i.e., $a$ divides $b$, i.e., $b = ax$ for some $x \in G).$
\vspace{-0.15cm}

\noi $(4)\; \se I b \subseteq \se I a$, where $\se I c = c \cdot T$ is the ideal of T generated by $c \in T$.
\efa

\vspace{-0.5cm}

\noindent {\bf Proof.} (4) $\Lra$ (3) $\Ra$ (1) are clear. For (2) $\Ra$ (3), scaling the equality (2) by $b$
yields $a^2b = b$, i.e., $b = a(ab)$.  For the implication (1) $\Ra$ (2) we use the separation 
theorem for TS's, \cite{DP1}, Thm. 1.9:  if $a^2b^2 \neq b^2$, there is $h \in \HTS T$ such 
that $h(a^2b^2) \neq h(b^2)$; hence $h(b^2) = 1$ and $h(a^2b^2) = h(a^2) = 0$, i.e., 
$h(a) = 0$ and $h(b) \neq 0$, i.e., $Z(a) \not\subseteq Z(b)$. \hfl $\Box$

Our next result gives alternative characterizations of the necessary condition in Fact \ref{zeros}.

\vspace{-0.4cm}

\bpr \label{charfan} Let $T$ be a ternary semigroup. The following conditions 
are equivalent:
\vspace{-0.15cm}

\noindent $(1)$ The family $\{Z(a)\, \vert \, a \in T \}$ is totally ordered 
under inclusion.
\vspace{-0.15cm}

\noindent $(2)$ For all \,$a,b \in T$, either $a^2b^2 = a^2$ or \,$a^2b^2 = b^2$.
\vspace{-0.15cm}

\noindent $(3)$ For all \,$a,b \in T$, either $b \mid a$ or $b \mid a$.
\vspace{-0.15cm}

\noindent $(4)$ Every proper ideal of \,$T$ is prime $($i.e., $ab \in I\; \Ra\; a \in I$ or $b \in I)$.
\vspace{-0.15cm}

\noindent $(5)$ The set of ideals of \,$T$ is totally ordered under inclusion.
\epr

\vspace{-0.45cm}

\noindent {\bf Proof.} The equivalence of (1) -- (3) follows immediately from \ref{zerosetsElts}.
\vspace{-0.15cm}

\noindent $(2) \Rightarrow (4)$. Let $I$ be an ideal of $T$, and suppose $ab \in I$; then $a^2b^2 \in I$ 
and, by (2), $a^2 \in I$ or $b^2 \in I$, which implies $a \in I$ or $b \in I$ (as $x = x^3 = x^2x$).
\vspace{-0.15cm}

\noindent $(4) \Rightarrow (5)$. If $J_1, J_2$ are incomparable ideals, then $J_1 \cap J_2$ is not prime 
(if $a \in J_2 \setminus J_1,\, b \in J_1 \setminus J_2$ then $ab \in J_1\, \cap J_2$ but 
$a,\, b \not\in J_1\, \cap J_2$).
\vspace{-0.15cm}

\noindent $(5) \Rightarrow (2)$. Given \,$a,b \in T$, by (5), either $\se I a \subseteq \se I b$ or 
$\se I b \subseteq \se I a$, and \ref{zerosetsElts} yields $b \mid a$ or $a \mid b$. \hfl $\Box$
\vspace{-0.5cm}

\paragraph{B. Real semigroups.}  \label{ReSe}

For easy reference we state the axioms defining real semigroups.  The
language for real semigroups, denoted $\se {\cal L} {\mathrm {RS}}$, is that of ternary 
semigroups ($\{ \cdot\,, 1, 0, -1 \}$) enriched with a ternary relation $D$. In agreement with 
standard notation (cf. \cite{M}, p. 99 ff.), we write $a \in D(b,c)$ instead of  $D(a,b,c)$. We set:
\vspace{-0.15cm}

\nhp [t-rep] \hfl{$a \in D^t(b,c) \Leftrightarrow a \in D(b,c) \wedge 
-b \in D(-a,c) \wedge -c \in D(b, -a).$} \hfl
\vspace{-0.15cm}

\nhp The relations $D$ and $D^t$ are called {\bf representation} and 
{\bf transversal representation}, respectively. 
\vspace{-0.45cm}

\bdf \label{RS}  A {\bf real semigroup} $($abbreviated RS$)$ is a ternary semigroup together 
with a ternary relation $D$ satisfying the following axioms:
\vspace{-0.15cm}

\nhp [RS0]\,\, $c \in D(a,b) \; \mbox{if and only if}\; c \in D(b,a)$.
\vspace{-0.15cm}

\nhp [RS1]\,\, $a \in D(a,b)$.
\vspace{-0.15cm}

\nhp [RS2]\,\, $a \in D(b,c) \; \mbox{implies}\;\; ad \in D(bd,cd)$.
\vspace{-0.15cm}

\nhp [RS3] ({\em Strong associativity}\,)\,\, If $a \in D^t(b,c) \; \mbox{and}\; c \in D^t(d,e)$, 
then there exists $x \in D^t(b,d)$ such that $a \in D^t(x,e)$.
\vspace{-0.15cm}

\nhp [RS4]\,\, $e \in D(c^2a, d^2b) \;\; \mbox{implies}\;\; e \in D(a,b)$.
\vspace{-0.15cm}

\nhp [RS5]\,\, If $ad = bd, ae =be, \; \mbox{and}\; c \in D(d,e), \;\; 
\mbox{then}\;\; ac = bc$.
\vspace{-0.15cm}

\nhp [RS6]\,\, $c \in D(a,b) \;\; \mbox{implies}\;\; c \in D^t(c^2a, c^2b)$.
\vspace{-0.15cm}

\nhp [RS7] ({\em Reduction}\,)\,\, $D^t(a,-b) \cap D^t(b,-a) \neq \emptyset \;\; \mbox{implies}\;\; a = b$.
\vspace{-0.15cm}

\nhp [RS8]\,\, $a \in D(b,c) \;\;\mbox{implies}\;\; a^2 \in D(b^2,c^2)$. \hfl $\Box$
\edf

\vspace{-0.35cm}

\noi For a detailed treatment of real semigroups the reader is referred to \cite{DP1}, \S\S\,2--4, pp. 106-119, 
or \cite{DP2}, \S\,2, pp. 57-59. Many of the properties of RSs and their duals, the {\it abstract real spectra} 
(abridged ARS) appearing in these and other references (e.g., \cite{M}, Chs. 6--8) are used below. For easy
reference we state the following fundamental:
\vspace{-0.45cm}

\bex \label{Ex3} (\cite{DP1}, Corollary 2.4, p. 109) The ternary semigroup ${\bf 3} = \{ 1,0, -1 \}$ has 
a {\it unique} structure of real semigroup, with representation given by:
\vspace{-0.2cm}

\nhp $\se D {\bf 3}(0,0) = \{ 0 \}; \hspace{1cm} \se D {\bf 3}(0,1) = \se D {\bf 3}(1,0) = \se D {\bf 3}(1,1) = \{ 0,1 \}$;
\vspace{-0.15cm}

\nhp $\se D {\bf 3}(0,-1) = \se D {\bf 3}(-1,0) = \se D {\bf 3}(-1,-1) = \{ 0, -1 \}$; 
\hspace{1cm} $\se D {\bf 3}(1,-1) = \se D {\bf 3}(-1,1) = {\bf 3}$;
\vspace{-0.1cm}

\nhp and transversal representation given by:
\vspace{-0.1cm}

\nhp $\se {D^t} {\bf 3}(0,0) = \{ 0 \}; \hspace{1cm} \se {D^t} {\bf 3}(0,1) = 
\se {D^t} {\bf 3}(1,0) = \se {D^t} {\bf 3}(1,1) = \{ 1 \}$;
\vspace{-0.15cm}

\nhp $\se {D^t} {\bf 3}(0,-1) = \se {D^t} {\bf 3}(-1,0) = \se {D^t} {\bf 3}(-1,-1) = \{ -1 \}$;
\hspace{0.95cm} $\se {D^t} {\bf 3}(1,-1) = \se {D^t} {\bf 3}(-1,1) = {\bf 3}$. \hfl $\Box$ 
\eex

\vspace{-0.35cm}

As we will see in \S\,2 below, the crucial axiom in the list above is the strong associativity axiom [RS3].
We register (\cite{M}, Prop. 6.1.1, p. 100, and Thm. 6.2.4, pp. 107-108) that [RS3] is equivalent 
to the conjunction of the weak associativity axiom [RS3a] obtained by replacing transversal 
representation by ordinary representation in [RS3] and

\nhp [RS3b]\;  For all $a,b,\; \Dt{}(a,b) \neq \0$.
\vspace{-0.15cm}

The following equivalent form of [RS3] turns out to be very useful at the time of checking
strong associativity in concrete examples:
\vspace{-0.4cm}

\bpr \label{axRS3'} In the presence of axiom {\em [RS2]}, the following is 
equivalent to axiom {\em [RS3]}:
\vspace{-0.15cm}

\nhp {\em [RS3$'$]}\;\; $\fa\,a,b,c,d\:(D^t(a,b) \cap D^t(c,d) \neq\, \0\;\; \Ra\;\; 
D^t(a,-c) \cap D^t(-b,d) \neq\, \0).$
\epr

\vspace{-0.45cm}

\nhp {\bf Proof.} [RS3] $\Ra$ [RS3$'$]. Let $x \in D^t(a,b) \cap D^t(c,d)$; 
the definition of $D^t$ yields $-b \in D^t(a,-x)$ (\cite{DP1}, Prop. 2.3\,(0), p. 107), 
and scaling by $-1$ ([RS2]) gives $-x \in D^t(-c,-d)$. By [RS3] there is 
$y \in D^t(a,-c)$ so that $-b \in D^t(y,-d)$. Again, the definition of 
$D^t$ and [RS2] yield $-y \in D^t(b,-d)$, and $y \in D^t(-b,d)$. Hence, 
$D^t(a,-c) \cap D^t(-b,d) \neq\, \0).$
\vspace{-0.15cm}

\nhp [RS3$'$] $\Ra$ [RS3]. Assume [RS3$'$] and let $x \in D^t(a,b)$ 
with $b \in D^t(c,d)$. By the definition of $D^t$, $-b \in D^t(a,-x)$, 
and by [RS2], $b \in D^t(-a,x)$, i.e., $D^t(-a,x) \cap D^t(c,d) \neq\, \0$. 
By [RS3$'$] there is $y \in D^t(-a,-c) \cap D^t(-x,d)$. By the same 
manipulation as above, we get $-y \in D^t(a,c)$ and $x \in D^t(-y,d)$. 
So, [RS3] is verified with witness $-y$. \hfl $\Box$

\nhp {\bf Remark.} Note that, while the weak associativity axiom [RS3a] 
 is a non-trivial property (in the sense that it {\it is not} a consequence of 
 the remaining axioms), the corresponding weak version of [RS3$'$], obtained
 by replacing $D$ for $D^t$, {\it does follow} from the remaining axioms for RSs, 
 as [RS1] and [RS4] imply $0 \in D(a,b)$ for all \y a, \y b (\cite{DP1}, Prop. 2.3\,(1), 
 p. 107). \hfl $\Box$

\noi {\bf The representation partial order on real semigroups.} Recall (\cite{DM1}, Cor. 4.4\,(c), p. 62) 
that for a reduced special group, $G$, the binary relation $a \leq b\; \Lra\; a \in \D G (1,b)$ is a partial 
order for which the operation ``multiplication by $-1$" is an involution. Further, this relation is induced 
from the partial order of the Boolean hull of $G$ (\cite{DM1}, Cor. 4.12, p. 69).
\vspace{-0.15cm}

In the context of RS's, none of the binary relations $a \in D(1,b)$ or $a \in D^t(1,b)$ 
defines a partial order for which the operation ``--" (multiplication by $-1$) is an involution. 
However, since every RS, $G$, is canonically embedded in a Post algebra (seen as a RS, 
its ``Post hull", \cite{DP2}, Prop. 4.1, p. 62), which is a distributive lattice, the 
latter induces a partial order on $G$ given by:

\vspace{-0.5cm}

\bdf \label{reprord} (\cite{DP2}, Rmk. 2.5, p. 59) Let $G$ be a RS, and let $a,b \in G$. 
We set:
\vspace{-0.15cm}

\nhp \centerline{$a\: \se {\leq} {G}\: b \;\;\;\;\mbox{iff}\;\;\;\; a \in \se D {G}(1,b)
\;\mbox{and}\; -b \in \se D {G}(1,-a)$.} 

\vspace{-0.15cm}

\nhp [Unless necessary we omit the subscript in $\se {\leq} {G}$.] This relation is a
partial order on $G$ (\ref{reprpo-RS}\,(1)) called the {\bf representation partial order}.
\hfl $\Box$ 
\edf

\vspace{-0.5cm}

\noi When\, $G = {\bf 3}$\, this definition gives $1\, \se {<}{\bf 3}\, 0\, \se {<}{\bf 3} -1$, 
the \underline{opposite}\vspace{0.05cm} of the order of these elements as integers. The binary 
relation just defined has the following properties:

\vspace{-0.45cm}

\bth\label{reprpo-RS} Let $G$ be a RS. For\, $a,b,x,y \in G$ we have:
\vspace{-0.15cm}

\nhp $(1)$ The relation $\leq$ is a partial order on $G$ such that $a \leq b\;
\Lra -b \leq -a$.
\vspace{-0.15cm}

\nhp $(2)$ For all $a \in G, \;\; 1 \leq a \leq -1$.
\vspace{-0.15cm}

\nhp $(3)$ \;$a \leq 0 \;\; \Lra \;\; a =\, a^2 \in \mbox{\em Id}(G)$,
\vspace{-0.2cm}

\nhp $\hspace{0.7cm} 0 \leq a \;\; \Lra \;\; a = -a^2 \in -\mbox{\em Id}(G)$.
\vspace{-0.15cm}

\nhp $[$\,\mbox{\em Id}$(G) = \{a^2 \sth a \in G\}$ is the set of {\em idempotents} of $G$.$]$
\vspace{-0.15cm}

\nhp $(4)$ Let $\se X G$ be the character space of $G$. For $a,b \in G$,
\vspace{-0.35cm}
$$
\begin{array}{ll} \hspace{0.5cm}  a\: \se {\leq} {G}\: b & \Lra \;\;\; \forall\, h \in
\se X {G}\ (h(a)\, \se {\leq} {\bf 3}\, h(b))\; \;\,\Lra\\ [0.1cm]
& \Lra \;\;\; \forall\, h \in \se X {G}\ [(h(b) = 1 \; \Ra \; h(a) = 1) \wedge (h(b) =
0 \; \Ra \; h(a) \in  \{ 0,1 \}) ]. 
\end{array}
$$

\vspace{-0.45cm}

\nhp $(5)$ The following are equivalent:
\vspace{-0.15cm}

\nhp \hmm\hem \mbox{\em (i)}\; $a^2 \leq b \leq -a^2$; \hmm\hem \mbox{\em (ii)}\; $Z(a) \subseteq Z(b)$; 
\hmm\hem \mbox{\em (iii)}\; $b = a^2b$.
\vspace{-0.2cm}
 
\nhp In particular,
\vspace{-0.2cm}

\nhp $(6)$\; $a^2 \leq ab \leq -a^2$\; $($hence\, $a^2 \leq \pm\, a \leq -a^2)$.
\vspace{-0.15cm}

\nhp $(7)$ If\; $a^2 \leq b \leq -a^2$\, and\, \y b\, is invertible, then\, \y a\, is invertible.
\vspace{-0.15cm}

\nhp $(8)$ $a \leq x,y\; \Ra\; a \leq -xy$. Hence, $x, y \leq a\; \Ra\; xy \leq a$.
\vspace{-0.15cm}

\nhp $(9)$ For all $a \in G$, the infimum and the supremum of $a$ and
$-a$ for the representation partial order $\leq$ exist, and $a\, \w -a
= a^2,\, a\: \jo -a = -a^2$. In particular,
\vspace{-0.2cm}

\nhp $(10)$ $a\, \w -a \leq 0 \leq b\: \jo -b$ for all $a,b \in G$.\,\footnote{\;Called
the {\it Kleene inequality}; cf. \cite{DP2}, Rmk. 1.2\,(b), p. 55.} \hfl $\Box$
\eth

\vspace{-0.35cm}

\noi The proof of Theorem \ref{reprpo-RS} appears in \cite{DP4}, Propositions I.6.4, I.6.5.
\vspace{-0.4cm}

\section{From ternary semigroups to\vspace{-0.4cm} real semigroups} \label{TStoRS}

In this short section we show that any non-empty set of TS-characters of a ternary semigroup
induces ternary relations that satisfy \und{all} RS-axioms except, possibly, the strong associativity
axiom [RS3] (cf. Definition \ref{RS}).   
\vspace{-0.4cm}

\bdf \label{ReprD_H} Given a ternary semigroup, $G$, and a set $\cH\, \sub\, \se X G = \HTS G$, 
we define a ternary relation $\D {G,\cal H}$ on $G$ ---abridged $\D {\cal H}$ if $G$ is clear 
from context--- as follows: for $a,b,c \in G$,

\vspace{-0.15cm}
\nhp $\se {[D]} {\cal H}$ \hfl $a \in \D {G,\cal H}(b,c)\;\; \Lra\;\;$ For all 
$h \in \cH,\; h(a) \in \Dthr(h(b), h(c))$. \hfl
\vspace{-0.2cm}

\nhp To avoid triviality we assume $\cH \neq \0$; to get best results we also make
the rather mild assumption that the set $\cH$ separates points in $G$: given $a \neq b$
in $G$, there is $h \in \cH$ such that $h(a) \neq h(b)$. In a similar way a corresponding 
``transversal'' relation is defined by:

\vspace{-0.15cm}
\nhp $\se {[D^t]} {\cal H}$ \hfl $a \in \Dt {G,\cal H}(b,c)\;\; \Lra\;\;$ For all 
$h \in \cH,\; h(a) \in \Dtthr(h(b), h(c))$. \hfl $\Box$
\edf

\vspace{-0.75cm}

\bre \label{IdenticalReprRelats} Given a TS,\, \y G, and a set\, $\cH\, \sub\: {\bf 3}^G$, 
in \cite{M}, p. 99, Marshall defines representation relations on\, \y G, as follows: for $a,b,c \in G$,
\label{RSrepr} 

\vspace{-0.45cm}

\begin{tabbing}
[R] $\;\;\; a \in \se D {\cal H}(b,c) \;\;$ \=iff \= $\;\; \fa \, h \in \cH \; [h(a) = 0 \; \jo \; 
(h(a) \neq 0 \; \w \; (h(a) = h(b) \, \jo \, h(a) = h(c)))]$.\\ [0.2cm]
\noindent [TR] $\,a \in \se {D^t} {\cal H}(b,c) \;\;$ \>iff \> $\;\; \fa \, h \in \cH \; 
[(h(a) = 0 \; $\=$ \w \; h(b) = -h(c)) \; \jo\; (h(a) \neq\, 0 \; \w$ \\ 
\> \> \> $\w\; (h(a) =  h(b) \, \jo \, h(a) = h(c)))]$.
\end{tabbing}

\vspace{-0.45cm}

The representation relation\, $\se D {\cal H}$\, defined by clause [R] 
\und{is identical} with the relation defined by clause\, $\se {[D]} {\cal H}$ in \ref{ReprD_H}. 
This is obvious by the fact that conditions\; $x \in \se D {\bf 3}(y,z)$\; and\; 
$x = 0\; \jo$ $(x \neq 0 \; \w \; (x = y \, \jo \, x = z))$\; are equivalent 
for all\, $x,y,z \in {\bf 3}$; this is straightforward checking using Example \ref{Ex3}.
\vspace{-0.15cm}

Likewise, the transversal representation relation\, $\Dt {\cal H}$\, defined by [TR] is 
identical to the transversal representation relation defined in terms of\, $\se D {G,{\cal H}}$\, 
by clause [t-rep], Section \ref{ReSe}, since conditions\; 
$x \in \se {D^t} {\bf 3}(y,z)$\; and\; $((x = 0 \; \w \; y = -z) \; \jo\; 
(x \neq\, 0 \; \w \;(x =  y \, \jo \, x = z)))$\; are equivalent, again by \ref{Ex3}. \hfl $\Box$
\ere

\vspace{-0.4cm}

\bth \label{RSaxioms} Let\, $G$\, be a ternary semigroup and let\, $\cH$\, be a non-empty subset 
of\, $\se X {G}$ separating points in $G$. The ternary relation\, $\se D {\cal H}$\, defined 
in \ref{ReprD_H} satisfies all axioms for real semigroups except, possibly, the axiom {\em [RS3]} 
of strong associativity.
\eth

\vspace{-0.45cm}

\nhp {\bf Proof.} The verification of axioms [RS0], [RS1], [RS2], [RS4] and [RS8] being 
straightforward, we deal only with the remaining axioms.
\vspace{-0.15cm}

\nhp [RS5]\; Let\, $a,b,c,d,e \in G$\, be such that\, $ad = bd,\; ae = be$\, and\, 
$c \in \D {\cal H}(d,e)$. Let us prove that\, $ac = bc$. Since $\cH$ separates points
in $G$, this boils down to proving $h(ac) = h(bc)$\, for all\, $h \in \cH$. This is clear 
if\, $h(c) = 0$. Let\, $h(c) \neq\, 0$. Since\, $c \in \D {\cal H}(d,e)$, either\, 
$h(c) = h(d)$\, or\, $h(c) = h(e)$. Since\, $ad = bd$\, and\, $ae = be$, invoking 
Definition \ref{ReprD_H}, in both cases we get the equality\, $h(ac) = h(bc)$. By\, 
$\se {[D]} {\cal H}$\, once again, we conclude that\, $ac = bc$, as required.
\vspace{-0.15cm}

\nhp [RS6]\; Let\, $a,b,c \in G$\, be such that\, $c \in \D {\cal H}(a,b)$, and take\, 
$h \in \cH$. Then,\, $h(c) \in \se D {\bf 3}(h(a),h(b))$. The real semigroup\, ${\bf 3}$\, 
verifies [RS6], and then\, $h(c) \in D^{t}_{\bf 3}(h(c)^{2}h(a),h(c)^{2}h(b))$. From the 
definition of\, $D^t$\, (cf. \S\,\ref{ReSe}, [t-rep]),  we have the following relations:
\vspace{-0.15cm}

\nhp $(i)$\;\; $h(c) \in \se D {\bf 3}(h(c^{2}a),h(c^{2}b)),$ \hmm $(ii)$\, 
$-h(c^{2}a) \in \se D {\bf 3}(-h(c),h(c^{2}b)),$\hmm and 
\vspace{-0.15cm}

\nhp $(iii)$ $\!-h(c^{2}b) \in \se D {\bf 3}(-h(c),h(c^{2}a)).$
\vspace{-0.15cm}

\noindent Since\, $h$ is arbitrary, from (i), (ii), (iii) and \ref{ReprD_H}\,$\se {[D]} {\cal H}$\, 
we get:
\vspace{-0.15cm}

\nhp $(i')$\;\, $c \in \D {\cal H}(c^{2}a,c^{2}b)$,\, \hmm $(ii')$\, 
$-c^{2}a \in \D {\cal H}(-c,c^{2}b),$
\vspace{-0.15cm}

\nhp $(iii')$ $-c^{2}b \in \D {\cal H}(-c,c^{2}a)$,
\vspace{-0.15cm}

\noindent which, together, amount to\, $c \in \se {D^{t}} {G/{\cal H}}(c^{2}a), c^{2}b)$.
\vspace{-0.15cm}

\noindent [RS7]\; Let\, $a,b \in G$\, be such that\, $\Dt {\cal H}(a,-b) \cap 
\Dt {\cal H}(b,-a) \neq\, \emptyset$. Take an element\, $c \in G$\, in this intersection. 
We must prove that\, $a = b$. By \ref{ReprD_H}\, $\se {[D]} {\cal H}$\, this boils down 
to showing that\, $h(a) = h(b)$\, for all\, $h \in \cH$. We consider the following cases:
\vspace{-0.15cm}

\noindent (i) $h(c) = 0$. If either\, $h(a) \neq  0$\, or\, $h(b) \neq 0$, from the 
relations\, $-a \in \D {\cal H}(-c,-b)$\, and\, $-b \in \D {\cal H}(-c, -a)$\, 
we obtain\, $h(-a) = h(-b)$, and then $h(a) = h(b)$. If\,$h(a) = h(b) = 0$, there is 
nothing to prove.
\vspace{-0.15cm}

\noindent (ii) $h(c) \neq 0$. Since\, $c \in \D {\cal H}(a,-b) \cap  \D {\cal H}(b,-a)$, 
we have\, $h(c) = h(a)$\, or\, $h(c) = -h(b)$, and\, $h(c) = h(b)$\, or\, $h(c) = -h(a)$. 
If\, $h(a) \neq\, h(b)$, these conditions yield either\, $h(c) = h(a) = -h(a)$\, or\, 
$h(c) = h(b) = -h(b)$; in both cases we have\, $h(c) = 0$, a contradiction. Hence, 
$h(a) = h(b)$. \hfl $\Box$

\noi {\bf Remark.} Theorem \ref{RSaxioms} reduces the question of checking whether
a TS with a ternary relation defined in terms of characters as in \ref{ReprD_H} above
is a RS, to checking whether the single axiom [RS3] holds. For example, when the set
$\cH$ of characters has 2 elements, or has 3 elements with a non-trivial specialization,
axiom [RS3] holds. However, there are examples of finite sets of characters $\cH$ for
which [RS3] \und{does not} hold in $\se G {\cal H}$ (\cite{DP4}, \S\,I.3). \hfl $\Box$ 
\vspace*{-0.45cm}

\section{Fans are real semigroups and abstract real spectra} \label{RSfan}
\vspace{-0.4cm}

Our main aim in this section  is to prove that (ARS-)fans, in any of the two senses 
considered in Definition \ref{fans-various}, are abstract real spectra.  The 
first step to achieve this is to work out the explicit form of the representation 
relations corresponding to the notion of ``q-fan" (under the assumption  that 
the necessary condition in \ref{zeros} is verified); this is done in Theorem 
\ref{fanrepr}. It follows that any TS verifying this necessary condition and 
endowed with the relations thus obtained is a real semigroup (Theorem \ref{fanRS}).
A number of results (\ref{equivfan} -- \ref{fanideals}) which determine   
to a large extent the structure of RS-fans follow from these theorems.

\vspace{-0.5cm}

\bth \label{fanrepr} Let \y G be a ternary semigroup verifying
\vspace{-0.15cm}
 
\nhp $[Z] \;\;\;\; \fa \, a,b \in G \, (Z(a)\, \sub\, Z(b)$ or $Z(b)\, \sub\, Z(a))$.
\vspace{-0.15cm}

\nhp Let  \y X $\sub$\ \HTS G be such that $(X,G)$ is a q-fan. With
$D = \se D {\!X}$ and $D^t = \se {D^t} {\!X}$ denoting the representation relations
defined by clauses $[R]$ and $[TR]$ in \ref{IdenticalReprRelats}, for $a,b \in G$ we have:

\nhp $[D^t] \;\;\;
D^t(a,b) =
\left\{
\begin{array}{ll}
\{ a \} &  \;\mbox{if}\;\;\; Z(a) \subset Z(b)\\
\{ b \} &  \;\mbox{if}\;\;\; Z(a) \subset Z(b)\\
\{ a,b \} &  \;\mbox{if}\;\;\; Z(a) = Z(b) \;\, and \;\, b \neq -a\\
\,a \cdot G \,(= b \cdot G) &  \;\mbox{if}\;\;\; b = -a.
\end{array}
\right.
$

\nhp $[D]  \;\;\;\;  D(a,b) = a \cdot {\mbox{$\mathrm{Id}$}}(G) \, \cup \, b \cdot
{\mbox{$\mathrm{Id}$}}(G) \, \cup \, \{ x \in G \, \vert \, xa = -xb \; \w \; x = a^2x 
\}$.
\eth

\vspace{-0.7cm}

\bre \label{RmkOn[D]PrecThm} The inclusion $\supseteq$ in item
$[D]$ of \ref{fanrepr} holds for an arbitrary TS, $G$, and any set 
$X\, \sub\, \HTS G$
(indeed, it follows from axioms [RS1] and [RS4]): use Definition \ref{ReprD_H}, and 
that $\bf 3$ is a real semigroup (Example \ref{Ex3}).  \hfl $\Box$
\ere
\vspace{-0.5cm}

\bth \label{fanRepresInterdefinable} Let \y G be a ternary semigroup verifying 
condition $[Z]$ of Theorem $\ref{fanrepr}$. Then, conditions $[D]$ and $[D^t]$ in 
$\ref{fanrepr}$ are interdefinable in the following sense:
\vspace{-0.15cm}

\nhp $(1)$ Assuming that a ternary relation $D$ on $G$ is \und{defined} as in
$[D]$ and the corresponding transversal representation is given by the clause
\vspace{-0.15cm}

\nhp \centerline{$a \in D^t(b,c) \Leftrightarrow a \in D(b,c) \wedge -b \in
D(-a,c) \wedge -c \in D(b, -a),$}

\vspace{-0.15cm}

\nhp then $D^t$ verifies condition $[D^t]$ of $\ref{fanrepr}$.

\nhp $(2)$ Conversely, if $D^t$ is \und{defined} as in $[D^t]$ and the
associated ternary representation relation $D$ is defined by the stipulation $a
\in D(b,c)  \Leftrightarrow a \in D^t(a^2b,\,a^2c)$, then $D$ verifies clause
$[D]$ of $\ref{fanrepr}$.
\eth

\vspace{-0.7cm}

\bth \label{fanRS} Let \y G be a ternary semigroup verifying condition $[Z]$ of
Theorem \ref{fanrepr}. With the ternary relation $D$ defined as in $\ref{fanrepr}$, 
$(G,D)$ is a real semigroup. 
\eth

\vspace{-0.5cm}

\noi This result is a natural extension of \cite{Li}, Prop. 1.1.14 (see 
\ref{remksAOSFans}\,(b)) to the theory of real semi\-groups.  

Before engaging in the proof of these theorems we draw some important
consequences of them.
\vspace{-0.4cm}

\bpr \label{equivfan} Let $G$ be a TS verifying condition $[Z]$ of Theorem
$\ref{fanrepr}$, and $X \; \sub\, \HTS G$. The following are equivalent:
\vspace{-0.1cm}

\nhp $(1)$ $(X,G) \models \se {\mathrm{fan}} {\,1}$ \, $($i.e., \y X = \HTS G$)$.
\vspace{-0.15cm}

\nhp $(2)$\; $i)$ $(X,G)$ is a q-fan.
\vspace{-0.2cm}

\noi \hspace*{0.5cm} $ii)$ For every subsemigroup \y S of\, \y G such that 
$S\, \cup \, -S = G$ and $S \, \cap \, -S$ is a $($proper$)$ prime 

\vspace{-0.2cm}

\noi \hspace*{1.05cm} ideal, there is $h \in X$ such that $S = h^{-1}[\,0,1]$.
\vspace{-0.2cm}

\nhp $(3)$ $(X,G) \models \se {\mathrm{fan}} {\,2}$.
\epr
\vspace{-0.4cm}

\nhp {\bf Proof.} (3) $\Ra$ (2). Assumption (3) implies that $(X,G)$ is an ARS. By axiom
[AX1] (\cite{M}, p.\, 99), $X$ separates points of $G$ and is closed under product of any
three of its members; so, (2.i) holds. By Theorem \ref{fanrepr}, the representation relation
$D$ is given by the equality [$D$] therein. Since $S \cup -S = G$ we have $S \supseteq {\mathrm{Id}}(G)$
and, by Corollary \ref{fanideals}\,(2), $S$ is saturated for $D$. So, axiom [AX2] in \cite{M}, p.\, 99
holds, and therefore $(X,G)$  verifies (2.ii).
\vspace{-0.15cm}

\nhp (2.ii) $\Ra$ (1). We must prove that \HTS G $\sub\, X $. Let $g \in$ \HTS G; 
set $S := g^{-1}[0,1]$. Clearly, this set verifies the assumptions in (2.ii). Hence, 
there is $h \in X$ so that $S = h^{-1}[0,1]$, whence $g^{-1}[0] = S \cap -S = h^{-1}[0]$. 
The equalities $g^{-1}[0,1] = h^{-1}[0,1]$ and $g^{-1}[0] = h^{-1}[0]$ entail $g = h$, 
and hence $g \in X$.
\vspace{-0.15cm}

\noi (1) $\Ra$ (3). Clearly, $X = \HTS G$ is closed under product of any three elements. To
prove it is an ARS we argue as follows. By the above and the separation theorem for TS,
\cite{DP1}, Thm. 1.9, pp. 103--104, $(X, \cl {G})$ is a q-fan\,\footnote{\;Recall that 
${\cl G} = \{{\h a} \sth a \in G\}$, where 
${\h a} \in {\bf 3}^{\,{\mathrm {Hom}_{\mathrm {TS}}}(G,{\bf 3})}$ is the map 
``evaluation at $a$'': for $\si \in \HTS G,\; {\h a}(\si) := \si(a)$.}. Since $G$ satisfies 
condition [Z], Theorem \ref{fanrepr} guarantees that the representation relation
$\se D X$ is given by the equality in \ref{fanrepr}\,(b). Then, Theorem \ref{fanRS}
proves that $(G, \se D X)$ is a RS. By the equivalence of (2) and (4) in Corollary 
\ref{ARSoffan}, its dual structure, $(X, \cl {G})$, is an ARS; hence, $(X,G)$ is a
$\se {\mathrm{fan}} {\,2}$. \hfl $\Box$

\vspace{-0.45cm}

\bdfn \label{def:fan}({\bf Fan})\, Henceforth we simply write ``{\bf fan}" 
(or ``{\bf ARS-fan}'') for either of the equivalent conditions 
{\bf $\se {\textrm{fan}} {1}$} or {\bf $\se {\textrm{fan}} {2}$}. In using 
the notation ``$(X,G) \models \textrm{fan}$'' we implicitly assume that the 
underlying ternary semigroup $G$ verifies condition [Z] in Theorem \ref{fanrepr}; 
this assumption is crucial and, in fact, distinguishes fans from many 
other classes of ARSs. We shall also say ``$G$ is a {\bf fan}" 
(or a ``{\bf RS-fan}"), tacitly assuming that its representation relations 
are those given in Theorem \ref{fanrepr}. \hfl $\Box$
\edfn
\vspace{-0.8cm}

\bco \label{fanmaps} Let $G$ be a TS verifying condition $[Z]$ of Theorem
$\ref{fanrepr}$. Let  $H$ be a real semigroup, and let $f : G \: \lra \, H$ be a
homomorphism of ternary semigroups. Then, $f$ preserves the representation
relation $D$ defined by clause $[D]$ of $\ref{fanrepr}$, and hence it is a
RS-homomorphism from $(G,D)$ into $H$. In other words, 
$\mathrm {Hom}{\raisebox{-5pt}{$\scriptstyle{\mathrm {RS}}$}}((G,D),H) =
\mathrm {Hom}{\raisebox{-5pt}{$\scriptstyle{\mathrm {TS}}$}}(G,H)$.
\eco
\vspace{-0.4cm}

\nhp {\bf Proof.} In view of the definition of $D$, the proof boils down to the
following obvious facts:
\vspace{-0.15cm}

\nhp (1) $f$ preserves products and idempotents, hence the clauses defining the
relation $D$.
\vspace{-0.15cm}

\nhp (2) For $a,b,c$ in an arbitrary RS, $H$, we have 
\vspace{-0.15cm}

\nhp (i) $\; a \cdot \textrm{Id}(H) \, \sub \, \se D H(a,b)$, \hem and \hem 
(ii) $\; ca = -cb \; \w \; c = a^2c \;\, \Ra \:\, c \in \se D H(a,b)$.
\vspace{-0.15cm}

\nhp Item (2.i) follows from axiom [RS4] of real semigroups. For the proof of (2.ii), 
observe that, for $h \in \se X H$, $h(c) \neq 0$ and $c = a^2c$ imply\vspace{0.05cm} 
$h(a) \neq 0$. If $h(c) \neq h(a)$, then \,$ca = -cb$\, yields $\ h(c) = h(b)$. This 
shows that $h(c) \in \se D {\bf 3}(h(a), h(b))$ for all  $h \in \se X H$; by the Separation 
Theorem for RSs (\cite{DP1}, Thm. 4.4, p. 116), we conclude that $c \in \se D H(a,b)$. 
\hfl $\Box$

\vspace{-0.15cm}

In particular we have:
\vspace{-0.45cm}

\bco \label{ARSoffan} Let $G$ be a TS verifying condition $[Z]$ of Theorem 
$\ref{fanrepr}$. Then,
\vspace{-0.1cm}

\noi $(1)$\; \HRS {(G,D)} = \HTS G.
\vspace{-0.15cm}

\noi Hence, 
\vspace{-0.15cm}

\noi $(2)$\; The ARS dual to the real semigroup $(G,D)$ is $($\HTS G, $\cl G)$.\!\!

\vspace{-0.15cm}

\noi $(3)$\; $($\HTS G, $\cl G)$ is a $\se {\mathrm {fan}} 1$ $($hence an ARS-fan, 
see $\ref{def:fan})$.
\eco
\vspace{-0.4cm}

\nhp {\bf Proof.} (1) is \ref{fanmaps} with $H = \bf 3$, and (2) comes from the 
definition of the ARS dual to any RS (see proof of the Duality Theorem 4.1, 
\cite{DP1}, p. 117). For (3), see \ref{fans-various}\,(1) and \ref{def:fan}. \hfl $\Box$
\vspace{-0.4cm}

\bco \label{unitsRS-fans} Let $(G,D)$ be a RS-fan. Then, the set 
$G^{\times} = \{a \in G \sth  a^2 = 1\}$ of invertible elements of $G$ 
with representation induced by restriction of $D$ to $G^{\times}$, is a RSG-fan, 
i.e., a fan in the category of reduced special groups.
\eco
\vspace{-0.45cm}

\nhp {\bf Proof.} Since $Z(a) = \0$ for $a \in G^{\times}$, only the last two
clauses in the characterization of $\Dt G$ given by Theorem \ref{fanrepr} apply, 
whenever $a, b \in G^{\times}$, and we have\,:
\vspace{-0.45cm}

$$\Dt G (a,b) =
\left\{
\begin{array}{ll}
\{ \, a,b \} &  \;\mbox{if}\;\;\; b \neq -a\,,\\
\;\;\; G &  \;\mbox{if}\;\;\; b = -a\,.
\end{array}
\right.
$$
\vspace{-0.45cm}

\nhp But this is exactly the definition of representation in a RSG-fan, cf. 
[RSG-fan] (Introduction). Axiom [RS6] for RSs implies 
$\Dt G (a,b)\, \cap\, G^{\times} = \DG (a,b)\, \cap\, G^{\times}$,
proving our contention. \hfl $\Box$
\vspace{-0.15cm}

\bco \label{fanideals} Let\, $G$ be a TS verifying condition $[Z]$ of Theorem
$\ref{fanrepr}$ and let $D$ be the ternary relation on $G$ defined by clause $[D]$
therein. Then,
\vspace{-0.15cm}

\nhp $(1)$ Every TS-ideal of\, \y G is a saturated prime ideal of the real semigroup
$(G,D)$. 
\vspace{-0.15cm}

\nhp $(2)$ A TS-subsemigroup \y S of\, \y G is saturated in $(G,D)$ iff it contains 
\mbox{$\mathrm {Id}(G) = \{ x^2 \, \vert \, x \in G \}$} and $S \cap -S$ is an ideal.
\eco

\vspace{-0.45cm}

\nhp {\bf Proof.} (1) Straightforward verification, using \ref{charfan}.
\vspace{-0.2cm}

\nhp (2) The implication ($\Ra$) is obvious. For the converse, write $I = S \,
\cap \, -S$; \y I is a  prime ideal (\ref{charfan}\,(3)). Let $a,b \in S$ and $c
\in D(a,b) = a \cdot \textrm{Id}(G) \, \cup $ 
$b \cdot \textrm{Id}(G) \, \cup \, \{ \, x \in\, G \, \vert \, xa = -xb \; \w \;
x = a^2x \, \}$. If $c \in a \cdot \textrm{Id}(G)$, then $c = ax^2$, whence $c
\in S$, since both \y a and $x^2$ are in \y S. The case $c \in b \cdot
\textrm{Id}(G)$ is similar. If $ca = -cb$ and $c = a^2c$, then $c^2a = -c^2b$,
which implies $c^2a \in I$. Since \y I is prime, either \y a or \y c are in \y
I; if $a \in I$, then $c = a^2c  \in I$; in either case we have $c \in I \, \sub
\,S$. \hfl $\Box$

\vspace{0.1cm}

\nhp {\bf Proof of Theorem \ref{fanrepr}.} First we prove:
\vspace{-0.15cm}

\nhp (A) If $(X,G) \models \textrm{q-fan}$, the relation $\se {D^t} {\!X}$ 
(defined by clause [TR], Remark \ref{IdenticalReprRelats}) verifies condition 
$[D^t]$ in the statement of the Theorem. 
\vspace{-0.15cm}

\noi Note that:
\vspace{-0.15cm}

\nhp (1) $Z(a) \, \sub \, Z(b) \; \Ra \; a \in \se {D^t} {\!X}(a,b)$ (immediate verification).
\vspace{-0.15cm}

\nhp  Next we prove: 
\vspace{-0.15cm}

\nhp (2) $Z(a) \, \sub \, Z(b) \, \w \; b \neq -a \;\, \Ra \,\, \se {D^t} {\!X}(a,b)
\, \sub \, \{ a,b \}$.
\vspace{-0.15cm}

\nhp  \und{Proof of (2)}. Suppose there is $c \in \se {D^t} {\!X}(a,b)$ such\vspace*{0.05cm} that $c
\neq a$ and $c \neq b$. Since $X$ separates points, these inequalities, together
with $b \neq -a$, give TS-characters $h_1,h_2,h_3 \in X$ whose images at the
points $a,b,c$ verify the corresponding inequalities in {\bf 3}. By assumption,
$h = h_1h_2h_3 \in X$, and we prove below that $h$ contradicts 
$c \in \se {D^t} {\!X}(a,b)$; more precisely, $h$ verifies either
\vspace{-0.15cm}

\nhp (*) \; $h(c) = 0$ \, and \, $h(a)\neq -h(b)$, \; or
\vspace{-0.15cm}

\nhp (**)\, $h(c) \neq 0$, \,  $h(c) \neq h(a)$\; and\; $h(c) \neq h(b)$.
\vspace{-0.1cm}

\nhp (I) Since $c \neq a$, there is $h_1 \in X$ so that $h_1(c) \neq h_1(a)$.
According to the values of $h_1(c) \in \{ 0,1,-1 \}$, conditions $c \in \se
{D^t} {\!X}(a,b)$ and $Z(a) \, \sub \, Z(b)$ yield the following alternatives:
\vspace{-0.15cm}

\nhp \hspace{0.5cm} I.a\,: \;\;\;\; $h_1(c) = 0$ \, and \, $h_1(a)h_1(b) = -1$\; or,
\vspace{-0.15cm}

\nhp \hspace{0.5cm} I.b\,: \;\;\;\; $h_1(c) = h_1(b) \neq 0$ \, and \, 
$h_1(a)h_1(b) = -1$.
\vspace{-0.1cm}

\nhp (II) Assumption $c \neq b$, yields a character $h_2 \in X$ so that
$h_2(c) \neq h_2(b)$. An analysis similar to that of (I) narrows the possible
values of $h_2$ at the points $a,b,c$ down to:
\vspace{-0.15cm}

\nhp \hspace{0.5cm} II.a\,: \;\;\; $h_2(c) = 0$ \, and \, $h_2(a)h_2(b) = -1$\; or,
\vspace{-0.15cm}

\nhp \hspace{0.5cm} II.b\,: \;\;\; $h_2(c) = h_2(a) \neq 0$ \, and \, $h_2(b) \in
\{ \, 0, -h_2(a)\}$.
\vspace{-0.1cm}

\nhp (III) The hypothesis $b \neq -a$ gives an $h_3 \in X$ such that $h_3(b)
\neq h_3(-a)$. An argument similar to that of (I) and (II), using the
assumptions $c \in \se {D^t} {\!X}(a,b)$ and $Z(a) \, \sub \, Z(b)$, shows that
$h_3$ can only take the following combination of values at $a,b,c$\,:
\vspace{-0.15cm}

\nhp \hspace{0.5cm} III.a\,: \;\; $h_3(a) = h_3(b) = h_3(c) \in \{ \pm 1 \}$\; or,
\vspace{-0.15cm}

\nhp \hspace{0.5cm} III.b\,: \;\; $h_3(b) = 0$ \, and \, $h_3(a) = h_3(c) \in 
\{\pm 1 \}$.
\vspace{-0.15cm}

\nhp With these data, a long, tedious, but straightforward checking of all
possible combinations of values of the characters $h_i \; (i = 1,2,3)$ at
the points $a,b,c,$ shows that $h = h_1h_2h_3$ has properties (*) and (**),
contradicting $c \in \se {D^t} {\!X}(a,b)$. This proves item (2).
\vspace{-0.1cm}

Next we show:
\vspace{-0.15cm}

\nhp (3) $Z(a)  \subset  Z(b) \; \Ra \; \se {D^t} {\!X}(a,b) = \{ a \}$.
\vspace{-0.15cm}

\nhp  \und{Proof of (3)}. $b \in \se {D^t} {\!X} (a,b)$ implies $Z(b) \, \sub \,
Z(a)$ (immediate verification); hence $b \not\in \se {D^t} {\!X} (a,b)$. Since
$Z(a)  \subset  Z(b)$ implies $b \neq -a$, items (1) and (2) give the
conclusion.
\vspace{-0.1cm}

The assertions (1) through (3) yield at once:
\vspace{-0.15cm}

\nhp (4) If $b \neq -a$, then $\; \se {D^t} {\!X}(a,b) = \{ a, b \} \;\; \Lra \;\;
Z(a) = Z(b) \;\; \Lra \;\; a^2 = b^2$.
\vspace{-0.15cm}

\nhp (5) $c \in \se {D^t} {\!X} (a,-a) \;\; \Lra \;\; c = a^2c  \;\; \Lra \;\; c =
a^2x$ for some $x \in G$.
\vspace{-0.15cm}

\nhp  \und{Proof of (5)}. The last equivalence is obvious: $c = a^2x$ implies
$a^2c = a^2(a^2x) = a^2x = c$. As\vspace{0.05cm} for the first, we have: 
\vspace{-0.15cm}

\nhp ($\Leftarrow$) Let $h \in X$. Obviously, $h(a) = -h(-a)$. The equality 
$c = a^2c$ implies $Z(a) \, \sub \, Z(c)$; hence $h(c) \neq 0$ implies  
$h(a) \neq 0$, and $h(c)$ equals either $h(a)$ or $h(-a)$, proving 
$c \in \se {D^t} {\!X} (a,-a)$.
\vspace{-0.15cm}

\nhp ($\Ra$) For $h \in X$, $c \in \se {D^t} {\!X} (a,-a)$ and $h(c) \neq 0$ imply
$h(c) = h(a)$ or $h(c) = -h(a)$; hence $h(a) \neq 0$. This shows that $Z(a) \,
\sub \, Z(c)$, which implies $a^2c^2 = c^2$ (cf. \ref{zerosetsElts}); scaling by 
\y c gives $c = a^2c$.
\vspace{-0.15cm}

This completes the proof of statement (A).
\vspace{-0.1cm}

Next we deal with the identity [$D$]. We shall, in fact, prove assertion (2) 
of Theorem \ref{fanRepresInterdefinable}:
\vspace{-0.15cm}

\nhp (B) If the ternary relation $D^t$ \und{is defined} as in [$D^t$] of
\ref{fanrepr} and the equivalence
\vspace{-0.15cm}

\nhp (\dag) \;\;\; $c \in D(a,b) \;\, \Lra \;\, c \in D^t(c^2a,\, c^2b)$
\vspace{-0.15cm}

\nhp holds for all $a,b,c \in G$, then $D$ verifies the equality [$D$]\,\footnote{\,The 
equivalence (\dag) is readily checked to hold for the relations $\se D X$ and $\se {D^t} {\!X}$, 
using their definitions [R] and [TR] in \ref{IdenticalReprRelats}.}. We write Id for Id(\y G).
\vspace{-0.15cm}

\nhp (6) $a \cdot \mbox{$\mathrm {Id}$} \; \sub \, D(a,b)$.
\vspace{-0.15cm}

\nhp  \und{Proof of (6)}. Let $x \in G$. By (\dag) it suffices to show:
\vspace{-0.15cm}

\nhp (\dag\dag) $\;\; ax^2 \in D^t(ax^2,\, a^2x^2b)$.
\vspace{-0.15cm}

\nhp Since $Z(ax^2) \, \sub \, Z(a^2x^2b)$, in case $a^2x^2b \neq -ax^2$, the
first and third clauses of [$D^t$] give (\dag\dag), and in case $a^2x^2b =
-ax^2$ the last clause in [$D^t$] proves (\dag\dag).
\vspace{-0.15cm}

\noi A similar argument gives $b \cdot \mbox{$\mathrm {Id}$} \; \sub \, D(a,b)$.
\vspace{-0.15cm}

\nhp (7) $ca = -cb \; \w \; c = a^2c \;\; \Ra \,\; c \in D(a,b)$.
\vspace{-0.15cm}

\nhp  \und{Proof of (7)}. Since  $c^2a = -c^2b$ and $c = a^2c = (c^2a^2)c$, the
last clause\vspace{0.05cm} in [$D^t$] yields \lbr 
$c \in D^t(c^2a,\, c^2b)$, whence, by (\dag), $c \in D(a,b)$.
\vspace{-0.1cm}

Items (6) and (7) prove the inclusion $\supseteq$ in [$D$]. Conversely, assuming
$c \in D(a,b)$, we have \lbr
$c \in D^t(c^2a,\, c^2b)$, by (\dag). An analysis according to the inclusions of
the zero-sets of $c^2a$ and $c^2b$ gives:
\vspace{-0.15cm}

\nhp (8) If $Z(c^2a) \, \sub \, Z(c^2b)$ and $c^2a \neq -c^2b$, then $c \in
D^t(c^2a,\, c^2b) \, \sub \, \{ c^2a,\, c^2b \}$, implying $c \in a \cdot
\textrm{Id} \, \cup \, b \cdot \textrm{Id}$\,.
\vspace{-0.15cm}

\nhp (9) If $c^2a = -c^2b$, scaling by \y c gives $ca = -cb$ and (by the last
clause in [$D^t$]) $c = c^2a^2x$ for some $x \in G$; this proves $a^2c =
a^2(c^2a^2x) = c^2a^2x =c$, as required. \hfl $\Box$


\nhp {\bf Proof of Theorem \ref{fanRepresInterdefinable}.} Item \ref{fanRepresInterdefinable}\,(2) 
has just been proved (item (B), proof of \ref{fanrepr}).
\vspace{-0.15cm}

\nhp  \und{Proof of \ref{fanRepresInterdefinable}\,(1)}. Assume $G$ as in the\vspace{0.05cm} statement, 
the ternary relation $D$ defined by clause [$D$] in \ref{fanrepr}, and $D^t$ given by:
\vspace{-0.15cm}

\nhp \centerline{$c \in D^t(a,b) \Leftrightarrow c \in D(a,b) \wedge -a \in
D(-c,b) \wedge -b \in D(a, -c).$} 

\vspace{-0.15cm}

\nhp The right-hand side of this equivalence amounts to:
\vspace{-0.15cm}

\nhp (I) \;\;\;\;\; $c \in \, a \cdot \textrm{Id}(G) \, \cup \, b \cdot
\textrm{Id}(G) \, \cup \, \{ x \in G \, \vert \, xa = -xb \; \w \, x = a^2x \}$.
\vspace{-0.15cm}

\nhp (II) \;\; $-a \in -c \cdot \textrm{Id}(G) \, \cup \, b \cdot \textrm{Id}(G)
\, \cup \, \{ x \in G \, \vert \, xc = xb \; \w \, x = c^2x  \}$.
\vspace{-0.15cm}

\nhp (III) \, $-b \in a \cdot \textrm{Id}(G) \, \cup \, -c \cdot \textrm{Id}(G)
\, \cup \, \{ x \in G \, \vert \, xa = xc \; \w \, x = a^2x  \}$.
\vspace{-0.15cm}

\nhp As above we write Id for Id(\y G). Remark that
\vspace{-0.15cm}

\nhp (*)\;\; $x \in y \cdot \textrm{Id} \;\, \Lra \;\, x = yx^2$, \;\;\;\;
and\;\;\;\; (**)\;\; $xy = xz \;\; \Ra \;\; xy^2 = xz^2\, (= xyz)$.
\vspace{-0.15cm}

We argue by cases, according to the various clauses in [$D^t$].
\vspace{-0.15cm}

\nhp  (1) $Z(a) \subset  Z(b) \;\, \Ra \;\, c = a$.
\vspace{-0.15cm}

\nhp The clauses $-a \in b \cdot \textrm{Id}$ and $-a = c^2(-a) = b^2(-a)$ (see
(**)) in (II) imply $Z(b) \, \sub \, Z(a)$, and hence are excluded; thus, (II)
reduces to $a \in c \cdot \textrm{Id}$. The following cases arise from (I) and (II):

\vspace{-0.15cm}

\nhp (1.i) $\;c \in a \cdot \textrm{Id}$  and  $a \in c \cdot \textrm{Id}$.
\vspace{-0.15cm}

\nhp By (*), $c = ac^2$ and $a =ca^2$. Hence, $c = ac^2 = (ca^2)c^2 = ca^2 = a$.
\vspace{-0.15cm}

\nhp (1.ii) $\;c \in b \cdot \textrm{Id}$  and  $a \in c \cdot \textrm{Id}$.
\vspace{-0.15cm}

\nhp Then, $c = bc^2$ and $a =ca^2$, implying $a =ca^2 = a^2c^2b$\,; it follows
that $Z(b) \, \sub \, Z(a)$, contrary to the assumption in (1).
\vspace{-0.15cm}

\nhp (1.iii) $\;ca = -cb \; \w \; c = a^2c = b^2c \; \w \; a \in c \cdot
\textrm{Id}$.
\vspace{-0.15cm}

\nhp The middle equality implies $Z(b) \, \sub \, Z(c)$ and the last $Z(c) \,
\sub \, Z(a)$. Hence $Z(b) \, \sub \, Z(a)$, and this case is also excluded.
\vspace{-0.15cm}

\nhp  (2) $Z(b) \subset  Z(a) \;\, \Ra \;\, c = b$.
\vspace{-0.15cm}

\nhp Same argument as in (1) interchanging \y a and \y b.
\vspace{-0.15cm}

\nhp  (3) $Z(a) = Z(b) \, \w \; b \neq -a \;\, \Ra \;\, c \in \{ a, b \}$.
\vspace{-0.15cm}

\nhp The first assumption gives $a^2 = b^2$. Each of the clauses $-a \in b \cdot
\textrm{Id}$ in (II) and $-b \in a \cdot \textrm{Id}$ in (III) yield $-a = ba^2
= b$ and hence are excluded. From (I)\,--\,(III) the following cases arise:
\vspace{-0.15cm}

\nhp (3.i) $\;\;c \in a \cdot \textrm{Id}$  and  $a \in c \cdot \textrm{Id}$.
\vspace{-0.15cm}

\nhp We have $c = ac^2$ and $a =ca^2$; as in (1.i) we get $c = a$.
\vspace{-0.15cm}

\nhp (3.ii) $\;\,c \in a \cdot \textrm{Id}$, $ac = ab$ and $a = c^2a$.
\vspace{-0.15cm}

\nhp The first term gives $c = ac^2$; hence $a = c$. The cases
\vspace{-0.15cm}

\nhp (3.iii) $\;c \in b \cdot \textrm{Id}$  and  $-b \in -c \cdot \textrm{Id}$,\, and
\vspace{-0.15cm}

\nhp (3.iv) $\;c \in b \cdot \textrm{Id}$, $ab = ac$  and  $b = c^2b$,
\vspace{-0.15cm}

\nhp are similar to (3.i) and (3.ii) ---with \y b replacing \y a---, and yield $c = b$.
\vspace{-0.15cm}

\nhp (3.v) $\;\; ca = -cb, \, c = a^2c = b^2c$  and  $-a \in -c \cdot
\textrm{Id}$.\vspace{-0.15cm}

\nhp As in (3.ii) this gives $c = a^2c = a$ (see (*)).
\vspace{-0.15cm}

\nhp (3.vi) $\; ca = -cb, \, c = a^2c = b^2c$  and  $-b \in -c \cdot
\textrm{Id}$.
\vspace{-0.15cm}

\nhp As in (3.v) we obtain $c = b$.
\vspace{-0.15cm}

\nhp (3.vii) $ca = -cb, \ c = a^2c = b^2c, \ ac = ab, \ a = c^2a, \ ab = bc$ 
and  $b = c^2b$ \, (the last disjunct from (I), (II) and (III)).
\vspace{-0.15cm}

\nhp We have $ac = -bc, \, ac = ab$ and $ab = bc$; hence $bc = -bc$. Scaling by
\y b, $b^2c = -b^2c$, whence $c = -c$. It follows that $c = 0$, 
which clearly implies $a = b =0$, i.e., $a,b,c$ are all 0.
\vspace{-0.15cm}
 
\nhp (4) $b = -a \;\, \Ra \;\, c = a^2c = b^2c$.
\vspace{-0.15cm}

\nhp Each disjunct in (I) implies $c = a^2c$. The third disjunct contains this
condition. If $c \in a \cdot \textrm{Id}$, then $c = ac^2$ (by (*)); hence $a^2c
= a^2(ac^2) = ac^2 = c$. Likewise, $c \in b \cdot \textrm{Id}$, implies $c =
b^2c = a^2c$. \hfl $\Box$

\vspace{-0.1cm}

\nhp {\bf Proof of Theorem \ref{fanRS}.} We check that the axioms for real semigroups, cf.
\ref{RS}, hold in $(G,D)$. By Theorem \ref{RSaxioms} we need only prove the associativity
axiom [RS3] or, instead, the equivalent statement [RS3$'$], see Proposition \ref{axRS3'}. 
To abridge we shall write:
 \vspace{-0.15cm}

\noi (*) for $\fa\,a,b,c,d\:(D^t(a,b) \cap D^t(c,d) \neq\, \0$, the hypothesis of [RS3$'$], and 
 \vspace{-0.15cm}

\noi (**) for $D^t(a,-c) \cap D^t(-b,d) \neq\, \0$, its conclusion. 
 \vspace{-0.15cm}

\noi The characterization of transversal representation in Theorem \ref{fanrepr} will be of 
constant use and will be referred to as ``$[D^t]$".

\vspace{-0.1cm}
We consider three cases and, for each of them, several subcases according to the 
mutual inclusions of the zero-sets of $a, b, c, d.$
\vspace{-0.15cm}

\nhp \und{Case A}.\; $a = -b$.
\vspace{-0.15cm}

\nhp (A.i)\, $Z(c) \subset Z(a)$.
\vspace{-0.15cm}

\nhp By the last clause in $[D^t]$ and (*) there is\, $x \in G$\, so that\,
$a^2\,x \in D^t(c,d)$. If\, $Z(c) \subset Z(d)$, we would have\, $D^t(c,d) =
\{c\}$, whence\, $a^2\,x = c$, which implies\, $Z(a)\, \sub\, Z(c)$, contrary to
(A.i). Hence,\, $Z(d)\, \sub\, Z(c) \subset Z(a)$, and this yields\, $a^2\,x
\neq c, d$, implying\, $D^t(c,d) \not\subseteq \{c, d\}$. By $[D^t]$ we then
have\, $d = -c$, and it follows that\, $D^t(a,-c) \cap D^t(-b,d) = D^t(a,-c) =
\{-c\} \neq\, \0.$ 
\vspace{-0.15cm}

\nhp (A.ii)\, $Z(a) \subseteq Z(c)$.
\vspace{-0.15cm}

\nhp From $[D^t]$ we have\, $a \in D^t(a,-c)$, and show that\, $a \in D^t(-b,d)
= D^t(a,d)$. Otherwise, by $[D^t]$ again, we must have\, $Z(d)\, \subset Z(a)\,
\sub\, Z(c)$, and assumption (*) gives an\, $x \in G$ such that\, $a^2\,x
\in D^t(c,d) = \{d\}$, implying\, $Z(a)\, \sub\, Z(d)$, a contradiction.
\vspace{-0.1cm}

\nhp \und{Case B}.\; $c = -d$. Argument similar to that of Case A.
\vspace{-0.1cm}

\nhp \und{Case C}.\; $a \neq -b$\, and\, $c \neq -d$.
\vspace{-0.15cm}

\nhp The first three clauses of $[D^t]$ show that\, $D^t(a,b)\, \sub\,
\{a,b\}$\, and\, $D^t(c,d)\, \sub\, \{c,d\}$. We consider the following
subcases:
\vspace{-0.15cm}

\nhp (C.i)\; $Z(a) \subset Z(b)$\, and\, $Z(c) \subset Z(d)$.
\vspace{-0.15cm}

\nhp In this case $[D^t]$ and (*) imply $a = c$, and (**) reduces to\,
$D^t(a,-a)\, \cap\, D^t(-b,d) \neq\, \0$. Since\, $Z(a) = Z(c) \subset Z(b)\, \cap\, Z(d)$, 
we get\, $b = a^2\,b$\, and\, $d = a^2\,d$ (\ref{zerosetsElts}). 
If $b \neq\, d$, the first three clauses of $[D^t]$ show that one
of\, $-b$\, or\, \y d\, is in\, $D^t(-b,d)$, and the preceding equalities give\,
$a^2\,x \in D^t(-b,d)$\, for some\, $x \in G$. By the last clause of $[D^t]$
this also holds if\, $b = d$, proving\, $D^t(a,-a)\, \cap\, D^t(-b,d) \neq\, \0$, 
as required.
\vspace{-0.15cm}

\nhp (C.ii)\; $Z(a) \subset Z(b)$\, and\, $Z(d) \subseteq Z(c)$. 
\vspace{-0.15cm}

\nhp From (*) we have\, $a \in D^t(c,d)$; then $Z(d) \subseteq Z(c)$ implies\, 
$a = c$\, or\, $a = d$. In either case,\, $Z(d) \subset Z(b)$, whence\, 
$D^t(-b,d) = \{d\}$, and we are reduced to prove\, $d \in D^t(a,-c)$.
\vspace{-0.1cm}

In case\, $a = c$\, we must show that\, $a^2\,x = c^2\,x = d$\, for some\, 
$x \in G$. If\, $Z(d) = Z(a)$\, this holds with\, $x = d$\, by \ref{zerosetsElts}.
If\, $Z(d) \subset Z(a) = Z(c)$, (*) gives\, $a \in D^t(c,d) = \{d\}$, leading 
to\, $a = c = d$, a contradiction.
\vspace{-0.1cm}

Finally, in case\, $a = d$, since\, $D^t(-b,d) = \{d\}$, we are reduced to
prove\, $d \in D^t(d,-c)$. Since we may assume\, $d \neq\, c$, this follows
from\, $Z(d) \subseteq Z(c)$, using $[D^t]$.
\vspace{-0.15cm}

\nhp (C.iii)\; $Z(a) = Z(b)$\, and\, $Z(c) \subset Z(d)$.
\vspace{-0.15cm}

\nhp Assumptions (*) and (3) imply\, $c = a$\, or\, $c = b$, and we have\, $Z(a)
= Z(b) = Z(c)$. Then, the conclusion (**) boils down to\, $-b \in D^t(a,-c)$.
If\, $c = b$\, this follows from\, $Z(a) = Z(-b)$\, by the third clause of
$[D^t]$. If\, $c = a$, then\, $Z(a) = Z(-b)$\, implies\, $-b = a^2(-b)$, and the
conclusion holds as well. 
\vspace{-0.15cm}

\nhp (C.iv)\; $Z(a) = Z(b)$\, and\, $Z(d) \subseteq Z(c)$.
\vspace{-0.15cm}

\nhp If\, $Z(d) \subset Z(c)$, $[D^t]$ and the assumptions (*) and (3) give\, $d
\in D^t(a,b) = \{a,b\}$. If\, $d = a$, the desired conclusion boils down to\, $a
\in D^t(-b,a)$, as\, $D^t(a,-c) = \{a\}$. Since\, $Z(a) = Z(-b)$, this holds by
the last two clauses of $[D^t]$. If\, $d = b$ (and\, $a \neq\, b)$, conclusion
(**) reduces to\, $a \in D^t(-b,b)$, since\, $D^t(a,-c) = \{a\}$. But\, $Z(a) =
Z(b)$\, implies\, $a = b^2\,a \in D^t(-b,b)$.
\vspace{-0.1cm}

If\, $Z(d) = Z(c)$, assumptions (*) and (3) give\, $\{a,b\}\, \cap\, \{c,d\}\,
\neq\, \0$. It follows that\, $Z(a) = Z(b) = Z(c) = Z(d)$, and then\,
$\{a,-c\}\, \sub\, D^t(a,-c)\,,\, \{-b,d\,\}\, \sub\, D^t(-b,d)$. If\, $a = c$,
then\, $Z(a) =  Z(c) = Z(d)$\, entails\, $d \in D^t(a,-c)$, and (**) follows.
If\, $d = b$, the same argument shows\, $a \in D^t(-b,d)$. If\, $a \neq\, c$\,
and\, $d \neq\, b$, then\, $a = d$\, and\, $b = c$; this obviously implies\,
$\{a,-c\}\, \cap\, \{-b,d\,\} \neq\, \0$, whence\, $D^t(a,-c)\, \cap\, D^t(-b,d)
\neq\, \0$.
\vspace{-0.15cm}

\nhp (C.v)\; $Z(b) \subset Z(a)$\, and\, $Z(c) \subset Z(d)$.
\vspace{-0.15cm}

\nhp In this case we have\, $D^t(a,b) = \{b\}$\, and\, $D^t(c,d) = \{c\}$,
whence\, $b = c$, by (*). Therefore\, $Z(c) \subset Z(a)$,\, $Z(b) \subset
Z(d)$, which yields\, $D^t(a,-c) = \{-c\}$,\, $D^t(-b,d) = \{-b\}$, and hence (**). 

\vspace{-0.15cm}

\nhp (C.vi)\; $Z(b) \subset Z(a)$\, and\, $Z(d) \subseteq Z(c)$.
\vspace{-0.15cm}

\nhp The first clause of $[D^t]$ gives\, $D^t(a,b) = \{b\}$, and hence\, $b =
c$\, or\, $b = d$, by (*) and (3). In the latter case we must show that\,
$b^2\,x \in D^t(a,-c)$\, for some\, $x \in G$. If\, $a \in D^t(a,-c)$, since\,
$b^2\,a = a$\, (as\, $Z(b) \subset Z(a)$), it suffices to take\, $x = a$. If\,
$a \not\in D^t(a,-c)$, then\, $Z(c) \subseteq Z(a)$, and the second and fourth
clauses of $[D^t]$ yield\, $-c \in D^t(a,-c)$; since\, $b^2(-c) = -c$\, ($Z(b)
\subseteq Z(c)$), we can take\, $x = -c$.
\vspace{-0.1cm}

Finally, if\, $b = c$, then\, $D^t(a,-c) = \{-c\,\}$. If\, $-c \not\in D^t(-c,d)
= D^t(-b,d)$, then\, $Z(d) \subset Z(c) = Z(b) \subset Z(a)$, and (*) yields\,
$b = d$, a contradiction. Thus,\, $-c \in D^t(-b,d)$, verifying (**) and
completing the proof of Theorem \ref{fanRS}. \hfl $\Box$
\vspace{-0.35cm}

\bct \label{qfans} {\bf A digression on q-fans.}
\vspace{-0.15cm}

The results proved above use in a crucial way the auxiliary 
---but nonetheless important--- notion of a q-fan, introduced in
\ref{fans-various}\,(3). In Corollary \ref{equivfan} we proved that q-fans verifying 
Marshall's axiom [AX2] for ARSs are the same thing as fans. The following example
shows that this notion is genuinely weaker than that of a fan. 

\vspace{-0.35cm}

\bex \label{exqfan} A q-fan that is not a fan.
\vspace{-0.15cm}

\noi Let $C$ be the set of all non-decreasing functions $f: \N \lra \{0,1\}$, where
$0 < 1$; thus, $f \in C$\; iff\; $f = 0$\: or\: there is $n \in \N$ such that $f(m) = 0$ 
if $m < n$ and $f(m) = 1$ if $m \geq n$. Then, $-C = \{(-1)\cdot f \sth f \in C\}$ is the 
set of non-increasing maps $g: \N \lra \{0,-1\}$ (with $-1 < 0$). Let $T = C \cup -C$. 
Straightforward checking shows that, with pointwise defined product, $T$ is a ternary 
semigroup, having as distinguished elements the constant functions with values $1, 0$ and 
$-1$. In addition, $T$ verifies:
\vspace{-0.15cm}

\noi (i)\; Id$(T) = C$;
\vspace{-0.15cm}

\noi (ii) If $f,g \in T$, then $f \cdot g = f$\, or\, $f \cdot g = g$.
\vspace{-0.15cm}

\noi For $n \in \N$, let $\pi_n: T \lra \{0,1, -1\}$ stand for the projection onto the
$n$-th coordinate: $\pi_n(f) = f(n)\; (f \in T)$. Straightforward checking, using 
that the functions in $C$ are non-decreasing and those in $-C$ are non-increasing,
proves:
\vspace{-0.14cm}

\noi (\dag)\; For $k,n,m \in \N$,\; $k \leq n \leq m\; \Ra\; \pi_k \cdot \pi_n \cdot \pi_m = \pi_k$.
\vspace{-0.14cm}

\noi In other words, the projections $\pi_n\; (n \in \N)$ form a subset $P$ of $X_T$ 
closed under the product of any three of its members. It is also clear that $P$ separates
points in $T$. Hence, $(P,T)$ is a q-fan.
\vspace{-0.12cm}

However, $(P,T)$ \und{is not} a fan. To see this, consider the map $h:T \lra \{-1,0,1\}$
defined by: for $f \in T$, 
\vspace{-0.75cm}

\nhp $$
h(f)=
\left\{
\begin{array}{ll}
\;1 & \;\mbox{if}\;\; f\in\, C \setminus \{0\}\\
\;0& \;\mbox{if}\;\; f = 0\\
-1 & \;\mbox{if}\;\; f\in\, -C \setminus \{0\}.
\end{array}
\right.
$$

\vspace{-0.2cm}

\noi The reader can easily check that $h$ is a TS-character of $T$. However, $h \not\in P$;
for, if $n \in \N$, the map $f_n$ given by, $f_n(k) = 0$ for $k \leq n$ and $f_n(k) = 1$ for $k > n$, 
is in $C \setminus \{0\}$, whence $h(f_n) \neq 0$, while $\pi_n(f_n) = 0$. \hfl $\Box$
\eex
\ect
\vspace{-0.35cm}

The following result gives an interesting topological characterization of q-fans.
\vspace{-0.4cm}

\bth \label{qfandense} Let\, \y G be a ternary\vspace*{0.05cm} semigroup and let\, 
$X\, \sub\, \se X G =$ \HTS G\, be a non-empty set of TS-characters closed
under product of any three of its members. Then,
\vspace{-0.15cm}

\noi \hfl $(X,G)$\, is a q-fan\; $\Lra$\; \y X\, is dense for the constructible 
topology of\, $\se X G$\,. \hfl
\vspace{-0.15cm}

\noi In particular, if\, \y X\, is proconstructible, i.e., closed in the
constructible topology, then\, $X = \se X G$\,, and hence $(X,G)$ is a fan.
\eth

\vspace{-0.4cm}

\nhp {\bf Proof.} We shall write $a\, \se {\equiv} X\, b$ to mean $g(a) = g(b)$ for 
all $g \in X$; cf. \ref{propsTS-congr}\,(c). Hence, the clause that $X$ separates 
points in $G$ translates as\; $a\, \se {\equiv} X\, b\; \Lra\; a = b$.
\vspace{-0.15cm}

\noi ($\La$) Since $X$ is assumed to be closed under products of any three of its 
members, we need only show that $X$ separates points in $G$. Let $a \neq b$. By 
the separation theorem for TS, \cite{DP1}, Thm. 1.9, pp. 103--104, the set 
$\{h \in \se X G \sth h(a) \neq h(b)\}$ is non-empty. It is easily checked
that this set is open in the constructible topology of $\se X G$, and hence it 
intersects $X$, i.e., there is $h \in X$ such that $h(a) \neq h(b)$.
\vspace{-0.15cm}

\noi ($\Ra$) By definition the sets 
\vspace{-0.15cm}

\nhp (*) \hfl $U = U(\se a {1}) \cap \ldots \cap U(\se a {n}) \cap Z(\se b {1}) \cap
\ldots \cap Z(\se b {k})$, \hfl
\vspace{-0.15cm}

\nhp with\, $\se a {1},\ldots,\se a {n},\se b {1},\ldots,\se b {k} \in G$, form
a basis\vspace*{0.05cm} for the constructible topology of\, $\se X G$. 
We must show:\, $U \neq\, \0\;\, \Ra\; U\, \cap\, X \neq\, \0$.
\vspace{-0.15cm}

Assume otherwise, and fix $h \in U$. We first show
\vspace{-0.15cm}

\noi {\bf Claim A.} Let $U$ be a basic open containing $h$ such that $U\, \cap\, X = \0$, 
with minimal $n$ having the form (*) for some $\se b {1},\ldots,\se b {k} \in G\;\,
(k \geq 0)$. Then, $n > 1$.
\vspace{-0.15cm}

\noi {\bf Proof of Claim A.} \und{Case 1}. $n = 0$ and $k = 1$. Thus, $Z(\se b {1})\, \cap\, X = 
\emptyset$; this means\, $\se {b^2} {1}\: \se {\equiv} X 1$, hence $\se {b^2} {1} = 1$. Since 
$h \in \HTS G$, we get $h(\se {b^2} 1) = 1$, contradicting $h \in Z(\se b {1})$.
\vspace{-0.15cm}

\nhp \und{Case 2}. $n = 0$ and $k \geq 2$. Let $k$ be the smallest integer such that 
there are elements $\se b {1},\ldots,\se b {k}$ so that 
$h \in \bigcap_{j=1}^k Z(\se b {j}) $ and $\bigcap_{j=1}^k Z(\se b {j})\, \cap\, X = \emptyset$. 
Then, $\bigcap_{j=1}^{k - 1} Z(\se b {j})\, \cap\, X \neq \emptyset$\, 
and $\bigcap_{j=2}^{k} Z(\se b {j})\, \cap\, X \neq\emptyset$; let\, $\se h {1},\se h {2}$ 
belong, respectively, to these sets. Since $X$ is closed under products of any three elements,
$\se {h^2} {1}\se h {2} \in X$; clearly, $\se {h^2} {1}\se h {2} \in \bigcap_{j=1}^k Z(\se b {j})$,
contradicting the choice of $k$.
\vspace{-0.15cm}

\nhp \und{Case 3}. $n = 1$ and $k = 0$. Then $U(\se a {1}) \cap X = \emptyset$; this implies 
$\se a {1}\, \se {\equiv} X \!-\se {a^2} {1}$ and hence $\se a {1} = -\se {a^2} {1}$, contradicting 
$h(\se a {1}) = 1$.
\vspace{-0.15cm}

\nhp \und{Case 4}. $n = 1$ and $k = 1$. Then,
\vspace{-0.1cm}

\noi (\dag)\; $U(\se a {1}) \cap Z(\se b {1}) \cap X = \emptyset$.
\vspace{-0.15cm}

\noi Thus, $h(\se a {1}) = 1, h(\se b {1}) = 0$; this yields $\se {a^2} {1} \neq \se {a^2} {1}\se {b^2} {1}$.
Since $X$ separates points, there is $\se h {1} \in X$, so that $\se h {1}(\se b {1}) = 0$
and $\se h {1}(\se a {1}) \neq 0$. From (\dag) we get $\se h {1}(\se a {1}) = -1$. On the 
other hand, the minimality assumption implies $U(\se a {1}) \cap X \neq \0$, i.e., there is 
$\se h {2} \in X$ such that $\se h {2}(\se a {1}) = 1$. Then,\, 
$\se {h^2} {1} \se h {2} \in X,\; \se {h^2} {1} \se h {2}(\se a {1}) = 1$\, and\, 
$\se {h^2} {1} \se h {2}(\se b {1}) = 0$, i.e., $\se {h^2} {1} \se h {2} \in 
U(\se a {1}) \cap Z(\se b {1}) \cap X = \emptyset$, contrary to (\dag).

\vspace{-0.15cm}

\nhp \und{Case 5}. $n = 1$ and $k \geq 2$. Let $k$ be the least integer such that
there are $\se a {1},\se b {1},\ldots,\se b {k} \in G$ with
$U(\se a {1}) \cap \bigcap_{j=1}^k Z(\se b {j}) \cap X = \emptyset$. 
Then, there are $\se h {1},\se h {2} \in X$ so that $\se h {1}(\se a {1}) = \se h {2}(\se a {1}) = 1$, 
$\se h {1}(\se b {1}) = \ldots = \se h {1}(\se b {k-1}) = 0$ and $\se h {2}(\se b {2}) = \ldots = 
\se h {2}(\se b {k}) = 0$. Hence, $\se {h^2} {1}\se h {2} \in X \cap U(\se a {1}) \cap 
\bigcap_{j=1}^k Z(\se b {j})$, contradiction. \hfl $\Box$
\vspace{-0.1cm}

Next, observe that setting $c = \prod_{i=1}^{n} \se {a^2} {i}$ we have 
$U(\se {c\,a} {i})\, \sub\, U(\se {a} {i})$ for all $i = 1, \ldots, n$, and 
$\bigcap_{i=1}^n U(\se {c\,a} {i}) = \bigcap_{i=1}^n U(\se {a} {i})$. Hence,
\vspace{-0.3cm}

\noi (**) \hspace*{0.9cm} $h \in\; \displaystyle{\bigcap_{i=1}^n U(\se {c\,a} {i}) \cap 
\bigcap_{j=1}^k Z(\se b {j})} \;\;$ and \;\;
$\displaystyle{\bigcap_{i=1}^n U(\se {c\,a} {i})  \cap \bigcap_{j=1}^k Z(\se b {j}) 
 \cap X = \emptyset.}$ \hfl
\vspace{-0.2cm}

Choosing $n$ minimal so that $(**)$ holds for some $\se a i, \se b j\; (k \geq 0)$, for
each index $i \in \{1, \ldots, n\}$ there is $\se g {i} \in \se X G$ such that
\vspace{-0.3cm}

\noi (\dag\dag) \hspace*{1cm} $\se g {i} \in \displaystyle{\bigcap_{\ell=1, \ell \neq i}^n 
U(\se {c\,a} {\ell})\, \cap\, \bigcap_{j=1}^k Z(\se b {j})\, \cap\, X}$. 
\vspace{-0.2cm}

\noi Since $n > 1$ (Claim A), we get $\se g {i}(c) = 1$, whence $\se g {i}(\se a {i}) \neq 0$;
from (\dag\dag) comes
\vspace{-0.1cm}

\noi (\dag\dag\dag) \hspace*{0.8cm} $\se g {i}(\se a {i}) = -1$\;\; and\;\; 
$\se g {i}(\se a {\ell}) = 1$ for $i \neq \ell$.
\vspace{-0.1cm}

Observe next that $h \in U(\se a {1} \!\cdot \se a {2}\! \cdot \ldots \cdot \se a {n})\, \cap\, 
\bigcap_{j=1}^k Z(\se b {j})$. From Cases 3 -- 5 in the proof of Claim A, we get 
$U(\se a {1}\!\cdot \se a {2}\! \cdot \ldots \cdot \se a {n})$ $\cap\, \bigcap_{j=1}^k Z(\se b {j})\, 
\cap\, X \neq \emptyset$. Let $\se g {n+1}$ belong to this intersection; in particular, 
$\se g {n+1}(\se a {1} \cdot \ldots \!\cdot \se a {n}) = 1$. From (**) follows\, 
$\se g {n+1}(\se a {i}) = -1$\, for some index $i$. Then, the set
$\{i \in \{1,\ldots,n\}\, \vert\,\se g {n+1}(\se a {i}) = -1\}$ has even cardinality $> 0$; 
let $\{\se i {1},\ldots,\se i {2m}\}$ be an enumeration of it, and let $g := \se g {n+1} 
\cdot \prod_{j=1}^{2m} \se g {i_{j}}$. The assumption that $X$ is closed under products
of any three elements implies that it is closed under under products of any odd number of
its elements; therefore, $g \in X$. Clearly, $g \in \bigcap_{j=1}^k Z(\se b {j})$. Claim B
below shows that $g(\se a i) = 1$ for all $i = 1,\ldots, n$, which yields 
$g \in \bigcap_{i=1}^n U(\se {c\,a} {i})  \cap \bigcap_{j=1}^k Z(\se b {j}) \cap X$,
contradicting (\dag\dag), and hence will complete the proof of Theorem \ref{qfandense}.

\noi {\bf Claim B.}\; With notation as above, $g(\se a i) = 1$\: for all $i = 1,\ldots, n$.
\vspace{-0.1cm}

\noi {\bf Proof of Claim B.} Let  $i \in \{1,\ldots,n\}$. If $i \not\in \{\se i {1},\ldots,\se i {2m}\}$, 
then $\se g {i_{j}}(\se a {i}) = 1$ for every index $j \in \{1,\ldots,2m\}$, and 
$\se g {n+1}(\se a {i}) = 1$, implying $g(\se a {i}) = 1$. If $i = \se i {j}$ 
for some index $j \in \{1,\ldots,2m\}$, then $\se g {i_{j}}(\se a {i_{j}}) = -1$,\; 
$\se g {i_{\ell}}(\se a {i_{j}}) = 1$ for $\ell \neq j$,\, and\, $\se g {n+1}(\se a {i_{j}}) = -1$,
whence $g(\se a {i}) = 1$. \hfl $\Box$
\vspace{-0.1cm}

\noi {\bf Remark.} For a result akin to Theorem \ref{qfandense} in the context of reduced
special groups, see \cite{DM1}, Prop. 3.16, p. 55. \hfl $\Box$

\vspace{-0.5cm}

\section{Examples} \label{fan-examp}

\vspace{-0.3cm}
With the aim of illustrating the notions introduced above, we present in this section some examples of 
(alas, finite) fans based on ternary semigroups with up to three generators. For each example we shall 
draw both the root-system of an ARS-fan ordered under specialization and the representation partial 
order of its dual real semigroup.

In order to determine the representation partial order in the examples below, we will need the following 
supplement to Theorem \ref{reprpo-RS}, valid in the case of fans. 
\vspace{-0.4cm}

\ble \label{reprpo-fans} Let \y G be a RS-fan, and let\, $\leq$ denote its representation partial order 
{\em (\ref{reprord})}. For a non-invertible element\, $x \in G$, and a unit\, $w \in G$, we have:
\vspace{-0.15cm}

\nhp $(1)$ If\, $x \leq w$, then\, $w = -1$. \hmm\hmm $(2)$ If\, $w \leq x$, then\, $w = 1$.
\vspace{-0.15cm}

\nhp In other words,\, \y x\, and\, \y w\, are\, $\leq$-incomparable, unless\, $w \in \{ \pm 1 \}$.
\vspace{-0.15cm}

\nhp $(3)$ If\, $v \in G$\, is invertible,\, $v \neq w$\, and\, $v, w \not\in \{ \pm 1 \}$, then\, \y v\, 
and\, \y w\, are\, $\leq$-incomparable.
\ele  

\vspace{-0.35cm}
\nhp {\bf Proof.} (2) By assumption, $w^2 = 1$, and by [RS6], $w \in \D F(1,x)$ implies
$w \in \Dt F(w^2, w^2x) = \Dt F(1,x)$. By Theorem \ref{fanrepr}, $w = 1$.
\vspace{-0.1cm}

\noi (1) follows at once from (2).
\vspace{-0.15cm}

\noi (3) By Corollary \ref{unitsRS-fans}, $G^{\times}$ with representation induced from $G$
is a RSG-fan (a fan in the category of reduced special groups). Hence, if $w \neq -1$, 
$v \in \D {G^{\times}}(1,w)$ implies $v = 1$ or $v = w$ (cf. [RSG-fan], Introduction), from
which (3) follows. \hfl $\Box$

\vspace{-0.3cm}

\bct {\bf The examples.} \label{ex-fans} 

\ect
\vspace{-0.4cm}

\nhp {\bf Example \ref{ex-fans}.\,A.}\; \und{Ternary semigroups on one generator}.

\nhp Call\, \y x\, the generator. We treat first the case where there are no additional relations 
(``free" case). The corresponding TS is:
\vspace{-0.15cm}

\nhp \hfl $\se F 1 = \{1, 0, -1, x, -x, x^2, -x^2\}.$ \hfl
\vspace{-0.15cm}

\nhp The necessary condition [Z] is trivially verified. Characters are determined by their value on\, \y x, and 
any value 1, 0 and $-1$ is possible; hence the dual ARS,\, $\se X {F_1}$\,, consists of three characters given by:\, 
$\se h 1(x) = 0,\, \se h 2(x) = 1,\, \se h 3(x) = -1$. Clearly,\, $\se h 1 = \se {h^2} 1 \cdot \se h i$, whence\, 
$\se h i \spez \se h 1$, for $i = 2, 3$\, (Lemma \ref{char-specializ}). So we get the specialization root-system 
below left. 

\nhp \parbox[t]{200pt} {

\vspace{0.05cm}

$$\hspace{-2cm}\xymatrix@R=15pt@C=20pt{
& & & & *{\txt{$h_1$ \\ $\bullet$}} \ar@{-}@<1.04ex>[dl] \ar@{-}@<-1ex>[dr] & \\
& & & *+<1.7ex,0.4ex>{h_2\:{\bullet}} & & *{{\bullet}\;h_3}} \\
$$ 

\vspace{0.85cm}

\nhp\hspace{1cm} Specialization root-system of $\se X {F_1}$ }
\hfl \parbox[t]{180pt} {

\vspace{-0.55cm}

$$\hspace{-1cm}\xymatrix@R=12pt@C=16pt{
 & & \node \rlab{\mbox{$\raisebox{4pt}{$-1$}$}} \edge{d} & \\
 & & \node \edge{d} \rlab{\mbox{$\raisebox{4pt}{$\;\;-x^2$}$}} & \\
   & \node \llab{$-x\;\;$} \edge{ur} \edge{dr}
   & \node \edge{d} \rlab{$\hspace{-0.17cm} 0$} 
   & \node \rlab{$\hspace{-0.17cm} x$} \edge{ul} \edge{dl} & \\
 & & \node \rlab{\mbox{$\raisebox{3pt}{$x^2$}$}} \edge{d} & \\
 & & \node \rlab{$\hspace{-0.1cm}1$} &
 } 
$$

\vspace{0.05cm}

\nhp Representation partial order of $\se F 1$ }

\nhp \centerline{\hspace{0.2cm} Figure 1} 

\nhp The representation partial order of the real semigroup \se F 1 ---illustrated in Figure 1, right--- 
follows straightforwardly from Theorem \ref{reprpo-RS}\,(2)--(4).

\nhp {\bf Remark.} Barring the case where the generator\, \y x\, becomes invertible (i.e.,\, $x^2 = 1$, 
which gives a four element RSG-fan with an added 0), the only possible additional relation is\, $x^2 = x$,
which eliminates the character\, \se h 3. Thus, we get the following diagrams for the specialization order 
(left) and the representation order (right):

\vspace{-0.2cm}

\hfl \parbox[t]{50pt} {

\vspace{-0.2cm}

$$\hspace{-4cm}\xymatrix@R=22pt{
  & \node \rlab{$h_1$} & \\
  & \node \rlab{$h_2$} \edge{u} & \\
}
$$
}  
\hfl \parbox[t]{100pt} {

\vspace{-0.7cm}

$$\hspace{-6cm}\xymatrix@R=12pt{
 & \node \rlab{$-1$} & \\
 & \node \rlab{$-x$} \edge{u} & \\
 & \node \rlab{$0$} \edge{u} & \\
 & \node \rlab{$x$} \edge{u} & \\
 & \node \rlab{$1$} \edge{u} & \\
} 
$$
 
\vspace{-0.2cm}
 
\nhp \hspace*{-5.3cm} Figure 2 
}

A more interesting example is:

\nhp {\bf Example \ref{ex-fans}.\,B.}\; \und{Ternary semigroups on three generators}.
\vspace{-0.15cm}

\nhp Generators:\, \y x, \y y, \y z. Condition [Z] gives raise to the following possible relations:
\vspace{-0.15cm}

\nhp 1.\; $x^2 = y^2 =\, z^2$. (\,[Z] is automatically verified in this case.)
\vspace{-0.15cm}

\nhp 2.\; $x^2 = y^2 \neq\, z^2$\, and\, $x^2z^2 =\, y^2z^2 \in \{x^2, z^2\}$.
\vspace{-0.15cm}

\nhp The two identities obtained from the last clause give rise to non-isomorphic cases, and, 
upon permutation, all cases where two of the three generators have equal squares (i.e., 
equal zero-sets) are isomorphic to these.
\vspace{-0.15cm}

\nhp 3.\; $x^2, y^2, z^2$\, are pairwise different, and\, $x^2y^2 \in \{x^2, y^2\},\;
x^2z^2 \in \{x^2, z^2\},\; y^2z^2 \in \{y^2, z^2\}$. 
\vspace{-0.15cm}

\nhp A case-by-case analysis of all eight combinations of these values shows that, up to 
isomorphism by permutation, the only surviving case is\, $x^2y^2 = x^2z^2 = x^2$\, and\, $y^2z^2 = y^2$.
\vspace{-0.15cm}

As an illustration we analyze the following configuration:
\vspace{-0.15cm}

\nhp (a)\; $x^2 = y^2 \neq\, z^2$\, and\, $x^2z^2 =\, y^2z^2 =\, x^2$.
\vspace{-0.15cm}

\nhp This amounts to\, $Z(z)\, \subset\, Z(x) = Z(y)$ (\ref{char-zeroset}). We focus on two alternatives:
\vspace{-0.15cm}

\nhp i)\; \und{No relations other than the above}.
\vspace{-0.1cm}

\nhp Routine checking shows that the following are all possible characters:
\vspace{-0.15cm}

\nhp --- \se h 1\, sends all three generators to\, 0;
\vspace{-0.15cm}

\nhp --- \se h 2, \se h 3\, send\, \y x, \y y\, to\, 0\, and, say,\, $\se h 2(z) = 1,\, \se h 3(z) = -1$;
\vspace{-0.15cm}

\nhp --- $\se h 4, \dots, \se h {11}$\, assign to the generators all possible combinations of values 
$\pm 1$, with, say,\, $\se h 4, \dots, \se h 7$\, sending\, \y z\, to\, 1, and\, 
$\se h 8, \dots, \se h {11}$\, sending\, \y z\, to $-1$.
\vspace{-0.15cm}

\nhp Call\, \se F 2\, the TS corresponding to this case. Using Lemma \ref{char-specializ} one sees 
at once that the specialization root-system of the ARS dual to\, \se F 2\, looks\vspace*{-1.3cm} as in Figure 3. 

$$\xymatrix@C=6pt{
 & & & & & & \\
 & & & & & \node \ulab{$\;h_1$} & & & & &\\
 & & \node \ulab{$h_2$} \edge{urrr}
 & & & & & &
     \node \ulab{$\;h_3$} \edge{ulll} & & \\
 \node \dlab{$h_4$} \edge{urr}
 & \node \dlab{$\;h_5$} \edge{ur}
 &
 & \node \dlab{$h_6$} \edge{ul}
 & \node \dlab{$\;h_7$} \edge{ull} 
 &  
 &  \node \dlab{$h_8$} \edge{urr}
 & \node \dlab{$h_9$} \edge{ur}
 &
 & \node \dlab{$\!\!h_{10}$} \edge{ul}
 & \node \dlab{$\;\;h_{11}$} \edge{ull} \\
 & & & & & & & & & & & } 
$$ 
\vspace{-1.1cm}

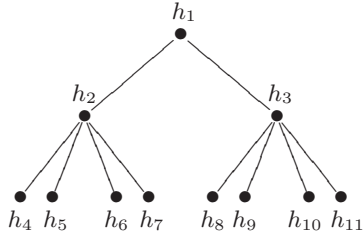
\nhp \hfl Figure 3. Specialization root-system of $\se X {F_2}$ \hfl
\vspace{-0.15cm}

Since\, $\se X {F_2}$\, has 11 elements, by \cite{DP5b}, Cor. 1.10, we must have\, card\,$(\se F 2) = 23$; 
the reader is invited to check that:
\vspace{-0.15cm}

\nhp \hspace{1.5cm} $\se F 2 = \{1, 0, -1, x, -x, y, -y, z, -z, x^2, -x^2, z^2, -z^2, xy, -xy, xz, -xz, $ 
\vspace{-0.25cm}

\nhp \hfl $yz, -yz, x^2z, -x^2z, xyz,-xyz\}.$\hspace{1.8cm} 

\nhp The Hasse diagram of the representation partial order of \se F 2 is drawn in Figure 4 below. 

Theorem \ref{reprpo-RS} is used in the computation of this diagram. 
For example, item (6) of the latter shows that\, $x^2 \leq x\,w \leq -x^2$ for all \y w. One uses 
item (4) to prove incomparability of elements of\, \se F 2\, as shown in Figure 4; 
for instance, to prove that\, $xyz$\, and\, \y y\, are $\leq$-incomparable, direct inspection of the 
characters of\, \se F 2\, described in item (i) above, shows that there are\, $h, h' \in \se X {F_2}$\, such 
that\, $h(y) = 1,\, h(xz) \in \{0, -1\}$ ---whence\, $h(xyz)\, \se {>} {\bf 3}\, h(y)$, i.e., 
$xyz \not\leq y$---, and\, $h'(y) = -1,\, h'(xz) \in \{0, -1\}$ ---hence\, $h'(y)\, \se {>} {\bf 3}\, h'(xyz)$, 
i.e., $y \not\leq xyz$. Further details are left to the reader.
\vspace{-1.2cm} 

$$\centerline{\hspace{0.1cm}\xymatrix@R=17pt@C=20pt{
 &  &  &  &  &  &  &  &  & \node \rulab{\mbox{$\raisebox{-4pt}{$\hspace{0.25cm}{-1}$}$}} &  &  &  & \\
 &  &  &  &  &  &  &  &  & \node \rulab{$\hspace{0.4cm}{-z^2}$} \edge{u} &  &  &  & \\
 &  &  &  &  &  &  &  &  & \node \rulab{$\hspace{0.4cm}{-x^2}$} \edge{u} &  &  &  & \\
 &  &  &  &  &  &  &  &  &  &  &  &  \\ 
 &  \node \sllab{$-z$} \edge{uuurrrrrrrr} \edge{dddrrrrrrrr}
   & \node \sllab{$-x$} \edge{uurrrrrrr}  \edge{ddrrrrrrr}
   & \node \sllab{$-y\!$} \edge{uurrrrrr}  \edge{ddrrrrrr} 
   & \node \sllab{$-xy\;$} \edge{uurrrrr} \edge{ddrrrrr}
   & \node \sllab{$-xz\;$} \edge{uurrrr} \edge{ddrrrr}
   & \node \sllab{$-yz\;$} \edge{uurrr} \edge{ddrrr}
   & \node \sllab{\mbox{$\raisebox{3pt}{$-x^2z\;\;\,$}$}} \edge{uurr} \edge{ddrr}
   & \node \sllab{$-xyz\;\;\;$} \edge{uur} \edge{ddr}
   & \node \srlab{$0\;\;$} \edge{uu}
   & \node \srlab{\mbox{$\raisebox{-5pt}{$\;\;xyz$}$}} \edge{uul} \edge{ddl} 
   & \node \srlab{\mbox{$\raisebox{3pt}{$\;\;x^2z$}$}} \edge{uull} \edge{ddll}
   & \node \srlab{\mbox{$\raisebox{-5pt}{$yz$}$}} \edge{uulll} \edge{ddlll}
   & \node \srlab{$xz$} \edge{uullll} \edge{ddllll}
   & \node \srlab{\mbox{$\raisebox{-6pt}{$xy$}$}} \edge{uulllll} \edge{ddlllll}
   & \node \srlab{\mbox{$\raisebox{-6pt}{$y\;\;$}$}} \edge{uullllll}  \edge{ddllllll}
   & \node \srlab{$\!\!x$} \edge{uulllllll}  \edge{ddlllllll}
      & \node \srlab{$z\;\;$} \edge{uuullllllll} \edge{dddllllllll} & \\ 
 &  &  &  &  &  &  &  &  &  &  &  &  \\ 
 &  &  &  &  &  &  &  &  & \node \rlab{\mbox{$\raisebox{-14pt}{$\!\!x^2$}$}} \edge{uu} & \\
 &  &  &  &  &  &  &  &  & \node \rlab{\mbox{$\raisebox{-11pt}{$\!\!z^2$}$}} \edge{u} & \\
 &  &  &  &  &  &  &  &  & \node \rlab{\mbox{$\raisebox{-14pt}{$\!1\;\;$}$}} \edge{u} & \\
 &  &  &  &  &  &  &  &  &  & \\
 }}
 $$
 
\vspace{-1cm} 

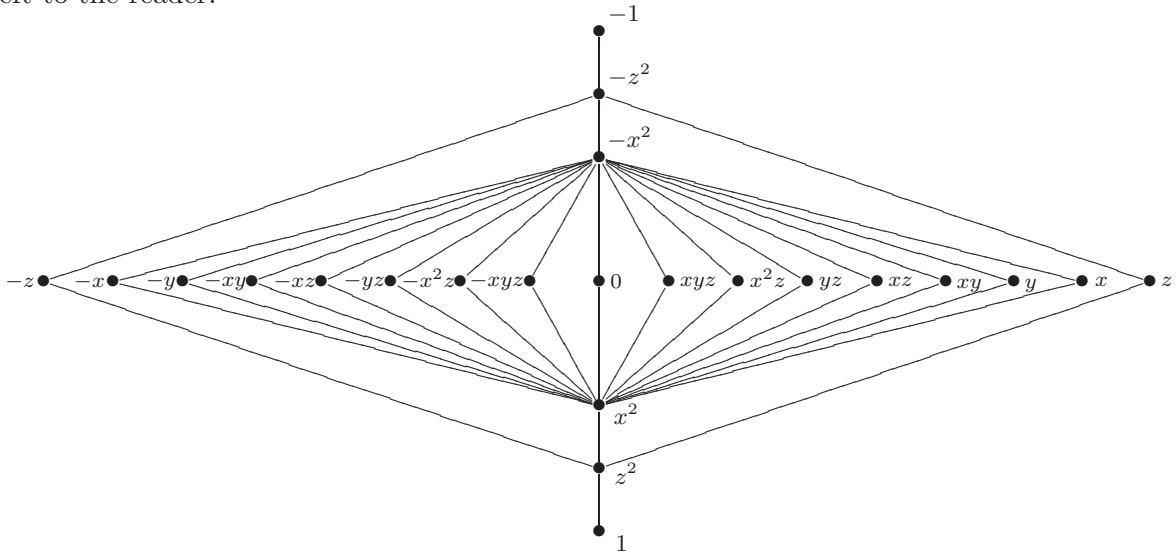
\nhp \hfl Figure 4. Representation partial order of\ \se F 2. \hfl

\vspace{-0.15cm}

One may also consider fans arising by adding relations between generators; as an example we describe
the fan obtained from the preceding one by adding:
\vspace{-0.15cm}

\nhp ii)\, \und{The extra relation\, $xz = x$}.
\vspace{-0.15cm}

\nhp Under this relation, each character sending\, \y z\, to\, $-1$\, must also send\, \y x\, to\, 0. Thus, 
the characters\, $\se h 8, \dots, \se h {11}$\, in the preceding example disappear. Specialization 
among the remaining characters does not change; the order of specialization of the resulting ARS 
is obtained by omitting these characters in Figure 3. The resulting RS-fan is:
\vspace{-0.2cm}

\nhp \hfl $\se {F} 3 = \{1, 0, -1, x, -x, y, -y, z, -z, x^2, -x^2, z^2, -z^2, xy, -xy\}.$ \hfl
\vspace{-0.15cm}

\nhp The diagram of its representation partial order is computed in much the same way as in the 
preceding example; details are left to the reader.
\vspace{-0.1cm}

\nhp iii)\, \und{The relations\, $xz = x$ and $z^2 = 1$.}
\vspace{-0.1cm}

\noi The additional relation $z^2 = 1$ makes $z$ invertible and hence excludes the character $\se h 1$ 
sending $z$ to $0$. This makes the characters\, \se h 2\,, \se h 3\, become ``disconnected"
(cf. \cite{DP5b}, Def. 2.14\,(a)); we obtain a two-component root-system:

\vspace{-1cm}

\nhp \hspace*{4.3cm} \parbox[t]{200pt} {

$\xymatrix{
 & & & \\
 & & \node \ulab{$h_2$} \\
   \node \dlab{$h_4$} \edge{urr}
 & \node \dlab{$h_5$} \edge{ur}
 &
 & \node \dlab{$h_6$} \edge{ul}
 & \node \dlab{$h_7$} \edge{ull} & \\
 & & & & & & } 
$ } 
 \hspace*{-3.5cm} \parbox[t]{100pt} {
$\xymatrix{
 & & & \\
 & & \node \ulab{$h_3$} \\ }
 $ } \hfl
 
\vspace{-0.8cm}

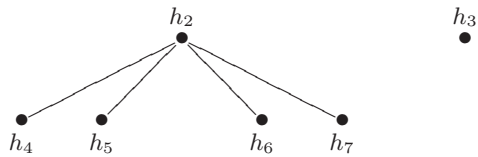
\nhp \hfl {Figure 5. Specialization root-system of\, $\se X {F_4}$\;\;\;\;} \hfl

\vspace{-0.15cm}

\nhp The dual RS-fan is:\;\; $\se {F} 4 = \{1, 0, -1, x, -x, y, -y, z, -z, x^2, -x^2, xy, -xy\},$\, 
with representation partial order:

\vspace{-0.95cm}

$$\xymatrix@R=18pt@C=22pt{
 &  &  &  &  & \node \rulab{\mbox{$\raisebox{-7pt}{$\;\;\;-1$}$}} &  &  &  & \\
 &  &  &  &  & \node \rulab{$\;\;\;\;\;-x^2$} \edge{u} &  &  &  & \\
 &  \node \llab{$-z\,$} \edge{uurrrr} \edge{ddrrrr}
   & \node \llab{$-x\;$} \edge{urrr}  \edge{drrr}
   & \node \llab{$-y$\,} \edge{urr}  \edge{drr} 
   & \node \llab{$-xy\;\;$} \edge{ur} \edge{dr}
   & \node \rlab{$0\;\;$} \edge{u}
   & \node \rlab{$xy\;\,$} \edge{ul} \edge{dl} 
   & \node \rlab{$\!\!y\;\;$} \edge{ull}  \edge{dll}
   & \node \rlab{$x\;\;\;$} \edge{ulll}  \edge{dlll}
   & \node \rlab{$z\;\;\;$} \edge{uullll} \edge{ddllll} & \\
 &  &  &  &  & \node \rlab{\mbox{$\raisebox{-14pt}{$\!\!x^2$}$}} \edge{u} & \\
 &  &  &  &  & \node \rlab{\mbox{$\raisebox{-14pt}{$\!1\;\;$}$}} \edge{u} & \\
 &  &  &  &  &  &  &  \\
 }
 $$
 
\vspace{-1cm}

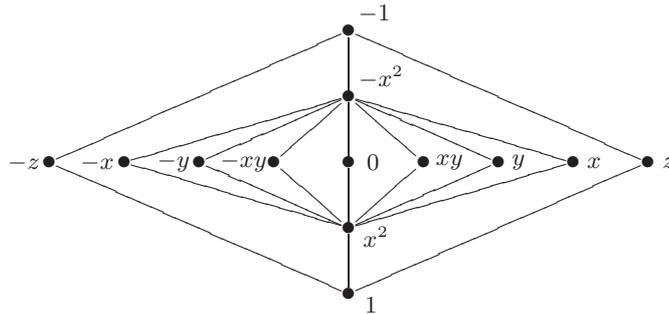
\nhp \hfl Figure 6. Representation partial order of\, $\se {F} 4$ \hfl 

Summarizing a common feature of the examples presented above, we shall prove that, under the representation 
partial order $\leq$\,, every RS-fan is a bounded lattice. 

\vspace{-0.4cm}

\ble \label{incomp} Let $F$ be a RS-fan and let $a, b \in F$. Then,
\vspace{-0.15cm}

\noi $(1)$ If $a < b$, then $a \in \mbox{$\mathrm {Id}$}(F)$ or $b \in -\mbox{$\mathrm {Id}$}(F)$.
\vspace{-0.15cm}

\noi In particular,
\vspace{-0.15cm}

\noi $(2)$ If $a,b \not\in \mbox{$\mathrm {Id}$}(F) \cup -\mbox{$\mathrm {Id}$}(F)$ ---i.e., 
$\pm\, a$ and $\pm\, b$ are not squares---, and $a \neq b$, then $a,b$ are incomparable under $\leq$.
\ele
\vspace{-0.35cm}

\nhp {\bf Proof.} (2) follows easily from (1).
\vspace{-0.15cm}

\noi (1) The assumption entails $a \in D(1,b)\,, -b \in D(1,-a)$ and $a \neq b$ (\ref{reprord}).
The definition of representation in a RS-fan (Theorem \ref{fanrepr}) yields
\vspace{-0.15cm}

\noi \hfl $a \in \mbox{$\mathrm {Id}$}(F)\: \jo\: a \in b\cdot \mbox{$\mathrm {Id}$}(F)\: 
\jo\: a = -ab$ \;\; and\;\; $-b \in \mbox{$\mathrm {Id}$}(F)\: \jo\: b \in a\cdot 
\mbox{$\mathrm {Id}$}(F)\: \jo\: b = ab$. \hfl
\vspace{-0.15cm}

Combining these alternatives gives:
\vspace{-0.15cm}

\noi --- If any one of the first disjuncts occur, we are done. 
\vspace{-0.15cm}

\noi --- If both the middle disjuncts occur, i.e., $a = bx^2$ and $b = ay^2$ for some 
$x,y \in F$, we get $a = ax^2y^2,\, b = bx^2y^2$ and $by^2 = ay^2$, whence 
$b = by^2x^2 = ay^2x^2 = a$, contradiction. 
\vspace{-0.15cm}

\noi --- If $a = bx^2$ and $b = ab$, then $a = bx^2 = abx^2 = a^2 \in \mbox{$\mathrm {Id}$}(F)$.
\vspace{-0.15cm}

\noi --- Likewise, $b \in a \cdot \mbox{$\mathrm {Id}$}(F)$ and $a = -ab$ yield 
$-b \in \mbox{$\mathrm {Id}$}(F)$.
\vspace{-0.15cm}

\noi --- The two last disjuncts imply $a = -b$, whence $a = -ab = a^2 \in \mbox{$\mathrm {Id}$}(F)$.
\hfl $\Box$

\vspace{-0.3cm}

\ble \label{upperbound} Let $F$ be a RS-fan and let $x, b \in F$. If $b \not\in \mbox{$\mathrm {Id}$}(F)$ 
$($i.e., $b \neq b^2)$, then $b \leq -x^2 \leq -b^2$ implies $b^2 = x^2$. That is, $-b^2$ is the smallest 
$y \in \mbox{$\mathrm {Id}$}(F)$ such that $b \leq -y$. Dually, $b^2$ is the largest 
$y \in \mbox{$\mathrm {Id}$}(F)$ such that $y \leq b$.
\ele
\vspace{-0.4cm}

\nhp {\bf Proof.} Assume $b \leq -x^2 < -b^2$. Since $b^2 \leq 0 \leq -x^2 \leq -b^2$, we get 
$b^2x^2 = x^2$ (and $Z(b)\, \sub\, Z(x))$ (Theorem \ref{reprpo-RS}\,(3),(5)). On the other hand, 
$b \leq -x^2$ yields $b \in D(1,-x^2)$, and then $b \in D^t(b^2,-b^2x^2) = D^t(b^2,-x^2)$. If 
$Z(b) \subset Z(x)$, the first clause in \ref{fanrepr} gives $b = b^2$, contrary to assumption. 
So, $Z(b) = Z(x)$, and we get $b^2 = x^2$. The dual assertion is obvious. \hfl $\Box$
\vspace{-0.5cm}

\bth \label{fanlattice} Let $F$ be a RS-fan and let\, $\leq$\, denote 
its representation partial order {\em (\ref{reprord})}. Then, $(F,\leq)$ is 
a lattice with smallest element\, $1$ and largest element $-1$.
\eth

\vspace{-0.4cm}

\noi {\bf Notation.} For elements $a,b$ in a RS, $G$, we write $a \perp b$
to mean that $a$ and $b$ are incomparable under the representation partial
order of $G$. \hfl $\Box$

\nhp {\bf Proof.} We must show that every pair of elements $a,b \in F$ 
has a least upper bound, $\jo$, and a greatest lower bound, $\w$\,, for 
the order $\leq$. If $a,b$ are comparable under $\leq$ there is nothing 
to prove; so, we may assume $a \perp b$.
\vspace{-0.15cm}

Since $F$ is a RS-fan, the zero-sets of $a$ and $b$ are comparable 
under inclusion. This, together with $a \perp b$, implies that one of 
$a$ or $b$ is not in $\mbox{\textrm Id}(F) \cup -\mbox{\textrm Id}(F)$; 
indeed:
\vspace{-0.15cm}

\nhp --- If $Z(a)\, \sub\, Z(b)$, \ref{reprpo-RS}\,(5) yields $a^2 \leq b \leq -a^2$; 
hence, $a \perp b$ implies $a \not\in\mbox{\textrm Id}(F) \cup -\mbox{\textrm Id}(F)$.
\vspace{-0.15cm}

\nhp --- Likewise, $Z(b)\, \sub\, Z(a)$ implies $b \not\in \mbox{\textrm Id}(F) \cup -\mbox{\textrm Id}(F)$.
\vspace{-0.15cm}

Since $-a^2$ and $-b^2$ are $\leq$-comparable (\ref{reprpo-RS}\,(5)), we may assume without 
loss of generality that $-a^2 \leq -b^2$. Further, Lemma \ref{upperbound} shows
\vspace{-0.15cm}

\nhp (*) \hfl $-b^2 =$ least $x \in -\mbox{\textrm Id}(F)$ such that $a,b \leq x$. \hfl
\vspace{-0.1cm}

\nhp {\bf Claim.}  $-b^2 = a\, \jo\, b$.
\vspace{-0.15cm}

\nhp {\bf Proof of Claim.} By assumption, $a,b \leq -b^2$; so we need only prove:
\vspace{-0.1cm}

\nhp \hfl $\fa\, c \in F\,(c \geq a,b\;\; \Ra\;\; c \geq -b^2)$. \hfl
\vspace{-0.1cm}

\nhp Note that $c \geq a,b$ and $a \perp b$ imply $c \neq a,b$. If 
$c \not\in \mbox{\textrm Id}(F) \cup -\mbox{\textrm Id}(F)$, since one of $a,b$ ---say $a$--- is not in 
$\mbox{\textrm Id}(F) \cup -\mbox{\textrm Id}(F)$, then, by Lemma \ref{incomp}, $c \neq a$ implies 
$c \perp a$, absurd; hence $c \in \mbox{\textrm Id}(F) \cup -\mbox{\textrm Id}(F)$. If 
$c \in \mbox{\textrm Id}(F)$, then $c \geq a,b$, implies $a,b \in \mbox{\textrm Id}(F)$, whence $a,b$ 
are $\leq$-comparable, contradiction. So, $c \in -\mbox{\textrm Id}(F)$, and (*) gives $c \geq -b^2$, as claimed.
\vspace{-0.15cm}

Under the current assumptions $-a^2 \leq -b^2$ and $a \perp b$; upon observing that
\vspace{-0.15cm}

\nhp \hfl $b^2 =$ largest $y \in \mbox{\textrm Id}(F)$ such that $y \leq a,b$, \hfl
\vspace{-0.15cm}

\nhp a similar argument yields\, $b^2 = a\, \w\, b$.  \hfl $\Box$
\vspace{-0.4cm}

\bres \label{rmks fanlattice} (a) A closer look at the examples presented above shows that 
the lattices $(F, \leq)$ are not modular ---hence not distributive either--- except in 
very special cases. In fact, most of these lattices contain the configuration
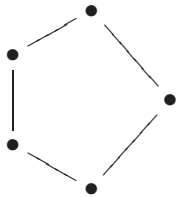
\vspace{-0.7cm}

$$\xymatrix@R=6pt@C=6pt{
& & {\bullet} \ar@{-}[dll] \ar@{-}[ddrr] & & \\
  {\bullet} \ar@{-}[dd] & & & & \\
& & & & {\bullet} \ar@{-}[ddll] \\
  {\bullet} \ar@{-}[drr] & & & & \\
& & {\bullet} & & \\
& & & & & \\
}
$$

\vspace{-0.7cm}

\nhp as a sublattice (cf. \cite{B}, Ch. V, \S\,2, Thm. 2, p. 66). For instance, in Figures 4 and 6 
above, the sublattices $\{z^2 < x^2 < -x^2 < -z^2; z\}$ and $\{1 < x^2 < -x^2 < 1; z\}$, respectively, 
form such a pentagon. Example \ref{ex-fans}.A is modular but not distributive. Note that a RSG-fan
(i.e., a reduced special group that is a fan, cf. [RSG-fan], Introduction) {\it is a modular lattice} 
under the order $a \leq b\; \Lra\; a \in D(1,b)$. We also register that in \cite{DP3}, Thm. 6.6, pp.  396-397, 
we proved that the representation partial order in {\it spectral} real semigroups is a {\it distributive} lattice.
\vspace{-0.15cm}

\nhp (b) Since $\mbox{\textrm Id}(F) \cup -\mbox{\textrm Id}(F)$ is a totally ordered subset of $(F,\leq)$,
the proof of Theorem \ref{fanlattice} shows that the lattice operations in $(F,\leq)$ satisfy the following identities: 

${\hspace{2cm} a\, \w\, b\; =\;} 
\left\{
\begin{array}{ll}
\mbox{min}_{\leq} \{a^2, b^2\} &  \;\;\;\mbox{if}\;\; a \perp b\\
\mbox{min}_{\leq} \{a,b\} &  \;\;\;\mbox{if}\;\; a, b \;\mbox{are} \leq\!\mbox{-comparable}, 
\end{array}
\right.
$ 
\vspace{-0.2cm}

\noi and
\vspace{-0.2cm}

${\hspace{2cm} a\, \jo\: b\; =\;} 
\left\{
\begin{array}{ll}
\mbox{max}_{\leq} \{-a^2, -b^2\} &  \;\;\;\mbox{if}\;\; a \perp b\\
\mbox{max}_{\leq} \{a,b\} &  \;\;\;\mbox{if}\;\; a, b\;\mbox{are} \leq\!\mbox{-comparable}.
\end{array} 
\right. 
$ 

\nhp Note that, if $a \perp b$, then $a\, \w\: b,\, a\, \jo\: b \in \mbox{\textrm Id}(F) \cup -\mbox{\textrm Id}(F)$.
\vspace{-0.15cm}

\nhp (c) The operation $x \mapsto -x\; (x \in F)$ {\it is not} a complement in the lattice-theoretic sense, 
but it verifies:
\vspace{-0.15cm}

\nhp $(\se {\mbox{\textrm c}} 1)$ The Kleene inequality\; $a\; \w -a\, \leq\, 0\, \leq\, b\; \jo -b$.  
(A special case of Theorem \ref{reprpo-RS}\,(10).) 
\vspace{-0.15cm}

\nhp $(\se {\mbox{\textrm c}} 2)$ The De Morgan laws:
\vspace{-0.15cm}

\nhp \hfl (i)\;\; $-(a\, \w\: b\,) = -a\: \jo -b$\,; \hspace{1.2cm} (ii)\;\; $-(a\, \jo\: b\,) = -a\: \w -b$\,. \hfl
\vspace{-0.15cm}

\nhp This is clear if $a$ and $b$ are comparable under $\leq$. If $a \perp b$, assuming, without 
loss of generality, $-a^2 \leq -b^2$, (i.e., $b^2 \leq a^2$), from (b) we get:
\vspace{-0.15cm}

\nhp (i)\; $-(a\, \w\: b\,) = -(a^2\, \w\: b^2) = -b^2$,\; and\; $-a\, \jo -b = -(-a)^2\, \jo -(-b)^2 = 
-a^2\, \jo -b^2 = -b^2$.
\vspace{-0.15cm}

\nhp (ii) $-(a\, \jo\: b\,) = -(-a^2\, \jo -b^2) = -(-b^2) = b^2$,\, and\, $-a\, \w -b = 
(-a)^2\, \w\, (-b)^2 = a^2 \w\: b^2 = b^2$.\hfl $\Box$
\eres

\vspace{-0.8cm}

\section {Quotients of fans} \label{fan-quot}

\vspace{-0.2cm}

\paragraph{A. Preliminaries.\vspace{-0.3cm} Congruences of ternary semigroups and real semigroups.} \label{congr}

\vspace{-0.2cm}

\bdf \label{TS-congr} A congruence of ternary semigroups $($abbreviated 
{\bf TS-congruence}$)$ \index{sub}{ternary semigroup (TS)!congruence} 
\index{sub}{congruence!of ternary semigroup}\index{sub}{TS-congruence} 
is an equivalence relation\, 
$\equiv$\, on a TS, $G$, compatible with the semigroup operation and such 
that the induced quotient structure $G/\!\!\equiv$ is a ternary semigroup. 
$[$This is equivalent to require $\equiv$\, to be proper, i.e. 
$\equiv\; \subset G \times G$, and for $x \in G$, $\;\;x \equiv -x\;\; \Ra\;\; 
x \equiv\, 0$.$]$  \hfl $\Box$
\edf

\vspace{-0.75cm}

\bres \label{propsTS-congr}(a) The condition that\, $\equiv$\, is proper ensures 
that $1 \not\equiv\, 0$, and hence, by the last requirement, $1 \not\equiv -1$.
\vspace{-0.15cm}

\noindent (b) Since the axioms for TSs are universal, the quotient map\, 
$\se {\pi} {\equiv}: G \longrightarrow G/\!\equiv$\, is automatically a TS-homomorphism.
\vspace{-0.15cm}

\noindent (c) For each non-empty set\, $\cal H \subseteq \HTS G$, the relation
\vspace{-0.1cm}

\noindent $\se {(\dag)} {\cal H}$ \hfl $a\; \se {\equiv} {\cal H}\, b$\;\; $\Leftrightarrow$ \;\; 
For all $h \in \cH,\; h(a) = h(b),$\;\;\;$(a,b \in G),$ \hfl

\vspace{-0.1cm}

\noindent defines a TS-congruence of $G$ (straightforward checking). 
We shall write $G/\cH$ for the quotient TS\; $G/\se {\equiv} {\cal H}.$ 
\vspace{-0.15cm}

\noi (d) The set\, $\cH\, \sub\, \HTS G$ can be identified with a subset $\h {\cH} =$ 
$\{\h h\, \vert\, h \in \cH\}$ of ${\bf 3}^{\,G/{\cal H}}$ by the map\, $h \mapsto \h h$, where\, 
$\h h: G/{\cal H}\; \lra\; {\bf 3}$ is defined by the functional equation\, $\h h \circ \pi = h$. 
By clause\, $\se {(\dag)} {\cal H}$ above, \h h is well-defined and the map\, $h \mapsto \h h$ is 
obviously injective. The reader can easily check that $\HTS {G/{\cal H}} = \h {\cH}$. Further, if
$\HTS G$ and $\HTS {G/{\cal H}}$ are endowed with their spectral topologies  and $\cH$ is a 
proconstructible subset of $\HTS G$ closed under the product of any three of its members, the 
map above induces a homeomorphism between $\HTS {G/{\cal H}}$ and $\cH$ (\cite{DP4}, Thm. I.1.26).
\vspace{-0.15cm}

\noi (e) In \cite{DP4}, Thm. I.1.26, it is shown that every TS-congruence of a 
ternary semigroup \y G is of the form\, $\se {\equiv} {\cal H}$\, for a 
suitable set\, $\cal H$\, of TS-characters. The set $\cal H$\, can even be
taken to be proconstructible in $\HTS G$. \hfl $\Box$
\eres

\vspace{-0.75cm}
\bdf \label{defi:coc} A {\bf (RS-)congruence} of a real semigroup\, \y G\, 
is an equivalence relation\, $\equiv$\, satisfying the following requirements:
\vspace{-0.15cm}

\nhp $($i$)$\; $\equiv$\, is a congruence of ternary semigroups (\ref{TS-congr}).
\vspace{-0.15cm}

\nhp $($ii$)$ There is a ternary relation\, $\se D {G/\equiv}$\, in the quotient {\it ternary semigroup}\, 
$(G/\!\equiv,\cdot,-1,0,1)$\, so that\, $(G/\!\equiv,\cdot,\se D {G/\equiv},-1,0,1)$\, is a real semigroup, 
and the canonical projection\, $\pi:G\, \lra\, G/\!\equiv$\, is a RS-morphism.
\vspace{-0.15cm}

\nhp $($iii$)$ (Factoring through\, $\pi$.) For every RS-morphism\, $f:G\, \lra\, H$\ into a real semigroup\, 
\y H\, such that\, $a \equiv b$\, implies\, $f(a) = f(b)$\, for all\, $a,b \in G$, there exists a RS-morphism 
$($necessarily unique$)$,\, $\h{f}:G/\!\equiv\, \longrightarrow\, H$,\, such that\, $\h{f} \circ \pi = f$, i.e. 
the following diagram commutes
\vspace{-0.25cm} 

\begin{center}
\unitlength=1mm
\begin{picture}(20,20)(0,0)
\put(8,16){\vector(3,0){5}} \put(2,10){\vector(0,-2){5}}
\put(10,4){\vector(1,1){7}}
\put(2,16){\makebox(0,0){$G$}} \put(20,16){\makebox(0,0){$H$}}
\put(2,0){\makebox(0,0){$G/\!\equiv$}}
\put(2,20){\makebox(17,0){$f$}} \put(-2,8){\makebox(-6,0){$\pi$}}
\put(18,2){\makebox(-4,3){$\h{f}$}}
\end{picture}
\end{center}

\edf

\vspace{-0.7cm} 

\bdf \label{ReprInQuot} With notation as in \ref{propsTS-congr}, if $G$ is a RS and 
$\cH\, \subû\, \HRS G = \se X G$, we define a ternary relation\, $\se D {G/{\cal H}}$\, on\, 
$G/{\cal H}$\, as follows: for\, $a,b,c \in G$,

\vspace{-0.15cm} 

\nhp $\se {(\dag\dag)} {\cal H}$ \hfl $\pi(a) \in \se D {G/{\cal H}}(\pi(b),\pi(c))$\;\; $\Lra$\;\; 
For all $h \in \cH,\; h(a) \in \se D {\bf 3}(h(b),h(c))$. \hfl

\edf

\vspace{-0.45cm} 

\noindent Cf. clause $\se {[D]} {\cal H}$ in Definition \ref{ReprD_H}. Obviously\, $\se D {G/{\cal H}}$\, 
is well-defined and every RS-congruence of real semigroups is obtained in this way. To be precise, we 
state the following result, a straightforward consequence of the separation theorems for RSs and for 
TSs (\cite{DP1}, Thms. 4.4, pp. 116--117 and 1.9, pp. 103-104): 
\vspace{-0.4cm}

\bpr \label{pr:iso} Given a real semigroup $G$ and a RS-congruence\, $\equiv$ of\, $G$, let
\vspace{-0.1cm}

\nhp \hfl $\se {\cH} {\equiv} = \{p \in \se X {G}\, \vert\, {\mbox{\textrm There\ exists\ }}\; \sigma \in 
\se X {G/\equiv}\ {\mbox{\textrm so\ that\ }}\;  p = \sigma \circ \pi\}.$ \hfl

\vspace{-0.3cm}

\noi Then, 
\vspace{-0.15cm}

\nhp $(1)$\;\: For\, $a,b,c \in G$,

\vspace{-0.2cm}

\hspace*{0.25cm}\parbox[t]{370pt} {\hspace*{0.13cm}$(i)$\;  $\pi(a) \in \se D {G/\equiv}(\pi(b),\pi(c))\;\; \Lra\;\;$ For all 
$p \in \se \cH {\equiv},\; p(a) \in \se D {\bf 3}(p(b),p(c))$. 

$(ii)$\; $a \equiv b\;\; \Lra\;\;$ For all $p \in \cH_{\equiv},\; p(a) = p(b)$.} 

\noi $(2)$ $(G/{\equiv}, \D {G/\equiv}) \cong (G/{\cal H}_{\equiv}, \D {G/{\cal H}_{\equiv}})$ 
as $\se {\cal L} {\mathrm {RS}}$-structures.
\vspace{-0.15cm}

\noi Hence,
\vspace{-0.15cm}

\noi $(3)$ $(G/{\cal H}_{\equiv}, \D {G/{\cal H}_{\equiv}})$ is a RS. \hfl $\Box$
\epr

\vspace{-0.8cm}

\paragraph{B. Congruences of fans.} We shall now consider the structure of congruences of 
RS-fans, giving an explicit description of them, and proving, in particular, that arbitrary 
quotients of fans are fans. We start with some preliminary observations used in the proof.
\vspace{-0.4cm}

\bfa \label{CondZ-Quot} Let $G$ be a TS satisfying condition $[$Z\,$]$ in $\ref{fanrepr}$ and
let $\0 \neq \cH\, \sub\, \HTS G$. Then, the quotient TS, $G/\cH$, also satisfies
condition $[$Z\,$]$.
\efa
\vspace{-0.45cm}

\nhp {\bf Proof.} Follows from (1) $\Lra (2)$ in Proposition \ref{charfan}. \hfl $\Box$
\vspace{-0.4cm}

\bfa \label{D-fanQuot} Let $G$ be a RS and let $\0 \neq \cH\, \sub\, \HTS G$. With the relation
$\D {G/\cal H}$ defined by clause $\,\se {(\dag\dag)} {\cal H}$ in $\ref{ReprInQuot}$, for
$a,b \in G$, we have 
\vspace{-0.5cm}

\begin{tabbing}
\hspace*{0.3cm}$\pi(a) \cdot {\mathrm {Id}}({G/\cal H}) \cup \pi(b) \cdot {\mathrm {Id}}({G/\cal H})) 
\cup \{\pi(c) \sth \pi(c)\pi(a) = -\pi(c)\pi(b)\; and $ \=$\pi(c) = \pi(a^2)\pi(c)\}\, \sub $ \\ [0.1cm]
\> $\sub\, \D {G/\cal H}(\pi(a), \pi(b))$.
\end{tabbing}
\efa
\vspace{-0.55cm}

\nhp {\bf Proof.} See Remark \ref{RmkOn[D]PrecThm}. \hfl $\Box$

\vspace{-0.4cm}

\bpr \label{pr:fanQuot-1} Let\, $F$\, be a RS-fan and let\, $\cH$\, be a non-empty subset 
of $\se X F$ \vspace{0.04cm} which is $3$-closed $($i.e., stable under product of any three 
of its elements$)$. Then the quotient\, $F/\cH$\, is a RS-fan\, $($and $\se {\equiv} {\cal H}$\, 
is a RS-congruence$)$.  
\epr

\vspace{-0.4cm}
\nhp {\bf Proof.} With $\pi: F\, \lra\, F/{\cal H}$ denoting the quotient map, we must show, 
for $a,b,c \in F$:
\vspace{-0.45cm}

\begin{tabbing}
\hspace{1.5cm} $\pi(c) \in \D {F/\cal H}(\pi(a), \pi(b))\;\; \Lra\;$ \= $\pi(c) \in \pi(a) \cdot 
{\mathrm {Id}}({F/\cal H})\; \jo$ \hspace*{4cm} \=[1]\\
\> $\pi(c) \in \pi(b) \cdot {\mathrm {Id}}({F/\cal H})\; \jo$ \>[2]\\[0.08cm]
\> $\pi(c)\pi(a) = -\pi(c)\pi(b)$ and $\pi(c) = \pi(a^2)\pi(c)$. \>[3]
\end{tabbing}
\vspace{-0.5cm}

\noi ($\La$) is the content of Fact \ref{D-fanQuot}. 
\vspace{-0.15cm}

\noi ($\Ra$) Assume $\pi(c) \in \D {G/\cal H}(\pi(a), \pi(b))$ and [1]--[3] false. By 
$\se {(\dag\dag)} {\cal H}$ in \ref{ReprInQuot}, the negation of [1] and [2] yield\vspace{0.05cm}
characters $\se h 1, \se h 2 \in \cH$ such that $\se h 1(c) \neq \se h 1(a) \se h 1(c^2)$
and $\se h 2(c) \neq \se h 2(b) \se h 2(c^2)$. The negation of [3] is equivalent to
[3.i] \jo\,[3.ii], where
\vspace{-0.15cm}

\noi [3.i]\;\; $\pi(c)\pi(a) \neq -\pi(c)\pi(b)$,
\vspace{-0.15cm}

\noi [3.ii]\; $\pi(c)\pi(a) = -\pi(c)\pi(b)$ and $\pi(c) \neq \pi(a^2)\pi(c)$.
\vspace{-0.15cm}

\noi Case [3.i] yields $\se h 3 \in \cH$ so that $\se h 3(c) \se h 3(a) \neq -\se h 3(c) \se h 3(b)$.
These inequalities obviously imply $\se h i(c) \neq 0$, i.e., $\se h i(c^2) = 1$, for $i = 1,2,3$,
and $\se h 1(c) \neq \se h 1(a), \se h 2(c) \neq \se h 2(b)$.
\vspace{-0.15cm}

Let $h := \se h 1 \se h 2 \se h 3$. Then, $h \in \se X F$ ($F$ is a fan), $h \in \cH$ ($\cH$ is
3-closed), and $h(c) \neq 0$. The representation assumption implies, then, $h(c) = h(a)$ or 
$h(c) = h(b)$; suppose the first equality holds. From $\se h 2(c) \neq \se h 2(b)$ we get
$\se h 2(c) = \se h 2(a)$, and hence
\vspace{-0.15cm}

\noi [4]\; $\se h 1(c) \se h 3(c) = \se h 1(a) \se h 3(a)$.
\vspace{-0.15cm}

\noi Since $\se h 1(c) \se h 3(c) \neq 0$, we have $\se h 1(a), \se h 3(a) \neq 0$; from
$\se h 1(c) \neq \se h 1(a)$ comes $\se h 1(c) = -\se h 1(a)$ and, by [4], 
$\se h 3(c) = -\se h 3(a)$. The representation assumption yields, then, $\se h 3(c) = \se h 3(b)$, 
wherefrom $\se h 3(b) = -\se h 3(a)$. Scaling by $\se h 3(c)$ we get
$\se h 3(c) \se h 3(b) = -\se h 3(c) \se h 3(a)$, contrary to assumption [3.i].
\vspace{-0.15cm}

The case $h(c) = h(b)$ is dealt with by a similar argument.
\vspace{-0.15cm}

\noi Case [3.ii] yields $\se h 3 \in \cH$ so that $\se h 3(c) \neq \se h 3(a^2) \se h 3(c)$.
Squaring the equality $\pi(c)\pi(a) = -\pi(c)\pi(b)$ and scaling by $\pi(c)$ we get
$\pi(c)\pi(a^2) = \pi(c)\pi(b^2)$, whence $\se h 3(c)\se h 3(a^2) = \se h 3(c)\se h 3(b^2)$
$\neq \se h 3(c)$. This implies $\se h 3(c) \neq 0, \se h 3(a) = 0$ and $\se h 3(b) = 0$.
Assumption $\pi(c) \in \D {G/\cal H}(\pi(a), \pi(b))$ and $\se h 3 \in \cH$ yield
$\se h 3(c) \in \Dthr(\se h 3(a), \se h 3(b)) = \Dthr(0,0)$, whence $\se h 3(c) = 0$,
contradiction. 
\vspace{-0.15cm}

That $\se {\equiv} {\cal H}$ is a RS-congruence follows from the fact, just proved,
that $F/{\cH}$ is a RS-fan (\ref{def:fan}). This completes the proof of Proposition 
\ref{pr:fanQuot-1}. \hfl $\Box$    

\vspace{-0.1cm}
 
Observe that {\it all}\, RS-congruences of a fan are obtained in the way given 
by the preceding Proposition:
\vspace{-0.45cm}

\bco \label{fancon-cor} Let\, $F$\, be a RS-fan and let\, $\equiv$\, be a 
RS-congruence of\, \y F. Then:
\vspace{-0.15cm}

\nhp $(a)$ $\equiv\; =\; \se {\equiv} {\cal H}$\, for some proconstructible, 
$3$-closed set\, $\cH\, \sub\, \se X F$. Hence,
\vspace{-0.15cm}

\nhp $(b)$ $F/{\equiv}$\, is a RS-fan.
\vspace{-0.15cm}

\nhp $(c)$ The correspondence $\cH \longmapsto \se {\equiv} {\cal H}$ 
establishes an inclusion-reversing bijection between proconstructible 
$3$-closed subsets of\, $\se X F$\, and the set of RS-congruences of\, \y F.
\eco

\vspace{-0.45cm}
\nhp {\bf Proof.} (a) The set\, $\cH =\, \se {\cal H} {\equiv}$\, is given by
Proposition \ref{pr:iso}. Item (b) follows from Proposition \ref{pr:fanQuot-1}, 
and item (c) is proved in \cite{DP4}, Thm. I.1.26. \hfl $\Box$

\noi {\bf Remark.} We register in passing that quotients of fans have a much stronger
property called {\it transversal $2$-regularity}, introduced and studied in \cite{DP4},
Ch. III, \S\,3. The proof of this property appears in \cite{DP4}, Thm. VI.11.3. \hfl $\Box$
\vspace{-0.5cm}

\paragraph{C. Quotients modulo ideals.} As a last point in this section we address the 
special case of quotients of RS-fans modulo ideals. Amongst the outstanding cases of 
congruences of a RS (cf. \cite{DP4}, Ch. II, \S\,3, \cite{M}, \S\S\,6.5, 6.6) one considers 
those determined by saturated prime ideals. 
\vspace{-0.15cm}

A saturated prime $I$ ideal of a RS, $G$, determines the set of characters 
$\se {\cH} I = \{h \in \se X G \sth$ $Z(h) = I\}$. The congruence $\se {\equiv} {{\cal H}_I}$ 
induced by $\se {\cal H} I$ will be denoted by $\se {\sim} I$, and the corresponding quotient
set by $G/I$.\,\footnote{\,Quotients of this type have been considered by Marshall in the dual 
category of abstract real spectra; cf. \cite{M}, p. 102 and Cor. 6.6.9.} In \cite{DP4}, Thm. 
II.3.15, we characterize the congruence $\se {\sim} I$ and both representation relations of 
$G/I$ solely in terms of the data carried by $G$. We also prove that the representation relation
$\D {G/I}$ induces on the set $\se G I := (G/I) \setminus \{\pi(0)\}$, obtained from $G/I$ by 
omitting zero, the structure of a {\it reduced special group}. Including proofs of these 
results in full generality will take us too far afield. However, in the special case of RS-fans this 
property follows from Proposition \ref{pr:fanQuot-1}. The following fact is used in the proof and 
elsewhere in this paper:

\vspace{-0.4cm}

\ble \label{congrfan} Let \y I be an ideal of a RS-fan \y F. Then, for\, $a,b \in
F \setminus I$: 
\vspace{-0.1cm}

\nhp \hfl $a \, \se {\sim} I \, b \;\, \Lra \;\, \ex z\, \not\in\, I \;(az = bz)$. \hfl 

\ele 

\vspace{-0.45cm}

\noi {\bf Proof.} The implication $(\La)$ is clear: if $az = bz$ with 
$z \not\in I$ and $h \in \se {\cH} I$, i.e., $Z(h) = I$, then $h(a)h(z) = h(b)h(z)$ 
entails $h(a) = h(b)$.
\vspace{-0.15cm}

\noi ($\Ra$) Assume $a \, \se {\sim} I \, b$; then $ab \: \se {\sim} I \, b^2$.
Since $b \not\in I$, we have $h(b) \neq 0$, i.e., $h(b^2) = 1$, for all
$h \in \se {\cH} I$, whence $ab \: \se {\sim} I \, b^2 \, \se {\sim} I \, 1$.
Set $p := ab$ and assume there is no $z \in F \setminus I$ such that $pz = z$,
i.e., $\{z \in F \sth pz = z\}\, \sub\, I$. Let $S$ be the subsemigroup of $F$
generated by $I \cup {\mbox{\textrm {Id}}}(F) \cup \{-p\}$. Clearly,
$S = I \cup {\mbox{\textrm {Id}}}(F) \cup -p \cdot {\mbox{\textrm {Id}}}(F)$.
Next, observe that $S \cap -S = I$. Indeed, if $x \in S \cap -S$, we have
$-x^2 \in S$. If $-x^2 \in I$, clearly $x \in I$. If $-x^2 \in {\mbox{\textrm {Id}}}(F)$,
then $x = 0 \in I$. Finally, if $-x^2 \in -p \cdot {\mbox{\textrm {Id}}}(F)$, then
$x^2 = pz^2$, whence $pz^2x^2 = z^2x^2$, yielding $z^2x^2 \in I$, and hence 
(by \ref{charfan}) $z \in I$ or $x \in I$. In either case we conclude $x \in I$, proving 
$S \cap -S\, \sub\, I$.
\vspace{-0.15cm}

Let $T$ be a subsemigroup of $F$ containing $S$ and maximal for $T \cap -T = I$. 
By Lemma 1.5, p. 102 of \cite{DP1} there is a character $h \in \se X F$ such
that $Z(h) = I$ and $T = h^{-1}[0,1]$. Since $-p \in S\, \sub\, T$ and $p \not\in I$, we get 
$h(p) = -1$, contradicting $p \: \se {\sim} I \, 1$. This proves that $pz = abz = z$ 
for some $z \in F \setminus I$. Scaling by $b^2$ we get $b(bz) = a(bz)$ with
$bz \in F \setminus I$, as required. \hfl $\Box$

\vspace{-0.35cm}

\bpr \label{fanQuot-2} Let \y F be a RS-fan. Let \y I be a proper ideal of\, \y F.
Let $\py = \se {\pi} I: F \, \lra \, F/I$ denote the canonical quotient map.
Then,\, $\se F I = (F/I) \setminus \{ \py(0) \}$ is a RSG-fan. 
\epr

\vspace{-0.45cm}

\nhp {\bf Proof.} It suffices to prove:
\vspace{-0.15cm}

\nhp (\dag) $\se F I$ is a group of exponent 2 with $1 \neq -1$.
\vspace{-0.15cm}

Indeed, a straightforward computation using clause [D] of \ref{fanrepr} and (\dag)
proves that $(\se F I, \D {F_I})$ satisfies condition [RSG-fan] (Introduction)
defining RSG-fans, i.e., for $a, b \in F$ so that \py(\y a), \py(\y b) $\neq 0$ 
and $\pi(b) \neq \pi(-1)$,
\vspace{-0.15cm}

\nhp \hfl $\pi(a) \in \se D {F/I}(\pi(1),\pi(b)) \;\, \Ra \;\, \pi(a) = \pi(1)
\; \jo \; \pi(a) = \pi(b)$. \hfl
\vspace{-0.15cm}

\noi \und{Proof of (\dag)}. We must prove: if $a \in F$ is\vspace{0.05cm} such that 
$\py(\y a) \neq 0$ (i.e., $a \not\in I$), then $\py(\y a^2) = 1$ (i.e., $a^2\, \se {\sim} I \, 1$).
By Lemma \ref{congrfan} the conclusion is equivalent to $\ex\,z \not\in I\, (a^2z = z)$; take $z = a$.
Note also that $\pi(1) \neq \pi(-1)$, i.e., $1\; {\se {\not\sim} I}\! -1$, since $1 \cdot z \neq
(-1) \cdot z$ holds for every $z \neq 0$. \hfl $\Box$
\vspace{-0.3cm}

\bres \label{CollapseInQuotients} (a) Quotients of real semigroups produce a considerable amount
of collapse. An example is the quotient of the RS-fan $\se F 3$ of 3.2.B\,(ii) by
the ideal $\{0\}$: since the relation $z \cdot x = x = 1 \cdot x$ holds in $\se F 3$, by 
\ref{congrfan} the generator $z \in \se F 3$ collapses to $1$ in $\se F 3/\{0\}$.
\vspace{-0.15cm}

\noi (b) There is also collapse ``from above'': any ideal $J \supset I$ collapses onto the
whole of $F/I$: $\se {\pi} I[J] = F/I$. 
\vspace{-0.15cm}

\noi \und{Proof}. Observe first that $\se {\pi} I(x^2) = 1$ for $x \in J \setminus I$: 
since $x^2 \cdot x = x = 1 \cdot x$ and $x \not\in I$, we have $x^2\, \se {\sim} I 1$. 
Now, let $a \in F$; we prove there is $y \in J$ so that $\se {\pi} I(a) = \se {\pi} I(y)$. 
If $a \in I$, take $y = 0$. If $a \not\in I$, let $x \in J \setminus I$; then, 
$y = ax^2 \in J \setminus I$, and
\vspace{-0.15cm}

\noi \hfl $\se {\pi} I(ax^2) = \se {\pi} I(a) \se {\pi} I(x^2) = \se {\pi} I(a) \se {\pi} I(1) =
\se {\pi} I(a)$. \hfl
\vspace{-0.15cm}

\noi (c) It can be proved that the inverse images ${\se {\pi} I}^{-1}[\De]$ of proper
(saturated) subgroups $\De$ of the RSG-fan $\se F I$ are exactly the\vspace*{0.03cm}
saturated subsemigroups $\Ga\, \sub\, F$ such that $\Ga\, \cap\, I = \0$ and 
$\Ga \supseteq {\se {\pi} I}^{-1}[1]$. The proof is omitted. \hfl $\Box$
\eres

\vspace{-0.8cm}

\section{Characterizations of fans} \label{sect:char-fans}

\vspace{-0.45cm}

The main result of this section is the following characterization of RS-fans:
 
\vspace{-0.45cm}
\bth \label{thm:char-fan} For a real semigroup $G$, the following are equivalent:

\vspace{-0.2cm}
\nhp $(1)$ $G$ is a RS-fan.

\vspace{-0.2cm}
\nhp $(2)$ $G$ satisfies the following conditions:

\vspace{-0.15cm}
\tab{\nhp $(2)$ }$(i)$\;\; $\fa a,b \in G\,(a^2b^2 = a^2$\: or\: $a^2b^2 = b^2)$.

\vspace{-0.2cm}
\tab{\nhp $(2)$ }$(ii)$\: Given $g,h \in \se X G$ such that $Z(g)\, \sub\, Z(h)$, there is
$h' \in \se X G$ such that $Z(h) = Z(h')$ and 

\vspace{-0.25cm}
\tab{\nhp $(2)$ $(ii)$\; }$g \spez h'$.

\vspace{-0.2cm}
\tab{\nhp $(2)$ }$(iii)$ For every saturated prime ideal $I$ of $G$, the quotient
reduced special group $(\se G I, \D {G_I})$

\vspace{-0.25cm}
\tab{\nhp $(2)$ $(iii)$ }is a RSG-fan.

\eth

\vspace{-0.4cm}

\nhp {\bf Remark.} The implication (1) $\Ra$ (2.iii) is Proposition \ref{fanQuot-2}.   
Condition (2.i) is assumption [Z] in the definition of a RS-fan (\ref{def:fan}); see also
Fact \ref{zeros} and Proposition \ref{charfan}. Therefore, to complete the proof of
\ref{thm:char-fan} we must only take care of (1) $\Ra$ (2.ii) and  (2) $\Ra$ (1), proved, 
respectively, in Propositions \ref{charfan:(1)implies(2.ii)}\,(2) and \ref{charfan:(2)implies(1)} 
below. \hfl $\Box$

\vspace{-0.4cm}

\bpr \label{charfan:(1)implies(2.ii)} Let $G$ be a RS-fan. Then: 

\vspace{-0.15cm}
\nhp $(1)$ For all elements $g,\,h \in \se X G$ such that $g \spez\, h$ $($hence
$Z(g) \, \sub \, Z(h))$ and every ideal \y I such that $Z(g) \, \sub \, I \,
\sub \, Z(h)$ there is $f \in \se X G$ such that $g \spez f \spez\, h$ and $Z(f) = I$.
\vspace{-0.1cm}

\nhp $(2)$ For every \, $g \in \se X G$ and every ideal \,$I \supseteq Z(g)$ there is 
a $($necessarily unique$)$ $f \in \se X G$  such that $g \spez f$ and
$Z(f) = I$. 
\vspace{-0.1cm}

\nhp $(3)$ For every ideal \y I of \y F there is an $f \in \se X G$ such that
$Z(f) = I$.
\epr

\vspace{-0.4cm}

\nhp {\bf Proof.} Since $G$ is a RS-fan, every TS-character $f:G\, \lra\, {\bf 3}$ 
is a RS-homomorphism. Thus, it suffices to construct TS-homomorphisms $f:G\, \lra\, {\bf 3}$ 
verifying (1) -- (3) of the statement.

\vspace{-0.15cm}

First we prove (1); the same proof, omitting item (c) below, also proves (2). 
Let $f: G \, \lra \; {\bf 3}$ be defined by:
\vspace{-0.25cm}

\nhp \hfl $f \, \lceil \, I = 0$ \;\; and \;\; $f \, \lceil \,(G \setminus I) =
g \, \lceil \,(G \setminus I)$. \hfl
\vspace{-0.25cm}

\nhp (a) $Z(f) = I$.
\vspace{-0.15cm}

\nhp By construction, $I \, \sub \; Z(f)$. Since $Z(g) \, \sub \; I$, $f(x) =
g(x) \neq 0$ for $x \in G \setminus I$, i.e., $Z(f) \, \sub \; I$.
\vspace{-0.15cm}

\nhp (b) $g \spez \, f$.
\vspace{-0.15cm}

\nhp Clear, from (a) and Lemma \ref{char-specializ}\,(4),
\vspace{-0.15cm}

\nhp (c) $f \spez \, h$.
\vspace{-0.15cm}

\nhp Clearly, $Z(f)\, \sub\, I\, \sub\, Z(h)$. If $h(a) \neq\, 0 $, then $a \not\in I$; since 
$g \spez \, h$, then $g(a) = h(a)$. Hence, $f(a) = g(a) = h(a)$, and we get $f \spez \, h$ by Lemma
\ref{char-specializ}\,(4).
\vspace{-0.15cm}

\nhp (d) \y f is a TS-homomorphism.
\vspace{-0.15cm}

\nhp Clearly $f(0) = 0$ and $f(\pm 1) = g(\pm 1) = \pm 1$. Let $a,b \in G$. 
If one of $a, b$ is in \y I, so is $ab$, and we have $f(a)f(b) = 0 = f(ab)$. 
If $a,b \not\in I$, then $ab \not\in I$, and \y f and \y g take the same 
value on \y a, \y b and $ab$; the result follows from the fact that 
\y g is a TS-character. Since \y G is a fan, \;$f \in \se X G$. 
\vspace{-0.15cm}

\nhp (3) This is Lemma 1.5, p. 102 (alternatively, Lemma 3.5, p. 114) in \cite{DP1}. \hfl $\Box$
\vspace{-0.35cm}

\bre \label{succ-1} The element \y f such that $g \spez  f$ and $Z(f) = I$ in 
\ref{charfan:(1)implies(2.ii)}\,(2) can also be obtained by taking any $h \in \se X G$ 
with $Z(h) = I$ (\ref{charfan:(1)implies(2.ii)}\,(3)) and setting $f = h^2g$.\hfl $\Box$
\ere

\vspace{-0.4cm}

The next Proposition proves the implication (2) $\Ra$ (1) in Theorem \ref{thm:char-fan}.  

\vspace{-0.45cm}
\bpr \label{charfan:(2)implies(1)} Let $G$ be a real semigroup verifying
conditions $(2.i) - (2.iii)$ of Theorem $\ref{thm:char-fan}$. Then, $G$ is a RS-fan.
\epr

\vspace{-0.45cm} 

\nhp {\bf Proof.} Item (2.i) is condition [Z] of the definition of RS-fan, \ref{def:fan}.  
It suffices to show that transversal representation in $G$ satisfies clause $[D^t]$ in 
\ref{fanrepr}, i.e., for $a,b,c \in G$:

\vspace{-0.1cm}

\nhp (I)\; $c \in \Dt G(a,b)$ and $Z(a) \subset Z(b)$ imply $c = a$.
\vspace{-0.15cm}

\nhp (II) $c \in \Dt G(a,b)$, $Z(a) = Z(b)$ and $a \neq -b$ imply $c = a$ or $c = b$.
\vspace{-0.1cm}

\nhp \und{Proof of (I)}. We first observe that the assumptions of (I) imply $Z(a) = Z(c)$. 

\vspace{-0.15cm}
Let $h \in \se X G$. If $h(a) = 0$, then $h(b) = 0$ (as $Z(a)\, \sub\, Z(b)$), and
$c \in \Dt G(a,b)$ yields $h(c) = 0$; hence, $Z(a)\, \sub\, Z(c)$. 

\vspace{-0.15cm}
If $Z(c)\, \sub\, Z(b)$, then $c \in \Dt G(a,b)$ yields $-a \in \Dt G(-c,b)$, and so
$Z(c)\, \sub\, Z(a)$. If $Z(b)\, \sub\, Z(c)$, then $-a \in \Dt G(-c,b)$ entails 
$Z(b)\, \sub\, Z(a)$, contrary to assumption. Hence, $Z(c)\, \sub\, Z(a)$, and
$Z(a) = Z(c)$.
\vspace{-0.15cm}

In order to prove $c = a$, let $h \in \se X G$. If $h(b) = 0$, then 
$h(c) \in \Dthr(h(a),0) = \{h(a)\}$, whence $h(c) = h(a)$. Henceforth, assume
$h(b) \neq 0$. Since $Z(a) \subset Z(b)$, there is $g \in \se X G$ so that
$g(b) = 0$ and $g(a) \neq 0$. Since the set of ideals of $G$ is totally ordered 
under inclusion, $h(b) \neq 0$ and $g(b) = 0$, we have $Z(h) \subset Z(g)$. By
(2.ii),\vspace{0.05cm} there is $g' \in \se X G$ so that $Z(g') = Z(g)$ and $h \spez g'$. 
Then, $g'(b) = 0$; from $c \in \Dt G(a,b)$ and $g(a) \neq 0$ comes $g'(a) = g'(c) \neq 0$.
From $h \spez g'$ we infer $h(a) = g'(a)$ and $h(c) = g'(c)$ (Lemma 
\ref{char-specializ}\,(4)), and from $g'(a) = g'(c)$ we conclude $h(a) = h(c)$,
and hence $a = c$.
\vspace{-0.1cm}

\nhp \und{Proof of (II)}. Assume $c \in \Dt G(a,b)$, $Z(a) = Z(b)$ and $a \neq -b$;
then, there is $g \in \se X G$ so that $g(b) = g(a) \neq 0$. First we claim:
\vspace{-0.1cm}

\nhp {\bf Claim 1.} Under the assumptions of (II), $Z(c) = Z(a) = Z(b)$.
\vspace{-0.15cm}

\nhp \und{Proof of Claim 1}. In fact, $c \in \Dt G(a,b)$ yields $Z(a) = Z(b)\, \sub\, Z(c)$.
Assume, towards a contradiction, that there is $h \in \se X G$ such that $h(c) = 0$
and $h(a) \neq 0$. From $c \in \Dt G(a,b)$ and $g(b) = g(a)$ we get 
$g(c) = g(b) = g(a) \neq 0$. Since the set of ideals of $G$ is totally ordered 
under inclusion, this and $h(c) = 0$ imply $Z(g) \subset Z(h)$. By (2.ii), there
is $h' \in \se X G$ such that $Z(h') = Z(h)$ and $g \spez h'$; it follows that
$h'(a) \neq 0$ and, since $Z(a) = Z(b)$, $h'(b) \neq 0$. Invoking Lemma 
\ref{char-specializ}\,(4), we get $h'(a) = g(a)$ and $h'(b) = g(b)$; from
$g(b) = g(a)$ we obtain $h'(b) = h'(a)$. On the other hand, $c \in \Dt G(a,b)$
and $h'(c) = h(c) = 0$ entail $h'(a) = -h'(b)$, whence $h'(a) = h'(b) = 0$,
contradiction. This proves $Z(c) = Z(a) = Z(b)$, as asserted.

If one of $a$ or $b$ is $0$, the equality of zero-sets in Claim 1 implies $c = a = b = 0$. 
So, assume, e.g., $b \neq 0$. Let $I$ be an ideal of $G$ ---necessarily prime and saturated---
maximal for $b \not\in I$. Let $\se {\sim} I$ be the congruence relation on $G$ determined by 
$I$, namely, for $x,y \in G$,
\vspace{-0.1cm}

\nhp \hfl $x\, \se {\sim} I\, y\;\; \Lra\;\; h(x) = h(y)$ for all $h \in \se X G$ 
such that $Z(h) = I$. \hfl \hfl (Cf. \S\,\ref{fan-quot}.C)

\nhp Note that the equality of zero-sets established in Claim 1, together with $b \not\in I$, 
implies that none of $a,b,c$ is in $I$.
\vspace{-0.1cm}

\nhp {\bf Claim 2.} $a\; \se {\not\sim} I\! -b$.
\vspace{-0.15cm}

\nhp \und{Proof of Claim 2}. Assume that $a\; \se {\sim} I\! -b$. Since $g(b) \neq 0$, i.e., 
$b \not\in Z(g)$, maximality of $I$ entails $Z(g)\, \sub\, I$. By (2.ii), there is $h \in \se X G$ 
such that $Z(h) = I$ and $g \spez h$. Since $h(b), h(a) \neq 0$, the specialization
$g \spez h$ yields $h(a) = g(a)$ and $h(b) = g(b)$ (\ref{char-specializ}\,(4)),
which, by $g(a) = g(b)$, entails $h(a) = h(b)$. On the other hand, $a\, \se {\sim} I\, -b$
and $Z(h) = I$ imply $h(a) = -h(b)$. Altogether, these equalities imply $h(a) = h(b) = 0$,
a contradiction, showing that $a\; \se {\not\sim} I\! -b$.
\vspace{-0.15cm}

By assumption (2.iii), the reduced special group $\se G I$ is a RSG-fan. With 
$\se {\pi} I:G\, \lra\, G/I$ denoting the canonical quotient map, we\vspace{0.05cm} have 
$\se {\pi} I(a) \neq \se {\pi} I(-b) = -\se {\pi} I(b)$. Note also that $a,b,c \not\in I$
implies $\se {\pi} I(a), \se {\pi} I(b), \se {\pi} I(c) \neq \se {\pi} I(0)$. From $c \in \Dt G(a,b)$ 
it follows $\se {\pi} I(c) \in \Dt {G/I}(\se {\pi} I(a), \se {\pi} I(b))$ which implies
(since $\se G I = (G/I) \setminus \{\se {\pi} I(0)\}$ is a RSG-fan) $\se {\pi} I(c) = \se {\pi} I(a)$ 
or $\se {\pi} I(c) = \se {\pi} I(b)$.
\vspace{-0.15cm}

\nhp {\bf Claim 3.} $\se {\pi} I(c) = \se {\pi} I(a)\;\; \Ra\;\; c = a$.
\vspace{-0.15cm}

\nhp \und{Proof of Claim 3}. Assumption $\se {\pi} I(c) = \se {\pi} I(a)$ means $c\; \se {\sim} I\, a$.
\vspace{-0.15cm}

Let $h \in \se X G$. Since the saturated prime ideals of $G$ are an inclusion chain, we consider 
two cases:
\vspace{-0.15cm}

\nhp --- $Z(h)\, \sub\, I$.

\vspace{-0.15cm}
\nhp Invoking (2.ii), let $h' \in \se X G$ be such that $Z(h') = I$ and $h \spez h'$. Since
$a,c \not\in I$, the specialization $h \spez h'$ entails $h(a) = h'(a)$ and $h(c) = h'(c)$; 
further, $c\; \se {\sim} I\, a$ gives $h'(c) = h'(a)$, whence $h(c) = h(a)$, for all $h \in \se X G$ 
such that $Z(h)\, \sub\, I$. 

\vspace{-0.15cm}
\nhp --- $Z(h)\, \supset\, I$.

\vspace{-0.15cm}
\nhp The maximality of $I$ implies $a,b,c \in Z(h)$, i.e., $h(a) = h(b) = h(c) = 0$.
\vspace{-0.15cm}

\nhp These two cases prove that $h(c) = h(a)$, for all $h \in \se X G$, i.e., $c = a$.
\vspace{-0.15cm}

A similar argument proves that $\se {\pi} I(c) = \se {\pi} I(b)\;\; \Ra\;\; c = b$, completing
the proof of (II), of Proposition \ref{charfan:(2)implies(1)}, and of Theorem \ref{thm:char-fan}. 
\hfl $\Box$
\vspace{-0.35cm}

\bre \label{ChainLength} (Chain length) There is a well-known characterization of fans in the
categories {\bf AOS} and {\bf RSG} in terms of {\it chain length}, i.e., the size of longest
strict inclusion chain of non-empty subbasic opens $U(a)$, cf. proof of \ref{qfandense}:  an AOS 
is a fan if and only if its chain length is $\leq 2$, see \cite{ABR}, Prop. 3.11, p. 74, or \cite{M},
Thm. 4.2.1\,(2), p. 65. This notion of chain length also makes sense for ARSs, cf. \cite{M}, 
p. 167. However, this characterization is no longer valid for ARSs or RSs; an easy computation
shows that the RS-fan $\se F 2$ in Example \ref{ex-fans}.\,B.(i), see Figure 4, has chain length 4.
With $\Spec (F)$ denoting the set of ideals of $F$, the integer $2\, \cdot\,$card$\,(\Spec (F))$ is 
an upper bound on the chain length of a RS-fan, $F$, with a finite spectrum; this is easily proved 
using Theorem \ref{thm:char-fan}; see also \cite{M}, Thm. 8.5.3, p. 167.\hfl $\Box$   
\ere
\vspace{-0.4cm}

The next two corollaries of Theorem \ref{thm:char-fan} give stylized (abstract) versions of the 
notion of a {\it trivial fan}, a basic concept in the theory of (pre-)orders on fields 
(see \cite{La}, Prop. 5.3, p. 39). Their translation in the case of preordered rings is given in 
Theorem \ref{TotPreordsAndFans} below, where it will be obvious that in the case of fields they 
boil down to the notion of a trivial fan.

\vspace{-0.4cm}
\bco \label{TotOrder=Fan} Let $G$ be a real semigroup such that the character space $\se X G$
is totally ordered under specialization. Then, $G$ is a RS-fan.
\eco

\vspace{-0.45cm}
 
\nhp {\bf Proof.} We check that conditions (2.i) -- (2.iii) of Theorem \ref{thm:char-fan} hold.
\vspace{-0.2cm}

Since every saturated prime ideal of $G$ is the zero-set of some character (\cite{DP1}, Lemma
3.5, p. 114) and $g \spez h\; \Ra$ $Z(g)\, \sub\, Z(h)$ for $g,h \in \se X G$ (\ref{char-specializ}\,(4)), 
the set of saturated prime ideals of $G$ is an inclusion chain, i.e., item (2.i) of \ref{thm:char-fan} 
holds.
\vspace{-0.2cm}

Further, every saturated prime ideal is the zero-set of a {\it unique} character: if 
$\se h 1, \se h 2 \in \se X G$ are such that $Z(\se h 1) = Z(\se h 2)$, then $\se {h^2} 1 = 
\se {h^2} 2$ (\ref{char-zeroset}\,(1)); if, say, $\se h 1 \spez\, \se h 2$, by Lemma 
\ref{char-specializ}\,(5), $\se h 2 = \se {h^2} 2\, \se h 1 = \se {h^2} 1\, \se h 1 = \se h 1$.
It follows that, for every saturated prime ideal $I$ the quotient $G/I$ has a unique character,
and hence $G/I\, \cong\, {\bf 3}$, which is a RS-fan, showing that condition (2.iii) of 
\ref{thm:char-fan} holds.
\vspace{-0.2cm}

Finally, to check item \ref{thm:char-fan}\,(2.ii), observe that the linearity assumption and 
the uniqueness proved in the preceding paragraph yield $Z(g)\, \sub\, Z(h)\; \Ra\; g \spez h$. 
\hfl $\Box$

\noi {\bf Remark.} We refer the reader to \cite{M}, Prop. 8.8.4, pp. 178-179, where he proves 
that the RSs in \ref{TotOrder=Fan} are spectral. Conversely, it is an easy exercise to prove that 
the RSs which are simultaneously spectral and fans are exactly those whose character space 
is totally ordered by specialization. \hfl $\Box$

\vspace{-0.35cm}
\bco \label{2chains=Fan} Let $G$ be a real semigroup satisfying the following requirements:
\vspace{-0.2cm}

\nhp $(1)$\; Condition $[$Z\,$]$ in $\ref{fanrepr}$.
\vspace{-0.2cm}

\nhp $(2)$\; The character space $\se X G$ of\, $G$ is the union of two specialization
chains, $\cC_0, \cC_1$.
\vspace{-0.2cm}

\nhp $(3)$\; For every saturated prime ideal $I$ of $G$ and for $i = 0,1$, there is $h_i \in \cC_i$
such that $Z(h_i) =$

\vspace{-0.23cm}

\tab{\nhp $(3)$\; }$= I$.
\vspace{-0.2cm}

\nhp Then, $G$ is a RS-fan.
\eco

\vspace{-0.45cm}

\nhp {\bf Remark.} The specialization chains in item (2) may not be disjoint, and the characters
$\se h 0, \se h 1$ in (3) may be identical.  Using condition (3) (and \ref{thm:char-fan}\,(2.ii)) it 
can be shown that the chains in (2) are maximal.
\vspace{-0.15cm}

\nhp {\bf Proof.} Again, we check that conditions (2.i) -- (2.iii) of Theorem \ref{thm:char-fan} hold.
Condition (1) is item \ref{thm:char-fan}\,(2.i).
\vspace{-0.15cm}

\nhp (2.ii) Let $g \in \se X G$ and let $I$ be a saturated prime ideal of $G$ such that
$Z(g)\, \sub\, I$. Condition (2) implies that either $g \in \cC_0$ or $g \in \cC_1$, say the
first. By (3) there is $\se h 0 \in \cC_0$ so that $Z(\se h 0) = I$. Since $\cC_0$ is a specialization 
chain and $g, \se h 0 \in \cC_0$ ((3)), the inclusion $Z(g)\, \sub\, I = Z(\se h 0)$ yields 
$g \spez \se h 0$, proving (2.ii).
\vspace{-0.15cm}

\nhp (2.iii) For every saturated prime ideal $I$ of $G$, the structure $\se G I = (G/I) \setminus 
\{\se {\pi} I(0)\}$, with representation induced by $\D G$\,, is a RSG (\cite{DP4}, Thm. II.3.15\,(d)) 
and $\se X {G_I} = \{h \in \se X G \sth Z(h) = I\}$ which, by 
assumptions (2) and (3), equals $\{\se h 0, \se h 1\}$. Since every reduced special group with 
at most two characters is a RSG-fan, so is  $\se G I$, as required. \hfl $\Box$
\vspace{-0.1cm}

\noi {\bf Remarks.} (a) There are examples satisfying conditions (2.i) and (2.iii) of Theorem 
\ref{thm:char-fan} but not condition (2.ii).
\vspace{-0.15cm}

\noi (b) The real semigroup $\se G {C(X)}$ associated to the ring $C(X)$ of continuous, 
real-valued functions on a topological space $X$ satisfies conditions (2.ii) and (2.iii) of 
Theorem \ref{thm:char-fan} but,\vspace{0.05cm} in general, not (2.i); cf. \cite{M}, 5.2\,(6), 
p. 87. 
\vspace{-0.15cm}

\noi (c) Even in the presence of conditions (1) and (3) of \ref{2chains=Fan}, if the character
space $\se X G$ of $G$ is the union of \und{more than two} maximal specialization chains, the situation 
becomes more involved, as illustrated by the examples in \cite{DP5b}, 2.18.  \hfl $\Box$

\vspace{-0.4cm}

\section{Fans and preordered rings} \label{fan-prings}

\vspace{-0.35cm}

In this section we prove a number of results about, and exhibit some 
examples of semi-real rings and preordered rings $($hereafter {\bf p-ring}s$)$ whose
associated real semigroups are fans.

\vspace{-0.5cm}

\paragraph{I. Properties of p-rings\vspace{-0.15cm} whose associated real semigroup is a fan.}

\

\nhp Throughout this subsection we assume that $\fm{A,T}$ is a p-ring. 
\vspace{-0.1cm}

\noi {\bf A. Basic correspondences.} Let $\fm{A,T}$ be a p-ring and let
$\se G {\!A,T}$ denote its associated real semigroup (\cite{DP2}, \S\,1, p. 51 or
\cite{DP3}, 9.1\,(A), pp. 406-407). Let $\mathrm{Sat}(G)$ denote the set of all saturated 
ideals of a real semigroup $G$. 

\vspace{-0.2cm}

\noi For an ideal $I$ of $A$, let ${\cl I} = \{\cl a \sth a \in I\}$, and for 
$J \in \mathrm{Sat}(\se G {\!A,T})$, set $\h J := \{a \in A \sth \cl a \in J\}$. The following facts 
are easily verified or their proofs briefly indicated:
\vspace{-0.4cm}

\bfa \label{FactsBasicCorresp} With notation as above and $J, \se J 1, \se J 2 \in \mathrm{Sat}(\se G {\!A,T})$, we have:
\vspace{-0.3cm}

\noi \;\:$(i)$\; $\cl I$ is a saturated ideal of $\se G {\!A,T}$. \hmm $(ii)$\; $\h J$ is an ideal of $A$. \hmm $(iii)$\; $\cl{\h J} = J$.
\vspace{-0.3cm}

\noi $(iv)$\; \y J prime\; $\Lra$\; $\h J$ prime. \hem $(v)$\; $\se J 1\, \sub\, \se J 2\;\; \Lra\;\; \h {\se J 1}\, \sub\, \h {\se J 2}$.
 \hem\, $(vi)$ The map\, $J\, \longmapsto\, \h J$ is injective.
\efa
\vspace{-0.45cm}

\noi {\bf Proof.} We only prove saturatedness in (i). Let $a,b \in I$ and $c \in A$ be\vspace{0.05cm} so that\, 
$\cl{c} \in \se D G(\cl{a},\cl{b})$. By \cite{M}, Prop. 5.5.1\,(5), p. 95, there are\, 
$\se t 0, \se t 1, \se t 2 \in T$\, so that\, $\se t 0 c = \se t 1 a + \se t 2 b$\, and\, 
$\cl{\se t 0 c} = \cl{c}$. From $a,b \in I$ follows $\se t 0 c \in I$, whence\, 
$\cl{c} =\, \cl{\se t 0 c} \in \cl{I}$. \hfl $\Box$

\noi The following notions from real algebra are used in the sequel:
\vspace{-0.35cm}

\bdf \label{ConvexIdeal} Given a (proper) preorder $T$ of $A$, an ideal $I\, \sub\, A$ is
\vspace{-0.15cm}

\noi (i)\;\; $T${\bf -radical} iff for all $a \in A$ and $t \in T$,\; $a^2 + t \in I\; \Ra\; a \in I$
(and hence $t \in I$).
\vspace{-0.15cm}

\noi (ii)\; $T${\bf -convex} iff for $\se t 1, \se t 2 \in T$,\; $\se t 1 + \se t 2 \in I\; \Ra\;
\se t 1, \se t 2 \in I$. \hfl $\Box$
\edf
\vspace{-0.4cm}

\noi Remark that an ideal is $T$-radical iff it is $T$-convex and radical (\cite{BCR}, Prop. 4.2.5, 
p. 87). We denote by ${\mathrm{PConv}}\,(A,T)$ the set of all $T$-convex prime (equivalently, $T$-radical
prime) ideals of $A$. For further properties of $T$-convexity (e.g., the definition of the $T$-radical 
of an ideal, $\sqrt[T]{I}$), the reader is referred to \cite{BCR}, \S\S\,4.2, 4.3. We prove:

\vspace{-0.35cm}

\bpr \label{conv-satthm} Let \y J be a saturated ideal of\, $G = \se G {\!A,T}$. 
Then $\h J$ is a \y T-radical ideal of \y A. 
\epr

\vspace{-0.65cm}

\nhp {\bf Proof.} Assume $a^2 +\, t \in \h J$, where $a \in A,\, t \in T$; 
we must show that $a \in \h J$. Write $j = a^2 +\, t$; then $\cl{j} \in J$ (definition of\; 
$\h {}$\;\,), and also $\cl{j} \in \se {D^t} G(\cl{a}\,^2, \cl{t})$ 
(cf. \cite{M}, p. 96). Recall that $\se X G = \Sper (A,T)$.
\vspace{-0.2cm}

Let $\al \in \Sper (A,T)$ be such that $\cl{j}(\al) = 0$; then, 
$\cl{a}\,^2(\al) = -\cl{t}(\al)$. Since $t \in T\, \sub\, \al$, we have 
$-\cl{t}(\al) \in \{0, -1\}$. On the other hand, $\cl{a}\,^2(\al) \in \{0, 1\}$, 
since $\cl{a}\,^2$ is a square. Thus, the equality above forces 
$\cl{a}(\al) = \cl{t}(\al) = 0$, proving that $Z(\cl{j})\, \sub\, 
Z(\cl{a})\, \cap\, Z(\cl{t})\, \sub\, Z(\cl{a})$. This inclusion is equivalent 
to $\cl{a}\,^2 = \cl{a}\,^2 \cdot \cl{j}^2$ (see \ref{char-zeroset}\,(2)). 
Then, $\cl{a}\,^2 \in J$, whence $\cl{a} \in J$, which proves $a \in \h J$. \hfl $\Box$

\vspace{-0.15cm}
We register the following consequences:
\vspace{-0.35cm}

\bco \label{fan-biject} For any ideal $I$ of $A$ we have:
\vspace{-0.2cm}

\noi $(i)$\; $\h {\cl I}$ is the smallest $T$-radical ideal containing $I$, i.e., 
$\h {\cl I} = \sqrt[T]{I}$.
\vspace{-0.2cm}

\noi $(ii)$\, $\cl{\sqrt[T]{I}}\; = \cl {\h {\cl{I}}} =\; \cl{I\,}$. \hfl $\Box$
\eco
\vspace{-0.45cm}

\noi {\bf Notation.} Given a p-ring $\fm{A,T}$, we denote by $\Spec(\se G {\!A,T})$  the set
of all ideals of the real semigroup $\se G {\!A,T}$. If $\se G {\!A,T}$ is a RS-fan, we know 
(\ref{charfan}\,(5)) that $\Spec(\se G {\!A,T})$ is totally ordered under inclusion. By 
 \ref{FactsBasicCorresp}\,(v) the set $\{\h J \sth J \in \Spec(\se G {\!A,T})\}$ of ideals of $A$ is totally ordered 
 under inclusion as well. \hfl $\Box$

\vspace{-0.45cm}
\bfa \label{hatfacts2} If\, $\se G {\!A,T}$ is a fan, then:

\vspace{-0.2cm}

\nhp $(i)$\; Every $T$-radical ideal of $A$ is prime.
\vspace{-0.15cm}

\nhp $(ii)$ $\{\h J \sth J \in \Spec(\se G {\!A,T})\} = {\mathrm{PConv}}\,(A,T)$.

\vspace{-0.2cm}
\nhp $(iii)$ The map $J \mapsto \h J$ is an order-preserving bijection from $\Spec(\se G {\!A,T})$ 
onto ${\mathrm{PConv}}\,(A,T)$.
\efa

\vspace{-0.55cm}
\nhp {\bf Proof.} (i) Let $I$ be a $T$-radical ideal of $A$. By \ref{fan-biject}\,(i), $I = \h {\cl{I}}$ and, 
\vspace{-0.08cm}by \ref{FactsBasicCorresp}\,(i), $\cl I \in \Spec(\se G {\!A,T})$. Then, \ref{charfan}\,(4) implies 
$\cl I$ prime, which, by \ref{FactsBasicCorresp}\,(iv), in turn yields $I = \h{\cl I}$ prime.
\vspace{-0.2cm}

\nhp (ii) Let $I \in {\mathrm{PConv}}\,(A,T)$. Then $I$ is a $T$-radical ideal of $A$, and by 
\ref{fan-biject}\,(i) we have $I = \h{\cl I}\; (= \sqrt[T]{I})$.

\vspace{-0.15cm}

\nhp (iii) Follows from (i) and \ref{FactsBasicCorresp}\,(vi, v). \hfl $\Box$

\vspace{-0.4cm}
\bres \label{RmksT-convPrimes} Assume, as above, that $\se G {\!A,T}$ is a RS-fan.
\vspace{-0.4cm}

\noi (i) $\cl {(0)}$ is an ideal of $\se G {\!A,T}$ (\ref{FactsBasicCorresp}\,(i)), hence 
prime\vspace{0.05cm} and saturated (\ref{charfan}\,(4); 
\ref{fanideals}); by \ref{FactsBasicCorresp}\,(iv),
$\h {\cl {(0)}}$ is prime, and by \ref{fan-biject}\,(i), $\h {\cl {(0)}} = \sqrt[T]{(0)} =
\{a \in A \sth -a^{2k} \in T \mathrm{\; for\; some\; } k \geq 0\}$ is the smallest element 
of ${\mathrm{PConv}}\,(A,T)$.
\vspace{-0.2cm}

\nhp (ii) The maximal element of $\Spec(\se G {\!A,T})$ is:

\vspace{-0.2cm}
\nhp \hfl $\fM$ = set of non-invertible elements of $\se G {\!A,T}$ = $\{\cl x \sth x \in A
\mathrm{\; and\;\, } {\cl x}^2 \neq 1\}.$ \hfl

\vspace{-0.2cm}
\nhp Then, the ideal

\vspace{-0.5cm}
\begin{tabbing}
\hspace*{2.5cm}$M= \h{\fM}\;$\= $= \{a \in A \sth \cl a \in \fM\} = \{a \in A \sth {\cl a}^2 \neq 1\} =$\\
\> $= \{a \in A \sth \ex \al \in \Sper(A,T) \mathrm{\; such\; that\;\, } {\cl a}(\al) = 0\} =$\\

\> $= \{a \in A \sth \ex \al \in \Sper(A,T) \mathrm{\; such\; that\;\, } a \in \supp(\al)\} =$\\
\> $= \,\bigcup\, \{\supp(\al) \sth \al \in \Sper(A,T)\},$
\end{tabbing}

\vspace{-0.4cm}

\noi is the maximal element of ${\mathrm{PConv}}\,(A,T)$.
\vspace{-0.15cm}

\nhp (iii) \und{Warning}. Even though the\vspace{0.08cm} ideal $M$ is maximal \und{in ${\mathrm{PConv}}\,(A,T)$}, 
it {\it may not} be a maximal ideal of $A$ (e.g., $\fm{\Zh, \sum {\Zh}^2})$; however, it is maximal in some important cases, e.g., when 
$\fm{A,T}$ is a {\it bounded inversion ring}; cf. \cite{DM2}, Prop. 7.2, p. 78. \hfl $\Box$
\eres

\vspace{-0.75cm}
\bct \label{TS-charsInP-rings} {\bf Ternary semigroup characters of  $\se G {\!A,T}$\,.}
The characterization of ARS-fans given in Proposition \ref{equivfan}\,(1), can be restated as follows:

\vspace{-0.8cm}
   
$$
(\dag) \;\; \left\{
\begin{array}{ll}
\textrm{A real semigroup $G$ is a RS-fan if and only if the set of its prime ideals is} \\ 
\textrm{totally ordered under inclusion and every character of {\it ternary} semigroup} \\ 
G\, \lra\, {\bf 3} \textrm{\; preserves representation.} 
\end{array}
\right.
$$
\ect

\vspace{-0.45cm}

In the case $G = \se G {\!A,T}$ we register, without proof, a description of the
TS-characters of $\se G {\!A,T}$ in terms of the p-ring $\fm{A,T}$.

\vspace{-0.35cm}

\bpr \label{TS-charInP-rings} Let $\fm{A,T}$ be a p-ring. There is a bijective correspondence 
between the TS-characters of $\se G {\!A,T}$ and the family of all subsets $S \, \sub\, A$
satisfying the following conditions:
\vspace{-0.2cm}

\noi $(1)\;\; T\, \sub\, S$. \;\; $(2)\;\; S$ is closed under product. \;\; $(3)\; -1 \not\in S$. \;\;
$(4)\;\; S \cup -S = A$. 
\vspace{-0.15cm}

\noi $(5)$ The set $S \cap -S$ $($not necessarily an ideal!\,$)$ is prime: for all 
$x,y \in A,\; xy \in S \cap -S$\, implies\, \hspace*{0.55cm} $x \in S \cap -S$ or $y \in S \cap -S$. 
\vspace{-0.15cm}

\noi $(6)\;\, S \cap -S$ is $T$-convex: for all $\se t 1, \se t 2 \in T$, if 
$\se t 1 + \se t 2 \in S \cap -S$, then $\se t 1, \se t 2 \in S \cap -S$. \hfl $\Box$
\epr
\vspace{-0.35cm}

By this characterization, the (necessary and sufficient) condition for $\se G {\!A,T}$ to be a
RS-fan given by \ref{TS-charsInP-rings}\,(\dag) translates as the conjunction of:
\vspace{-0.45cm}

\bct\hspace*{-0.1cm}{\bf (i)}\label{Cond1G_(A,T)Fan} The set ${\mathrm{PConv}}\,(A,T)$ of $T$-convex prime ideals of $A$
is totally ordered under inclusion, 
\vspace{-0.25cm}

\noi and
\vspace{-0.15cm}

\noi {\bf \ref{Cond1G_(A,T)Fan}\,(ii)}\label{Cond2G_(A,T)Fan} Every subset $S\, \sub\, A$ satisfying conditions 
(1) -- (6) of \ref{TS-charInP-rings} is closed under addition (i.e., is an element of $\Sper(A, T)$).
\ect
\vspace{-0.45cm}

However,  we show that in the present case we can dispense with \ref{Cond1G_(A,T)Fan}\,(i):
\vspace{-0.4cm}

\bpr \label{T-convPrimesTotOrd}  Let $\fm{A,T}$ be a p-ring. With notation as above, 
\vspace{-0.15cm}

\noi $(1)$ Condition $\ref{Cond2G_(A,T)Fan}\,(\mathrm{ii})$ implies $\ref{Cond1G_(A,T)Fan}\,(\mathrm i)$.
Hence, 
\vspace{-0.15cm}

\noi $(2)$ $\se G {\!A,T}$ is a RS-fan if and only if every subset $S\, \sub\, A$ satisfying conditions 
$(1) - (6)$ of Proposition {\em{\ref{TS-charInP-rings}}} is closed under addition.
\epr
\vspace{-0.45cm}

\nhp {\bf Proof.} We need only prove (1). Let $I,J \in {\mathrm{PConv}}\,(A,T)$; let $\al \in \Sper(A,T)$ be such that 
$I = \supp(\al)$ (cf. \cite{BCR}, Prop. 4.3.8, p. 90). Set $S = J \cup \al$.

\vspace{-0.2cm}
We first observe that $S$ satisfies conditions (1) -- (6) of \ref{TS-charInP-rings}. Conditions (1) -- (3) are obvious.

\vspace{-0.15cm}
\nhp (4). Since $-S = J \cup -\al$ and 
$\al \cup -\al = A$, we have $S \cup -S = J \cup \al \cup -\al = A$.

\vspace{-0.15cm}
\nhp (5). From the previous item we have  
$S \cap -S = (J \cup \al) \cap (J \cup -\al) = J \cup (\al \cap -\al) = J \cup I$. Since both 
$I,J$ are prime (ideals), we get $xy \in S \cap -S$ implies $x \in S \cap -S$ or $y \in S \cap -S$.

\vspace{-0.15cm}
\nhp (6). Again, since $S \cap -S = J \cup I$ and both $I,J$ are $T$-convex, we get the desired conclusion.

\vspace{-0.15cm}
By assumption, $S$ is additively closed. Assume, towards a contradiction, that there are $a,b \in A$
such that $a \in I \setminus J$ and $b \in J \setminus I$. In particular, $a \in I\, \sub\, \al\, \sub\, S$ 
and $b \in J\, \sub\, S$, whence $a + b \in S$. If $a + b \in J$, since $-b \in J$ we get 
$a = (a + b) + (-b) \in J$, contradiction. Then, $a + b \in \al$, and from $a \in I\; \sub\! -\al$,
we get $b = (a + b) + (-a) \in \al$. Next, since $-a \in I \setminus J$ and $-b \in J \setminus I$, 
the preceding argument can be carried out with $a,b$ replaced with $-a,-b$, respectively,
to conclude that $-b \in \al$. Thus, $b \in \al \cap -\al = I$, contradiction. \hfl $\Box$ 

\vspace{-0.4cm}
\bre \label{rmk:multSemiOrd+Fans} Proposition \ref{T-convPrimesTotOrd}\,(2) is {\it the exact} 
ring-theoretic analog of the definition of a fan (as a preorder) in a field, due to \cite{BK}; namely:
{\it A preorder $T$ of a field $F$ is a fan iff for any set $S \supseteq T$ such that $-1 \not\in S$ and 
$S^{\times} = S \setminus \{0\}$ is a subgroup of index 2 in $F^{\times}$, then $S$ is closed
under addition} (\cite{La}, Def. 5.1, p. 39).  \hfl $\Box$
\vspace{-0.25cm}

\medskip
The characterization in \ref{T-convPrimesTotOrd}\,(2) yields a first batch of natural examples 
of p-rings whose associated real semigroup is a fan.

\vspace{-0.5cm}
\bco \label{corol:ValRingsOfFans} Let $K$ be a field and $T$ be a preorder of $K$ which is
a fan. Let $A$ be a subring of $K$ whose field of fractions is $K$. Then, the real semigroup 
$\se G {\!A,T\, \cap\, A}$ is a fan. In particular, if $A = A_v$ is the valuation ring of a
$T$-compatible valuation $v$ of $K$, the real semigroup $\se G {\!A_v,T\, \cap\, A_v}$ is a fan.
\eco

\vspace{-0.5cm}
\nhp{\bf Proof.} According to Proposition \ref{T-convPrimesTotOrd}\,(2) it suffices to check
that any set $S \, \sub\, A)$ satisfying conditions \ref{TS-charInP-rings}\,(1) -- (6) in $\fm{A, T \cap T}$
is closed under addition. Let $S' = \{\fr a b \sth a,b \in A, b \neq 0 \mbox{\; and\; } ab \in S\}\, \sub\, K$.
We first show:

\vspace{-0.15cm}
\nhp --- $S' \setminus \{0\}$ is a subgroup of $K^{\times}$, $T\, \sub\, S'$ and $-1 \not\in S'$.

\vspace{-0.15cm}
\nhp Clearly, $S\, \sub\, S'$ and (by (3)) $-1 \not\in S'$. Since $K$ is the field of fractions of $A$, 
any element of $T$ can be written as $\fr a b$, with $a,b \in A, b \neq 0$. Then, 
$ab = \fr a b \cdot b^2 \in T \cap A$. Since $T \cap A\, \sub\, S$ ((1)), we get $\fr a b \in S'$, 
hence $T\, \sub\, S'$. Condition (2) implies that $S' \setminus \{0\}$ is a subgroup of $K^{\times}$.

\vspace{-0.15cm}
Since, by assumption, $T$ is a fan in the field $K$, the set $S'$ is closed under addition in $K$ (see
\ref{rmk:multSemiOrd+Fans}), which clearly implies that $S$ is additively closed in $A$. \hfl $\Box$

\vspace{-0.5cm}

\paragraph{B. Total preorders and trivial fans in rings.\vspace{-0.35cm}}

\bno \label{notat:Localizations} For a p-ring $\fm {A,T}$ and a prime ideal $I$ of $A$, we let 

\vspace{-0.35cm}
\begin{itemize}
\vspace{-0.2cm}
 \item $\se A I$ denote the localization of $A$ at $I$,
\vspace{-0.2cm}
 \item $\se M I = I \cdot\, \se A I$ denote the maximal ideal of $\se A I$, and
\vspace{-0.2cm}
 \item $\se T I = T \cdot (A \setminus I)^{-2}$ denote the preorder induced by $T$
 on $\se A I$. \hfl $\Box$
\end{itemize}
\eno

\vspace{-0.8cm}
\bfa \label{QuotPreord} If $I \in {\mathrm{PConv}}\,(A,T)$, then $\se T I/\se M I$ is a proper preorder 
of the field $\se A I/\se M I$ $($and $\se T I$ a proper preorder of $\se A I)$.
\efa

\vspace{-0.5cm}
\nhp {\bf Proof.} It suffices to prove the first assertion. Clearly $\se T I/\se M I$ is a preorder 
of $\se A I/\se M I$. To show it is proper, assume, on the contrary, that 
$-1 \in \se T I/\se M I$, i.e., $-1 = \left ( {\fr{t}{x^2}} \right )/\se M I$\,,
with $t \in T$ and $x \in A \setminus I$; that is, ${\fr{t}{x^2}} + 1 \in \se M I = I \cdot \se A I$, i.e.,
${\fr{t + x^2}{x^2}} = \fr{i}{y}$, for some $i \in I$ and $y \in A \setminus I$. Since $I$ is prime,
we get $y \cdot (t + x^2) \equiv x^2 i$ (mod $I$), whence $y \cdot (t + x^2) \in I$, and $t + x^2 \in I$. 
Since $t, x^2 \in T$ and $I$ is $T$-convex and radical, we obtain $x \in I$, contradiction. \hfl $\Box$

\vspace{-0.45cm}

\bdf \label{def:totalPreorder} A {\bf total preorder} in a ring $A$ is a (proper) preorder $T$ 
such that $T \cup -T$ $= A$. \hfl $\Box$ 
\edf

\vspace{-0.85cm}

\bfa \label{SuppTotPreordIsIdeal} For a total preorder $T$ of a ring $A$, $T \cap -T$ is a
 proper $T$-convex ideal of $A$. Any $T$-convex ideal of $A$ contains $T \cap -T$. \hfl $\Box$
\efa

\vspace{-0.75cm}

\bres \label{rmks:SuppTotPreord} (i) The ideal $T \cap -T$ {\it may not} be prime (see Example
\ref{ex:Tcap-TNotPrime}). When it is, the notion of ``total preorder''coincides with 
``prime cone'', i.e., element of $\Sper(A)$.

\vspace{-0.15cm}
\nhp (ii) When $T \cap -T = \{0\}$ the total preorders are just the total orders of $A$. \hfl $\Box$
\eres

\vspace{-0.75cm}
\bex \label{ex:Tcap-TNotPrime} Let $A := \re[X]/(X^2)$; the elements of $A$ are uniquely
representable in the form $aX + b$ with $a,b \in \re$. Clearly, the zero ideal of $A$ is
not radical, hence not prime either: $X \neq 0$ but $X^2 = 0$. We define a total (pre)order 
$T$ in $A$ by the stipulation:

\vspace{-0.15cm}
\nhp \hfl $aX + b \in T$\; iff\; $b > 0$\; or\; ($b = 0$\, and\, $a \geq 0$). \hfl

\vspace{-0.15cm}

\nhp Checking that $T$ is a total (pre)order of $A$ is routine, left to the reader.
However, the ideal $T \cap -T = \{0\}$ is not prime (not even radical). \hfl $\Box$
\eex

\vspace{-0.4cm}
Proposition \ref{LiftTotPreord} shows that total preorders are preserved by localization at, 
and lifting by convex prime ideals.

\vspace{-0.4cm}
\bpr \label{LiftTotPreord} Let $I$ be a prime\vspace*{0.05cm} ideal of a ring $A$, $T$ be a preorder of $A$, and $Q$ 
be a preorder of the localisation $\se A I$. Let $\se {\iota} I: A\, \lra\, \se A I$ be the canonical map 
$a \mapsto \fr a 1\;\; (a \in A)$. Then,
\vspace{-0.2cm}

\noi $(i)$\;\;\; If $T$ is a total preorder, then so is $\se T I$.
\vspace{-0.15cm}

\nhp $(ii)$\;\; $P := {\se {{\iota}^{-1}} I}[Q]$ is a preorder of $A$.
\vspace{-0.15cm}

\nhp $(iii)$\, $\se P I = Q$.
\vspace{-0.15cm}

\nhp $(iv)$\; $Q$ is total if and only if $P$ is total.
\vspace{-0.15cm}

\nhp $(v)$\;\, The maximal ideal $\se M I$ of $\se A I$ is $Q$-convex if and only if $I$ is $P$-convex.
\epr

\vspace{-0.4cm}
\nhp {\bf Proof.} (i), (ii) and (iv) are straightforward checking.

\vspace{-0.15cm}

\noi (iii) We show:
\vspace{-0.35cm}

\nhp --- $\se P I\, \sub\, Q$\,.\; Let $z \in \se P I$, i.e., $z = \fr p {x^2}$, with 
$p \in P, x \not\in I$. Then, $\fr p 1 \in Q$, $\fr {x^2} 1$ is invertible in $\se A I$, 
and $\fr 1 {x^2} \in Q$. It follows that $z = \fr p 1 \cdot \fr 1 {x^2} \in Q$\,.

\vspace{-0.25cm}
\nhp --- $Q\, \sub\, \se P I$\,.\; Let $z \in Q$; then, $z = \fr x y$, with $x,y \in A, y \not\in I$, 
which\vspace{-0.1cm} implies $z = \fr {xy} {y^2}$; this gives ${\fr {y^2} 1} \cdot z = \fr {xy} 1$. 
Clearly, $\fr {y^2} 1 = \left ( \fr y 1 \right )^2 \in Q$, whence ${\fr {y^2} 1} \cdot z \in Q$,
and $\fr {xy} 1 \in Q$, which shows that $xy \in P$. Hence, $z = \fr {xy} {y^2} \in \se P I$.
\vspace{-0.3cm}

\nhp (v) ($\Ra$). Let $\se p 1, \se p 2 \in P$ be such that $\se p 1 + \se p 2 \in I$. Then, 
$\fr {p_i} 1 \in Q\; (i = 1,2)$, and\; $\fr {p_1 + p_2} 1 \in I \cdot \se A I =
\se M I$. By the convexity assumption, $\fr {p_1} 1, \fr {p_2} 1 \in I \cdot \se A I$.
For $i = 1,2$, we have $\fr {p_i} 1 = \fr j x$, with $j \in I, x \not\in I$. It follows
that $x p_i - j \in I$; since $x \not\in I$, we get $p_i \in I$, as required. 
\vspace{-0.1cm}

\noi ($\La$). Let $\fr {x_i} {y_i} \in Q\; (x_i \in A, y_i \not\in I; i = 1,2)$ be such that 
$\fr {x_1} {y_1} + \fr {x_2} {y_2} = \fr {x_1 y_2 + x_2 y_1} {y_1 y_2} \in \se M I$.
Then, $ \fr {x_1 y_2 + x_2 y_1} {y_1 y_2} = \fr z w$, with\vspace*{0.05cm} $z \in I, w \not\in I$. 
We get\: $w (x_1 y_2 + x_2 y_1) = z y_1 y_2 \in I$; since $w \not\in I$, we have\: 
$x_1 y_2 + x_2 y_1 \in I$. By (iii)\vspace*{0.05cm} we get $\fr {x_i} {y_i} = \fr {p_i} {s^2_i}$, with 
$p_i \in P, s_i \not\in I$, whence 
\vspace{-0.35cm}

\noi (\dag)\; $x_i s^2_i = y_i p_i\;\; (i = 1,2)$.
\vspace{-0.1cm}

\noi Scaling\: $x_1 s^2_1 = y_1 p_1$\: by\: $y_2 s^2_2$\: yields\: 
$x_1 y_2 s^2_1 s^2_2 = y_1 y_2 p_1 s^2_2$. Likewise, we obtain
$x_2 y_1 s^2_1 s^2_2 = y_1 y_2 p_2 s^2_1$. Adding these\vspace*{0.05cm} terms gives\: 
$s^2_1 s^2_2 (x_1 y_2 + x_2 y_1) = y_1 y_2 (p_1 s^2_2 + p_2 s^2_1)$. Since 
$x_1 y_2 + x_2 y_1 \in I$ and $y_1 y_2 \not\in I$, we get $p_1 s^2_2 + p_2 s^2_1 \in I$. By 
$P$-convexity of $I$, $p_1 s^2_2\,,\, p_2 s^2_1 \in I$. From $s_1, s_2 \not\in I$ comes $p_1, p_2 \in I$;
whence, by (\dag), $x_is^2_i \in I$. Since $s_i \not\in I$, we get $x_i \in I$, wherefrom 
$\fr {x_i} {y_i} \in I\,\se A I = \se M I$.  \hfl $\Box$

\vspace{-0.35cm}
\bre Even if $Q$ is a total \und{order} of $\se A I$, $P$ {\it may not} be a total order of $A$.
In fact, 

\vspace{-0.2cm}
\nhp \hfl $P \cap -P = {\se {\iota} I}^{\!\!-1}[Q \cap -Q] = {\se {\iota} I}^{\!\!-1}[0]$, \hfl

\vspace{-0.2cm}
\nhp which, in general is not $\{0\}$. Note that, for $x \in A$,

\vspace{-0.15cm}
\nhp \hfl $x \in {\se {\iota} I}^{\!\!-1}[0]\;\, \Lra\;\, {\se {\iota} I}(x) = 0$ (in $\se A I)
\;\, \Lra\;\, \ex\, z \not\in I\, (zx = 0)$; \hfl 
\vspace{-0.15cm}

\nhp in particular, $x$ is a zero-divisor. Thus, $P$ \und{is} a total order when $A$ is an
integral domain. \hfl $\Box$
\ere
\vspace{-0.3cm}

The following result proves two important properties of total preorders in rings:
\vspace{-0.4cm}

\bth \label{TotPreordsAndFans} $(i)$ Let $T$ be a total preorder of a ring $A$. Then, the real
semigroup $\se G {\!A,T}$ is a fan and $\Sper(A,T)$ is totally ordered by specialization.

\vspace{-0.2cm}
\nhp $(ii)$ Let $\se T 0, \se T 1$ be total preorders of a ring $A$, and let 
$T = \se T 0 \cap\, \se T 1$. Assume\vspace{0.05cm} that the set ${\mathrm{PConv}}\,(A,T)$ of\, 
$T$-convex prime ideals of $A$ is totally ordered under inclusion. Then, the real 
semigroup $\se G {\!A,T}$ is a fan.
\eth

\vspace{-0.4cm}
\nhp {\bf Remark.} In case the ring $A$ is a field, $K$, a total preorder is just a (total)
order of $K$. Thus, Theorem \ref{TotPreordsAndFans} is a ring-theoretic analog of the 
well-known fact that the intersection of at most two total orders of a field is a fan, 
namely the {\it trivial fans}, cf. \cite{La}, Prop. 5.3, p. 39. \hfl $\Box$

\nhp {\bf Proof.} (i) By Corollary \ref{TotOrder=Fan} it suffices to check that $\Sper(A,T)\;
(= \se X {G_{A,T}})$ is totally ordered under inclusion (= specialization); the proof is
identical to that showing that the real spectrum of a ring is a root system: let 
$\al, \bt \in \Sper(A,T)$ and assume that $\al \not\subseteq \bt$ and $\bt \not\subseteq \al$;
let $a \in \al \setminus \bt$ and $b \in \bt \setminus \al$; since the preorder $T$ is total,
either $a - b \in T\, \sub\, \bt$ or $b - a \in T\, \sub\, \al$; hence, $a = b + (a - b) \in \bt$
or $b = a + (b - a) \in \al$, absurd. 

\nhp (ii) We check that assumptions (1) -- (3) of \ref{2chains=Fan} are verified by
$\se G {\!A,T}$\,. 
\vspace{-0.15cm}

Since the saturated prime ideals of $\se G {\!A,T}$ are in a bijective, inclusion-preserving 
correspondence with the $T$-convex prime ideals of $A$ (cf. Fact \ref{hatfacts2}), the 
argument proving Proposition \ref{T-convPrimesTotOrd}\,(2) shows that $\se G {\!A,T}$ verifies
condition [Z] in Theorem \ref{fanrepr}, i.e., assumption \ref{2chains=Fan}\,(1). 
\vspace{-0.15cm}

Assumption \ref{2chains=Fan}\,(2) follows from the proof of (i) and:
\vspace{-0.1cm}

\nhp (*)\;\;\; $\Sper(A,T) = \Sper(A,\se T 0) \cup\, \Sper(A,\se T 1)$. 
\vspace{-0.1cm}

\nhp \und{Proof of (*)}. Clearly, $\Sper(A,T_i)\, \sub\, \Sper(A,T)$ for $i = 0,1$.
Assume there\vspace{0.05cm} is $\al \in \Sper(A)$ such that $T\, \sub\, \al$ but 
$\se T 0\,, \se T 1 \not\subseteq \al$; for $i = 0,1$, let $t_i \in T_i \setminus \al$.
Then, $-\se t 0 \in \al$ and $\se t 0 \not\in \se T 1$ (otherwise, 
$\se t 0 \in \se T 0 \cap \se T 1\, \sub\, \al$). Since $\se T 1$ is a total
preorder, $\se t 0 \in -\se T 1$. Likewise, $-\se t 1 \in \al$ and $\se t 1 \in -\se T 0$.

From $\se t 0 \in \se T 0$ and $-\se t 1 \in \se T 0$ we get $-\se t 0 \se t 1 \in \se T 0$;
from $\se t 1 \in \se T 1$ and $-\se t 0 \in \se T 1$ we get $-\se t 0 \se t 1 \in \se T 1$;
hence, $-\se t 0 \se t 1 \in \se T 0 \cap \se T 1\, \sub\, \al$. From 
$-\se t 0, -\se t 1 \in \al$ comes $\se t 0 \se t 1 = (-\se t 0)(-\se t 1)\in \al$. Hence,
$\se t 0 \se t 1 \in \al \cap -\al = \supp(\al)$. Since this is a prime ideal, 
$t_i \in \supp(\al)\, \sub\, \al$ for $i = 0$ or $i = 1$, contradiction.
\vspace{-0.15cm}

In order to prove assumption (3) of \ref{2chains=Fan} we first show:
\vspace{-0.1cm}

\nhp (**)\;\; Every $T$-convex prime ideal $I$ of $A$ is {\it both} $\se T 0$-convex {\it and} 
$\se T 1$-convex.
\vspace{-0.1cm}

\nhp \und{Proof of (**)}. From \cite{BCR}, Prop. 4.2.8\,(ii), p. 87, we know that $I$ is either 
\vspace{0.08cm}$\se T 0$-convex or $\se T 1$-convex. Assume, towards a contradiction, that $I$ is 
$\se T 0$-convex but not $\se T 1$-convex. Then, there are elements $\se t 0, \se t 1 \in \se T 1$ 
such that $\se t 0 + \se t 1 \in I$, but $\se t 0, \se t 1 \not\in I$. Since $I$ is $T$-convex, 
we have $\se t 0, \se t 1 \not\in \se T 0$ and, since $\se T 0$ is a total preorder, 
$-\se t 0, -\se t 1 \in \se T 0$. As we have $-(\se t 0 + \se t 1) \in I$, $\se T 0$-convexity
yields $-\se t 0, -\se t 1 \in I$, whence $\se t 0, \se t 1 \in I$, contradiction.

Now, \cite{BCR}, Prop. 4.3.8, p. 90, finishes the proof: for $i = 0,1$, there is 
$\al_i \in \Sper(A,T_i)$ so that $\supp(\al_i) = I$. \hfl $\Box$

\vspace{-0.4cm}

\bre \label{ContrexToPrecedingThm(ii)} The following example shows that the requirement in item 
(ii) of Theorem \ref{TotPreordsAndFans} does not hold automatically. Let $A = C(\re)$ be the ring 
of real-valued continuous functions on the reals. For $i = 0,1$, let $T_i = \{f \in A \sth f(i) \geq 0\}$ 
and $M_i = \{f \in A \sth f(i) = 0\}$. The (maximal) ideal $M_i$ is $T_i$-convex; hence, with 
$T = \se T 0 \cap \se T 1$, both $\se M 0$ and $\se M 1$ are $T$-convex; however, 
$\se M 0$ and $\se M 1$ are incomparable under inclusion. \hfl $\Box$
\ere

\ere

\vspace{-0.2cm}

\nhp \parbox[t]{220pt}{M. Dickmann\\
Institut de Math\'ematiques de Jussieu --\\
\hspace*{3.3cm} Paris Rive Gauche\\
Universit\'es Paris 6 et 7\\ 
Paris, France\\
e-mail: dickmann@math.univ-paris-diderot.fr} 
\nhp \hfl \parbox[t]{200pt}{A. Petrovich\\
Departamento de Matem\'atica\\
Facultad de Ciencias Exactas y Naturales\\
Universidad de Buenos Aires\\ 
Buenos Aires, Argentina\\
e-mail: apetrov@dm.uba.ar}

\end{document}